\documentclass[12pt]{amsart}

\def\Cal{\mathcal}

\def\b1{\text{\bf 1}}

\def\BC{{\Bbb C}}

\def\BZ{{\Bbb Z}}
\def\BR{{\Bbb R}}
\def\BH{{\Bbb H}}
\def\CA{{\Cal A}}
\def\CI{{\Cal I}}
\def\CJ{{\Cal J}}
\def\CB{{\Cal B}}
\def\CC{{\Cal C}}

\def\CF{{\Cal F}}

\def\CH{{\Cal H}}

\def\CM{{\Cal M}}
\def\CN{{\Cal N}}
\def\CO{{\Cal O}}

\def\CT{{\Cal T}}
\def\CV{{\Cal V}}

\def\CW{{\Cal W}}

\def\fg{{\frak g}}

\def\Vir{{\Cal V}ir}

\def\iso{ \buildrel\sim\over\longrightarrow }

\def\plus{\buildrel .\over +}

\begin{document}

\title[Lagrangian approach to sheaves of vertex algebras]
{\bf Lagrangian approach to sheaves of vertex algebras}
\thanks{partially supported by the National Science Foundation}

\author{Fyodor Malikov}

\maketitle

\begin{abstract}

We explain how sheaves of vertex algebras are related to
mathematical structures encoded by a class of Lagrangians. The
exposition is focused on two examples: the WZW model and the
(1,1)-supersymmetric $\sigma$-model. We conclude by showing how to
construct a family of vertex algebras with base the
Barannikov-Kontsevich moduli space.

\end{abstract}
\bigskip
{\bf Contents}

\bigskip
Introduction

1 Diffieties and homotopy presymplectic strcutures

\indent\indent 1.1 The jets

\indent\indent 1.2 Local formulas

\indent\indent 1.3 De Rham complex

\indent\indent 1.4 Differential equations

\indent\indent 1.5 A homotopy presymplectic structure

\indent\indent 1.6 Calculus of variations and integrals of motion.
Bosonic $\sigma$-model

2 Vertex Poisson algebras

\indent\indent 2.1 Definition

\indent\indent 2.2 Lemma

\indent\indent 2.3 Tensor products

\indent\indent 2.4 From vertex Poisson algebras to Courant
algebroids

\indent\indent 2.5 Symbols of vertex differential operators

\indent\indent 2.6 A sheaf-theoretic version

\indent\indent 2.7 A natural sheaf of SVDOs

\indent\indent 2.8 Lagrangian interpretation

\indent\indent 2.9 An example: WZW model

3 Supersymmetric analogues

\indent\indent 3.1 Bits of supergeometry

\indent\indent 3.2 Homotopy presymplectic structures

\indent\indent 3.3 Calculus of variations

\indent\indent 3.4 An example: (1,1)-supersymmetric $\sigma$-model

\indent\indent 3.5 Vertex Poisson algebra interpretation. Witten's
models

\indent\indent 3.6 Quantization. B-model moduli.

\newpage

\centerline {\bf Introduction}

\bigskip

More than anything else, the present notes are a report on what we
have been able to make out of the recent papers [Kap,W4]. Even after
the tremendous effort  [QFS], much of mathematical literature
treating various aspects of string theory and related topics is
conspicuously lacking any mention of the Lagrangian, an object that
is at the heart of a physical theory. We would like to make precise
the relation of  sheaves of vertex algebras [MSV] to a Lagrangian
field theory.

What  vertex algebras and Lagrangians have in common is that both
produce infinite dimensional Lie algebras. If $V$ is a vertex
algebra, then, in particular, it is a vector space with
multiplications $_{(n)}$, $n\in\BZ$, and a derivation $T$. The
corresponding Lie algebra is
$$
\text{Lie}(V)=(V/T(V), _{(0)}).\eqno{(0.1)}
$$
We will be concerned with the class of vertex algebras, or rather of
sheaves thereof, introduced in [MSV]. To make our life easier, we
will, first, mostly consider their quasiclassical limits, i.e., the
corresponding vertex Poisson algebras and, second, work in the
$C^{\infty}$-setting. This class comprises  vertex algebra analogues
of sheaves of symbols of differential operators, and their natural
habitat is different versions of $\infty$-jet spaces. To give an
example, let $\Sigma'$ be a 1-dimensional real manifold and consider
$J^{\infty}(T^{*}M_{\Sigma'})$, the space of $\infty$-jets of
sections of the trivial bundle
$$
T^{*}M_{\Sigma}=T^{*}M\times\Sigma'\rightarrow\Sigma'.
$$
We will find it convenient to work with families of such jet-spaces,
to be denoted $J^{\infty}(T^{*}M_{\Sigma/\Sigma''})$, that are are
naturally and similarly attached to a ``time'' fibration
$\tau:\Sigma\rightarrow\Sigma''$ with fiber $\Sigma'$.

The push-forward of the structure sheaf
$\CO_{J^{\infty}(T^{*}M_{\Sigma/\Sigma''})}$ on $M_{\Sigma}$ carries
a structure of a sheaf of vertex Poisson algebras, a vertex analogue
of the Poisson algebra of functions on $T^{*}M$.  This sheaf is
natural in that the assignment $M\rightarrow
\CO_{J^{\infty}(T^{*}M_{\Sigma/\Sigma''})}$ is functorial in $M$.

Such sheaves of {\it symbols of vertex differential operators},
SVDO, can be defined axiomatically and classified [GMS1].  A simple
quiasiclassical version of the classification obtained in [GMS1]
shows that locally (on $M$) all SVDOs are isomorphic and the set of
isomorphisms classes is identified with $H^{3}(M,\BR)$. In
particular, for each closed 3-form $H$, an SVDO
$\CO_{J^{\infty}(T^{*}M_{\Sigma/\Sigma''})}\plus H$ is defined. To
each such SVDO construction (0.1) attaches a sheaf of Lie algebras
on $M$
$$
\text{Lie}(\CO_{J^{\infty}(T^{*}M_{\Sigma/\Sigma''})}\plus
H)=(\CO_{J^{\infty}(T^{*}M_{\Sigma/\Sigma''})}\plus
H)/T(\CO_{J^{\infty}(T^{*}M_{\Sigma/\Sigma''})}\plus H).\eqno{(0.2)}
$$
One of the more interesting examples arises when $M$ is a compact
simple Lie group. In this case the family of isomorphism classes of
SVDOs is 1-dimensional and we denote the representatives by
$SD_{G,k}$, $k\in\BC$. If we let $\fg$ be the corresponding Lie
algebra and $V(\fg)_{k}$ the corresponding vertex Poisson algebra,
then there arise two Poisson-commuting subalgebras [F,FP,AG,GMS2]
$$
V(\fg)_{k}\buildrel j_{l}\over\hookrightarrow \Gamma (G,
SD_{G,k})\buildrel j_{r}\over\hookleftarrow V(\fg)_{-k}.\eqno{(0.3)}
$$
engendered by the left/right translations of $G$ by itself.
 This implies
the existence of 2 copies of the affine Lie algebra (especially when
$\Sigma'$ is a circle)
$$
\hat{\fg}_{k}\buildrel j_{l}\over\hookrightarrow \Gamma (G,
\text{Lie}(SD_{G,k}))\buildrel j_{r}\over\hookleftarrow
\hat{\fg}_{-k}.\eqno{(0.3)}
$$
The case of $k\neq 0$ is exceptional; in this case the
$\hat{\fg}_{k}\times \hat{\fg}_{-k}$-module structure of the SVDO
has a  form familiar from the WZW
$$
\Gamma (G,
\text{Lie}(SD_{G,k}))=\bigoplus_{\lambda}V_{\lambda,k}\otimes
V_{\lambda^{*},-k}.\eqno{(0.4)}
$$

\bigskip

 On the Lagrangian side one deals
with a somewhat different jet-space, $J^{\infty}(M_{\Sigma})$,
where, recall, $\Sigma$ is 2-dimensional. Defined on it there is a
sheaf of variational bi-complexes,
$(\Omega^{\bullet,\bullet}_{J^{\infty}(M_{\Sigma})},\delta,d)$.  An
action, $S$, is defined to be
$$
S\in\Gamma(J^{\infty}(M_{\Sigma}),
\Omega^{2,n}_{J^{\infty}(M_{\Sigma})}/d\Omega^{1,n}_{J^{\infty}(M_{\Sigma})}),\;
n=\text{dim}M.\eqno{(0.5)}
$$
and can be represented by a collection
$$
\{L^{j}\in \Omega^{2,n}_{J^{\infty}(M_{\Sigma})}\}\eqno{(0.6)}
$$
of locally defined Lagrangians, equal to each other modulo $d$-exact
terms on double intersections. Action $S$ also produces a Lie
algebra, $\CI_{S}$, the algebra of integrals of motion arising by
virtue of the N\"other theorem. Let us relate this algebra to (0.2).

Each action (0.5) defines a space, $Sol_{S}$, often referred to as
the solution space, which for our purposes is better be chosen to be
the infinite dimensional submanifold of $J^{\infty}(M_{\Sigma})$
defined by the Euler-Lagrange equations. $Sol_{S}$ and the
jet-spaces considered are examples of a {\it diffiety}, the notion
introduced by A.M.Vinogradov [V].

Attached to each action is what is usually called a variational
2-form on $Sol_{S}$, $\omega_{S}$. Generally speaking, it is not a
form but a global section
$$
\omega_{S}\in
\Gamma(Sol_{S},\Omega^{1,2}_{Sol_{S}}/d(\Omega^{0,2}_{Sol_{S}}))\eqno{(0.7)}
$$
annihilated by $\delta$.

Special  diffiety properties allow to attach to each such form, not
necessarily coming from a Lagrangian, a Lie algebra structure on a
certain subsheaf of ``functions'' on $Sol_{S}$. We call it  a
homotopy presymplectic structure and denote the corresponding sheaf
of Lie algebras by $\CH^{\omega_{S}}_{Sol_{S}}$. One has
$$
\CI_{S}\subset\Gamma(Sol_{S},
\CH^{\omega_{S}}_{Sol_{S}}).\eqno{(0.8)}
$$
We have come to the point. For a class of Lagrangians, comprising
those that are convex and of order 1, $Sol_{S}$ is diffeomorphic to
$J^{\infty}(TM_{\Sigma/\Sigma''})$ and there is a version of the
Legendre transform
$$
g: J^{\infty}(TM_{\Sigma/\Sigma''})\rightarrow
J^{\infty}(T^{*}M_{\Sigma/\Sigma''})
$$
that delivers a Lie algebra sheaf isomorphism
$$
g^{\#}: \CH^{\omega_{S}}_{Sol_{S}}\iso
g^{-1}\text{Lie}(\CO_{J^{\infty}(T^{*}M_{\Sigma/\Sigma''})})\eqno{(0.9)}
$$
in the case of a single globally defined Lagrangian.

Classification of SVDOs is also reflected in the Lagrangian world:
given a globally defined Lagrangian $L$ and a closed 3-form $H$,
known as an $H$-flux,   one defines, following e.g. [GHR,W1], a
collection of Lagrangians  $L^{H}$, such as in (0.6), so that there
is an isomorphism
$$
g^{\#}: \CH^{\omega_{S}}_{Sol_{S}}\iso
g^{-1}\text{Lie}(\CO_{J^{\infty}(T^{*}M_{\Sigma/\Sigma''})}\plus
H),\eqno{(0.10)}
$$
hence an embedding
$$
\CI_{S}\subset \Gamma(M,
\text{Lie}(\CO_{J^{\infty}(T^{*}M_{\Sigma/\Sigma''})}).\eqno{(0.11)}
$$

 As an illustration,  let us relate the WZW model [W1, GW] to $SD_{G,k}$.
 It is conformally invariant meaning that 2 copies, left and right moving, of the
 Virasoro algebra are among the integrals of motion $\CI_{WZW}$.
 Both are embedded into $\text{Lie}(SD_{G,k})$ by virtue of (0.11).
  Precisely when the
 level, $k$, is non-zero, they coincide with the Virasoro algebra
 defined via the Sugawara construction inside the corresponding
 copies of the affine Lie algebra, see (0.3). The existence of the
 left/right moving Virasoro algebra allows to define the left/right
 moving subalgebra of $\text{Lie}(SD_{G,k})$, or indeed, the
 left/right moving subalgebra of $SD_{G,k}$. Again, precisely when
 $k\neq 0$,  the spaces of global sections of these coincide with
 resp. copies of the affine Poisson vertex algebra
 $j_{l}(V(\fg)_{k})$ or $j_{k}(V(\fg)_{-k})$, see (0.3). This
 is an easy consequence of decomposition (0.4).

 A disclaimer is in order. No attempt at originality is made. But it
 is also true that we failed to find an exposition of this material
 suitable for our purposes. Our main source of fact and inspiration was
 [Di], Ch.19; see also [DF,V,Z]. Needless to say,  much of
 what has just been discussed is contained in one form or another in
 [BD], e.g. sect. 2.3.20, 3.9, but note that the meaning is somewhat
 different: we work in the $C^{\infty}$-setting and our constructions
 are not necessarily chiral. In fact, one of our wishes was to understand ``left
 and right movers'', whose ubiquity in physics literature
 bedevils some of us in the mathematics community.

 All of the above has a more or less straightforward superanalogue.
 As an example, we analyze the (1,1)-supersymmetric $\sigma$-model
 arising on a Riemannian manifold $M$. It is similarly governed by
 a super SVDO, $\Omega^{poiss}_{M}$, which is a quasiclassical limit
 of the $C^{\infty}$-version of the chiral de Rham complex [MSV]. The Lie
 superalgebra of integrals of motion contains two copies of the N=1
 superconformal algebra, and we write explicit formulas for their
 embeddings into $\Gamma(M,\text{Lie}(\Omega^{poiss}_{M}))$.

 If, in addition, $M$ is a K\"ahler manifold, then   this
 symmetry algebra is enlarged to include 2 copies of the N=2
 superconformal algebra, a remarkable fact known since [Z,A-GF].
 If so, a quadruple of operators, $Q^{\bullet,\bullet}$,
 $\bullet=\pm$, arises, appropriate combinations of which are
 differentials on $\Omega^{poiss}_{M}$. The 3 cohomology sheaves,
 $$
 H_{Q^{--}+Q^{++}}(\Omega^{poiss}_{M}),\;
 H_{Q^{--}+Q^{+-}}(\Omega^{poiss}_{M}),\;
 H_{Q^{--}}(\Omega^{poiss}_{M}).\eqno{(0.12)}
 $$
 are versions of the quasiclassical limit of Witten's A-, B- and
 half-twisted models [W2]. Their cohomology can be computed
 using versions of the de Rham complex and $\bar{\partial}$ resolution;
 the result is $H^{*}(M,\BC)$,
 $H^{*}(M,\Lambda^{*}\CT_{M})$, and $H^{*}(M,\Omega^{poiss,an}_{M})$ resp.,
 where the latter is a purely holomorphic version of
 $\Omega^{poiss}_{M}$, the quasiclassical limit of the chiral de
 Rham complex [MSV].

 Note that in the present situation the left/right moving subalgebra
 can also be defined by analogy with the WZW-case considered above. Of the
  3 cohomology vertex Poisson algebras, $H^{*}(M,\BC)$,
 $H^{*}(M,\Lambda^{*}\CT_{M})$, and $H^{*}(M,\Omega^{poiss,an}_{M})$, the last one is
 an infinite dimensional
  subquotient of the left moving algebra and is often referred to as
  Witten's chiral algebra [W3,W4]. The first two are also appropriate
  subquotients but are finite dimensional and of topological nature.

\begin{sloppypar}
 The construction of the sheaf complexes used, such as
 $(\Omega^{poiss}_{M}, Q^{--}+Q^{++})$, $(\Omega^{poiss}_{M},Q^{--}+Q^{+-})$,
 $(\Omega^{poiss}_{M}, Q^{--})$, is easily quantized to produce the complexes
 of vertex algebra sheaves, $(\Omega^{vert}_{M}, Q^{--}+Q^{++})$, $(\Omega^{vert}_{M},Q^{--}+Q^{+-})$,
 $(\Omega^{vert}_{M}, Q^{--})$. The sheaf cohomology
 of the first two  remains the same, that of the 3rd is
 quite different and equals the cohomology of the chiral de Rham
 complex. When formulated in the physics language, the relation of this quantization to the genuine
 supersymmetric model is that the latter equals the former
 ``perturbatively'';
  this is the main result of Kapustin [Kap]. It seems, however, that the
 emphasis made in [Kap] on the ``infinite metric limit'' is somewhat
 disingenuous: the above  supersymmetries depend on a genuine
 metric and its K\"ahler form..
 \end{sloppypar}

 By skillfully applying the techniques of SUSY vertex algebras [HK],
 Ben-Zvi, Heluani, and  Szczesny have recently solved [B-ZHS] a harder problem of
 finding quantum versions of the above mentioned embeddings of the
 various
 superconformal algebras in the chiral de Rham complex. Our discussion
 seems to indicate that their quantization is related to
 physics (apart from non-perturbative effects) in about the same way
 the known quantization of (0.3) [FP,F,AG,GMS2] is related to the
 quantum WZW: it mixes the chiral and anti-chiral sectors. At
 this point we can only ask, following [FP], if there is a physical
 model of interest whose chiral algebra is as in [B-ZHS].

The differential $Q^{--}$ of the (quantum) complex
$(\Omega^{vert}_{M}, Q^{--})$ can be deformed. One of the ways to
think about the complex $(\Omega^{vert}_{M}, Q^{--})$  is that it is
a vertex algebra version of the the $\bar{\partial}$-resolution of
the algebra of polyvector fields.
 We conclude by showing that the Barannikov-Kontsevich construction [BK] has a
 vertex algebra analogue: there is a family of vertex algebras with
 base the Barannikov-Kontsevich moduli space. The conformal weight
 zero subspace of this family encodes precisely  the
 Frobenius manifold introduced in [BK]. Our construction amounts to
 defining a morphism of the deformation functor of the Lie algebra of
 polyvector fields to that of
 $\Gamma(M,\text{Lie}(\Omega^{vert}_{M}))$.

 This furnishes the B-model moduli for Witten's half-twisted model.
 The instanton effects seem to be out of reach, to us; see, however,
 an intriguing sentence in [W4] and [FL] for an interesting new approach
 based on marginal operators.

 \bigskip

 {\bf Acknowledgments.} The author thanks V.Gorbounov,
 A.Kapustin, and B.Khesin for illuminating discussions. Parts of this
 work were done while the author was visiting the Fields Institute,
 IHES,  Max-Planck-Institut in Bonn, and Erwin Schr\"odinger Institut in Vienna.
  It is a pleasure to acknowledge
 the support, hospitality, and stimulating atmosphere of these
 institutions.

\bigskip

\bigskip
\bigskip
\centerline{{\bf 1. Diffieties and Homotopy Presymplectic Structures
}}

\bigskip

{\bf 1.1. The jets.}

\bigskip

Assume given a $d$-dimensional $C^{\infty}$-manifold $\Sigma$, the
``world sheet'', and a smooth fiber bundle
$$
\tau :\Sigma\rightarrow\Sigma''\eqno{(1.1.1)}
$$
 with base $\Sigma''$, an open subset of $\BR$,
and fiber $\Sigma'$, a $(d-1)$-dimensional manifold. This is the
minimal requirement; most of the time, we will have the Cartesian
product
$$
\Sigma'\times \Sigma''\eqno{(1.1.2)}
$$
 with a fixed coordinate (``
intrinsic time'')
$$
\tau=\sigma^{0}:\;\Sigma''\rightarrow\BR,\eqno{(1.1.3)}
$$
on $\Sigma''$  and \'etale coordinates
$$
\sigma=(\sigma^{1},...,\sigma^{d-1}):\;\BR^{d-1}\rightarrow\Sigma'\eqno{(1.1.4)}
$$
on $\Sigma'$, so that (1.1.4) is a universal cover -- \'etale
because we would like to include the case of a torus; furthermore,
we will be mostly interested in  $d=1$.

Let $M$,  ``space-time'', be an an $n$-dimensional
$C^{\infty}$-manifold and $M_{\Sigma}=M\times \Sigma$. There arises
a fiber bundle $M_{\Sigma}\rightarrow\Sigma$, and we denote by
$J^{k}(M_{\Sigma})$ the space of $k$-jets of its sections.

Each $J^{k}(M_{\Sigma})$ is a finite dimensional
$C^{\infty}$-manifold, and the natural projections
$$
\pi_{k,l}: J^{k}(M_{\Sigma})\rightarrow J^{l}(M_{\Sigma}),\;k\geq l
$$
organize the collection $\{J^{k}(M_{\Sigma}), k\geq 0\}$ in a
projective system. The space of $\infty$-jets,
$J^{\infty}(M_{\Sigma})$, is the projective limit,
$\lim_{\leftarrow}J^{k}(M_{\Sigma})$. Its sheaf of smooth functions
is the direct limit of those on $J^{k}(M_{\Sigma})$, $\geq 0$. Let
$\CO_{J^{k}(M_{\Sigma})}$ be the sheaf of all smooth functions on
$J^{k}(M_{\Sigma})$ that are polynomials in positive order jets and
define
$$
\CO_{J^{\infty}(M_{\Sigma})}=\lim_{\rightarrow}\CO_{J^{k}(M_{\Sigma})}.
$$
The fiber bundle $J^{\infty}(M_{\Sigma})\rightarrow \Sigma$ carries
the well-known flat connection, or equivalently,
$\CO_{J^{\infty}(M_{\Sigma})}$ is a sheaf of $D_{\Sigma}$-algebras,
where $D_{\Sigma}$ is the sheaf of differential operators on
$\Sigma$. Denote by
$$
\rho: \CT_{\Sigma}\rightarrow\CT_{J^{\infty}(M_{\Sigma})}
\eqno{(1.1.5)}
$$
the corresponding morphism of the sheaves of vector fields. We will
often refer to this situation by calling
$J^{\infty}(M_{\Sigma})\rightarrow \Sigma$ a $D_{\Sigma}$-manifold,
thus mimicking [BD].

In particular, attached to any tangent vector $\xi\in T_{t}\Sigma$,
there is a tangent vector $\rho(\xi)\in
T_{(x,t)}(J^{\infty}(M_{\Sigma}))$. Thus there arises an integrable
 $d$-dimensional distribution
$$
M_{\Sigma}\ni (x,t)\mapsto\text{ span}\{\hat{\xi},\;\xi\in
T_{t}\Sigma\}\subset T_{(x,t)}J^{\infty}(M_{\Sigma})
$$
known as the {\it Cartan distribution}.

$J^{\infty}(M_{\Sigma})$ is a simple example of what A.M.Vinogradov
calls a {\it diffiety}. Since the Cartan distribution is an
important structure ingredient, by an infinitesimal automorphism of
$J^{\infty}(M_{\Sigma})$ one means a vector field that preserves the
Cartan distribution. This defines a Lie algebra sheaf, which
contains the Lie algebra subsheaf of vertical (i.e. tangent to the
fibers of the projection $J^{\infty}(M_{\Sigma})\rightarrow\Sigma$)
vector fields. We will let $\text{Evol}(J^{\infty}(M_{\Sigma}))$
denote this Lie algebra sheaf and call any, locally defined, vector
field $\xi\in \text{Evol}(J^{\infty}(M_{\Sigma}))$ {\it
evolutionary}.

All of this admits a relative version: if one defines
$J^{\infty}(M_{\Sigma/\Sigma''})$ to be the space of jets of
sections $\Sigma\rightarrow M_{\Sigma}$ in  the direction of fibers
of the bundle $\Sigma\rightarrow\Sigma''$, then the definitions of
the connection
$$
\rho:
\CT_{\Sigma/\Sigma''}\rightarrow\CT_{J^{\infty}(M_{\Sigma/\Sigma''})},
\eqno{(1.1.6)}
$$
(where $\CT_{\Sigma/\Sigma''}$ is the sheaf of vector fields on
$\Sigma$ tangent to the fibers of the projection
$\Sigma\rightarrow\Sigma''$), Cartan distribution, sheaf of
evolutionary vector fields
$\text{Evol}(J^{\infty}(M_{\Sigma/\Sigma''}))$, etc., are immediate.

{\em From now, unless otherwise stated, we will be working over base
$S$, $S$ being either $\Sigma''$ or a point.}

If (1.1.2-4) are valid, then the fibers of the projection
$J^{\infty}(M_{\Sigma})\rightarrow\Sigma$ are canonically
identified, and
$$
J^{\infty}(M_{\Sigma})\iso J^{\infty}_{\Sigma}(M)\times\Sigma,\;
J^{\infty}(M_{\Sigma/\Sigma''})\iso
J^{\infty}_{\Sigma/\Sigma''}(M)\times\Sigma
 \eqno{(1.1.7)}
$$
for some infinite dimensional manifolds, $J^{\infty}_{\Sigma}(M)$
and  $J^{\infty}_{\Sigma/\Sigma''}(M)$, whose definition is easy to
reconstruct from 1.2 below.

\bigskip

{\bf 1.2. Local Formulas}

\bigskip
 As an illustration, and for future use,
let us show what all of this means in terms of local coordinates.
Let $x^{1},...,x^{n}$ be local coordinates on  $M$, those on
$\Sigma$ being defined by (1.1.2-4). For a multi-index
$(m)=(m_{0},...,m_{d-1})$, let
$$
x^{j}_{(m)}=\partial_{\sigma^{0}}^{m_{0}}\cdots\partial_{\sigma^{d-1}}^{m_{d-1}}x^{j},
$$
where $\partial_{\sigma^{m}}=\partial/\partial\sigma^{m}$ and
$x^{j}$ is regarded, formally, as a function of
$\sigma^{0},...,\sigma^{d-1}$. Then
$$
\{\sigma^{i},x^{j}_{(m)}: \;0\leq i\leq d-1, 1\leq j\leq n,
(m)\in\BZ_{+}^{d}\}
$$
are local coordinates on $J^{\infty}(M_{\Sigma})$
$$
\{\sigma^{i},x^{j}_{(m)}: \;0\leq i\leq d-1, 1\leq j\leq n,
(m)\in\BZ_{+}^{d},m_{0}=0\}
$$
are local coordinates on $J^{\infty}(M_{\Sigma/\Sigma''})$,
 and sections of
$\CO_{J^{\infty}(M_{\Sigma/S})}$ are smooth functions in
$\sigma^{i},x^{j}$ and polynomials in $x^{j}_{(m)}$, $(m)\neq 0$,
with $m_{0}=0$ if $S=\Sigma''$.

Let $\delta/\delta x^{j}_{(m)}$ denote the vertical vector field
$\partial/\partial x^{j}_{(m)}\in \CT_{J^{\infty}(M_{\Sigma})}$.
Morphism  (1.1.5) is defined by
$$
\rho(\partial_{\sigma^{i}})=\partial_{\sigma^{i}}+\sum_{j=1}^{n}\sum_{(m)\in\BZ_{+}^{d}}x^{j}_{(m+e_{i})}\frac{\delta}{\delta
x^{j}_{(m)}},\eqno{(1.2.1)}
$$
where $e_{i}=(0,...,0,1,0,...,0)$, 1 appearing at the $i$-th
position.

Evolution vector fields are in 1-1 correspondence with $n$-tuples of
functions $F^{1},...,F^{n}\in \CO_{J^{\infty}(M_{\Sigma})}$ (called
the characteristic of a vector field) and are written
$$
\xi=\sum_{j=1}^{n}\left( F^{j}\frac{\delta}{\delta
x^{j}}+\sum_{(m)\neq
0}\left(\rho(\partial_{\sigma^{0}})^{m_{1}}\cdots\rho(\partial_{\sigma^{m_{d-1}}})^{m_{d-1}}F^{j}\right)\frac{\delta}{\delta
x^{j}_{(m)}}\right)\eqno{(1.2.2)}
$$
The relative analogs of (1.2.1,2) are obvious.

\bigskip

{\bf 1.3. De Rham Complex.}

\bigskip
Reflecting the product structure of $M_{\Sigma}$, the de Rham
complex on $J^{\infty}(M_{\Sigma/S})$ is bi-graded:
$$
\Omega_{J^{\infty}(M_{\Sigma/S})}^{p}=\oplus_{i+j=p}\Omega_{J^{\infty}(M_{\Sigma/S})}^{i,j}.
$$
It carries 2 anti-commuting differentials:
$$
\delta:\Omega_{J^{\infty}(M_{\Sigma})}^{i,j}\rightarrow
\Omega_{J^{\infty}(M_{\Sigma})}^{i,j+1},\;
d_{\rho/S}:\Omega_{J^{\infty}(M_{\Sigma/S})}^{i,j}\rightarrow
\Omega_{J^{\infty}(M_{\Sigma/S})}^{i+1,j},\; [\delta,d_{\rho/S}]=0,
$$
defined as follows. The space, or rather the sheaf of spaces,
$\Omega_{J^{\infty}(M_{\Sigma/S})}^{*,*}$ is naturally an
$\CO_{\Sigma}$-module, and $\delta$ is the vertical de Rham
differential, i.e., the one that is $\CO_{\Sigma}$-linear. The flat
connection $\rho$ gives rise to a differential, $d_{\rho/S}$, on
$\Omega_{J^{\infty}(M_{\Sigma/S})}^{*,0}$ in the standard manner.
For example, in terms of local cordinates,
$$
d_{\rho/S}x^{j}_{(m)}=\sum_{i=\epsilon}^{d-1}x^{j}_{(m)+e_{i}}d\sigma^{i},
\eqno{(1.3.1)}
$$
where $\epsilon=0$ if $S=\emptyset$ and $\epsilon=1$ if
$S=\Sigma''$.

 Then the condition $[\delta, d_{\rho/S}]$ allows one to extend
$d_{\rho/S}$ to the entire $\Omega_{J^{\infty}(M_{\Sigma/S})}^{*,*}$
unambiguously. Thus
$$
\delta F(\sigma,x_{(m)})=\frac{\partial F(\sigma,x_{(m)})}{\partial
x_{(m)}}\delta x_{(m)};\; \delta d\sigma=0,\eqno{(1.3.2)}
$$

$$
d_{\rho/S}\delta x^{j}_{(m)}=-\delta d_{\rho/S}
x^{j}_{(m)}=-\sum_{i=\epsilon}^{d-1}\delta x^{j}_{(m)+e_{i}}\wedge
d\sigma^{i}, \eqno{(1.3.3)}
$$
There is a mapping of bi-complexes
$$
\Omega_{J^{\infty}(M_{\Sigma})}^{*,*}\rightarrow
\Omega_{J^{\infty}(M_{\Sigma/\Sigma''})}^{*,*},\eqno{(1.3.4)}
$$
that sends a form to its restriction to the fibers of the composite
projection
$$
J^{\infty}(M_{\Sigma})\rightarrow\Sigma\buildrel\tau\over\rightarrow\Sigma''.
\eqno{(1.3.5)}
$$
As a practical matter, (1.3.4) amounts to
$$
d\tau\mapsto 0.\eqno{(1.3.6)}
$$

Let $\iota_{\xi}$ be the operator of contraction with a vector field
$\xi$. A straightforward computation proves the following.

\bigskip

{\bf 1.3.1. Lemma.} {\it A vertical vector field $\xi$ on
$J^{\infty}(M_{\Sigma/S})$ is evolutionary iff}
$$
[\iota_{\xi},d_{\rho/S}]=0. \eqno{(1.3.7)}
$$
\bigskip

{\bf 1.3.2. Corollary.} {\it If $\xi$ is evolutionary, then}
$$
[d_{\rho/S},\text{Lie}_{\xi}]=0. \eqno{(1.3.8)}
$$
\bigskip

Indeed,
$$
[d_{\rho/S},\text{Lie}_{\xi}]=[d_{\rho/S},[\delta,\iota_{\xi}]]=-[\delta,[d_{\rho/S},\iota_{\xi}]]=0.
$$

\bigskip\bigskip

{\bf 1.4. Differential Equations.}

\bigskip
Let $\CJ\subset \CO_{J^{\infty}(M_{\Sigma})}$ be a sheaf of ideals
closed under the connection $\rho$.  Let
$$
Sol\subset J^{\infty}(M_{\Sigma})
$$
be the zero locus of this ideal. If some regularity conditions hold,
then his submanifold delivers another example of a diffiety. For
example, one can, and we will, assume that $\CJ$ is  locally
pseudo-Cauchy-Kovalevskaya, i.e., there is a distinguished
coordinate on $\Sigma$, say, $\tau$, and for any point in $M$ there
is a coordinate system $x^{1},...,x^{n}$ s.t. the ideal is
generated, around the pre-image of this point on the jet space, by
the functions $E^{1},...,E^{n}$ satisfying
$$
E^{j}=\rho(\partial_{\tau})^{l_{j}}x^{j}+\cdots, \eqno{(1.4.1)}
$$
where $\cdots$ stand for the terms that do not involve jets of $x$'s
of degree $\geq l_{j}$ in the direction of $\partial_{\tau}$.

Here are some of the structure properties $Sol$ shares with the
ambient jet space:

$Sol$ is fibered over $\Sigma$, hence over $\Sigma''$ via
$\Sigma\buildrel\tau\over\rightarrow\Sigma''$;

the algebra of functions $\CO_{Sol}$ is a $D_{\Sigma}$-algebra
(because the flat connection preserves $Sol$), hence a
$D_{\Sigma/\Sigma''}$-algebra, where $D_{\Sigma/\Sigma''}$ is the
subalgebra of $D_{\Sigma}$ that commutes with
$\tau^{-1}\CO_{\Sigma''}$; we will write $Sol_{\Sigma/\Sigma''}$ if
we wish to emphasize the $D_{\Sigma/\Sigma''}$-algebra structure;

the de Rham complex $\Omega_{Sol/S}^{*}$ is bi-graded and carries
two commuting differentials, $\delta$, the vertical differential,
and $d_{\rho/S}$, the $D_{\Sigma/S}$-module differential.

If (1.4.1) is valid, then solving the equation $E^{j}=0$ for
$\rho(\partial_{\tau})^{l_{j}}x^{j}$ one sees that $Sol$ looks like
$\infty$-jets in the direction of the fiber of the bundle
$\Sigma\buildrel\tau\over\rightarrow\Sigma''$ to something finite
dimensional. In particular, if $l_{1}=l_{2}=\cdots =l_{n}=2$, then
$$
Sol\iso J^{\infty}(TM_{\Sigma/\Sigma''}), \eqno{(1.4.2)}
$$
as $D_{\Sigma/\Sigma''}$-manifolds.

Any evolutionary vector field on the ambient jet space that
preserves the ideal $\CJ$ descends to a vector field on $Sol$, which
still satisfies (1.3.7). We emulate this situation by making the
following definition.

\bigskip

{\bf 1.4.1. Definition.} {\it A vertical vector field $\xi$ on $Sol$
is called evolutionary (relative to $S$) if}
$$
[\iota_{\xi},d_{\rho/S}]=0. \eqno{(1.4.3)}
$$
\bigskip

{\bf 1.4.2. Lemma.} {\it If $\xi$ is evolutionary, then}
$$
[d_{\rho/S},\text{Lie}_{\xi}]=0. \eqno{(1.4.4)}
$$
\bigskip

{\bf 1.4.3.} Let $\text{Evol}(Sol)_{\Sigma/S}$ denote the sheaf of
all evolutionary vector fields on $Sol$. The relation
$[\text{Lie}_{\xi},\iota_{\eta}]=\iota_{[\xi,\eta]}$ combined with
(1.4.3,4) implies that $\text{Evol}(Sol)_{\Sigma/S}$ is a Lie
algebra.

Note that
$$
\text{Evol}(Sol)_{\Sigma}\subset \text{Evol}(Sol)_{\Sigma/\Sigma''}.
\eqno{(1.4.5)}
$$

\bigskip

\bigskip\bigskip

{\bf 1.5. A  homotopy pre-symplectic structure.}

\bigskip

A symplectic form, that is, a non-degenerate closed 2-form gives
rise to a Poisson algebra structure on the structure sheaf of a
manifold. A degenerate closed 2-form similarly gives rise to a
Poisson algebra structure on a certain, {\it admissible}, subalgebra
of the structure sheaf. This subalgebra consists of functions
constant along the leaves of the foliation tangent to the kernel of
the form. We would like to explain, in the spirit of [DF,Di], that
in the case of a diffiety, such as $Sol$, a pre-symplectic structure
gives rise to a Poisson structure, which may  be just as good for
all practical purposes as the symplectic one.

{\bf 1.5.1.} The following is a list of standard symplectic geometry
notions adjusted to the case where  vector fields are replaced with
evolutionary vector fields and   equalities are  valid up to
$d_{\rho/S}$-exact terms.

{\it From now on $\CM$ is a diffiety, such as $Sol\subset
J^{\infty}(M_{\Sigma})$ or $J^{\infty}(M_{\Sigma/S})$.}

 We will call $\omega\in H^{0}(\CM,\Omega_{\CM/S}^{d-1,2})$ a {\it
pre-symplectic  homotopy form} if
$$
\delta\omega\in H^{0}\left(\CM,
d_{\rho/S}\left(\Omega_{\CM/S}^{d-2,1}\right)\right), \eqno{(1.5.1)}
$$
where $d_{\rho/S}\left(\Omega_{\CM/S}^{d-2,1}\right)$ is understood
as the sheaf associated with the presheaf
$\text{Im}d_{\rho/S}\subset \Omega_{\CM/S}^{d-2,1}$.

 An evolutionary vector field
$\xi\in\text{Evol}(\CM)_{\Sigma/S}$ is called {\it
 homotopy Hamiltonian} if
$$
\text{Lie}_{\xi}\omega\in
d_{\rho/S}\left(\Omega_{\CM/S}^{d-2,1}\right).\eqno{(1.5.2)}
$$
Call  $F\in\Omega^{d-1,0}_{\CM/S}$ {\it admissible} if there is an
evolutionary vector field $\xi$ such that
$$
\delta
F\in\iota_{\xi}\omega+d_{\rho/S}\left(\Omega_{\CM/S}^{d-2,1}\right).
\eqno{(1.5.3)}
$$
Let $\Omega^{adm}_{\CM/S}$ be the space of all admissible
$(d-1,0)$-forms.

 Note that (1.5.3) implies that $\xi$ is homotopy Hamiltonian and
 that $d_{\rho/S}(\Omega^{d-2,0}_{\CM/S})\subset\Omega^{adm}_{\CM/S}$.

Introduce
$$
\CF=\{(F,\xi)\text{ such that (1.5.2) holds}\}, \eqno{(1.5.4)}
$$
Let
$$
J=\{(F,\xi)\in\CF\text{ s.t. } F\in
d_{\rho/S}(\Omega^{d-2,0}_{\CM/S})\}. \eqno{(1.5.5)}
$$
Note that
$$
\CF/J=\Omega^{adm}_{\CM/S}/d_{\rho/S}(\Omega^{d-2,0}_{\CM/S}).\eqno{(1.5.6)}
$$
\bigskip

{\bf 1.5.2. Proposition.}

(i) {\it The formula
$$
\{(F,\xi),(G,\eta)\}=(\xi G,[\xi,\eta]) \eqno{(1.5.7)}
$$
defines an algebra structure on $\CF$.

(ii) $J\subset\CF$ is a 2-sided ideal w.r.t. to algebra structure
(1.5.5).

(iii) Thus $\{.,.\}$ descends to the quotient
$\Omega^{adm}/d_{\rho/S}(\Omega^{d-2,0}_{\CM/S})$ and makes the
latter into a Lie algebra.}

\bigskip

{\bf Proof.}

(i) We have
$$
\xi G=\iota_{\xi}\delta G
=\iota_{\xi}\left(\iota_{\eta}\omega+d_{\rho/S}(...)\right)=
\iota_{\xi}\iota_{\eta}\omega+d_{\rho/S}(...),\eqno{(1.5.8)}
$$
where the 2nd equality uses the admissibility of $G$ and the 3rd
uses (1.4.3). This shows, in particular, that $\{.,.\}$ is
anti-commutative modulo $d_{\rho/S}$.

To verify that $(\xi G,[\xi,\eta])\in\CF$ compute further
$$
\delta\xi
G=\delta\left(\iota_{\xi}\iota_{\eta}\omega+d_{\rho/S}(...)\right) =
\text{Lie}_{\xi}\iota_{\eta}\omega-\iota_{\xi}\text{Lie}_{\eta}\omega+
\iota_{\xi}\iota_{\eta}\delta\omega+d_{\rho/S}(...)=
$$
$$
\iota_{[\xi,\eta]}\omega+\iota_{\eta}\text{Lie}_{\xi}\omega-\iota_{\xi}\text{Lie}_{\eta}\omega+
\iota_{\xi}\iota_{\eta}\delta\omega+d_{\rho/S}(...)=\iota_{[\xi,\eta]}\omega+d_{\rho/S}(...),
\eqno{(1.5.9)}
$$
where the 2nd equality follows from the identity
$[\delta,\iota_{\xi}]=\text{Lie}_{\xi}$, the 3rd from the identity
$[\text{Lie}_{\xi},\iota_{\eta}]=\iota_{[\xi,\eta]}$;  and the 4th
follows, because $\text{Lie}_{\xi}\omega$,
$\text{Lie}_{\eta}\omega$, $\delta\omega$ being $d_{\rho/S}$-exact,
(1.5.1-2), the rest of the terms are also  by virtue of the
omnipresent (1.4.3). This proves (i).

(ii) That $\{\CF,J\}\subset J$ follows from definition (1.5.7) and
 (1.4.4). That $\{J,\CF\}\subset J$ follows now from
the anti-commutativity of $\{.,.\}$ modulo $d_{\rho/S}$, see (1.5.8)
and the sentence that follows it.

(iii) The anti-commutativity has already been discussed. The Jacobi
identity is proved as follows:
$$
\{\{F,G\},H\}=[\xi,\eta]H=\xi\eta H-\eta\xi
H=\{F,\{G,H\}\}-\{G,\{F,H\}\}.
$$
\bigskip

{\bf 1.5.3. Definition.} {\it Let $\CH^{\omega}_{\CM/S}$ denote the
sheaf of Lie algebras
$\Omega^{adm}/d_{\rho/S}(\Omega^{d-2,0}_{\CM/S})$ constructed in
Proposition 1.5.2 (iii).}

\bigskip
The absolute and relative versions of this construction can be
compared. Indeed,  there is an embedding
$$
\Omega^{adm}_{\CM}\hookrightarrow\Omega^{adm}_{\CM/\Sigma''},\eqno{(1.5.10)}
$$
which arises by virtue of (1.4.5). Therefore, morphism of
bi-complexes (1.3.4) induces a morphism of the Lie algebra sheaves
$$
\CH^{\omega}_{\CM}\rightarrow
\CH^{\omega}_{\CM/\Sigma''},\eqno{(1.5.11)}
$$
which in terms of local coordinates amounts to
$$
d\tau\mapsto 0,\eqno{(1.5.12)}
$$
cf. (1.3.6).

{\bf 1.5.3.1. Remark.} Note that replacing $\omega$ with $\omega
+d_{\rho/S}\alpha$ gives rise to the sheaf of Lie algebras
$\CH^{\omega+d_{\rho/S}\alpha}_{\CM/S}$ canonically isomorphic to
$\CH^{\omega}_{\CM/S}$. Hence $\CH^{\omega}_{\CM/S}$ is really
associated to
$$
\aligned
 \omega&\in\text{
Ker}\left(\delta:H^{0}\left(\CM,\Omega_{\CM/S}^{d-1,2}/d_{\rho/S}\Omega_{\CM/S}^{d-2,2}\right)
\rightarrow
H^{0}\left(\CM,\Omega_{\CM/S}^{d-1,2}/d_{\rho/S}\Omega_{\CM/S}^{d-2,2}\right)\right)\\
&\subset
H^{0}\left(\CM,\Omega_{\CM/S}^{d-1,2}/d_{\rho/S}\Omega_{\CM/S}^{d-2,2}\right).
\endaligned
$$

\bigskip

\bigskip

{\bf 1.5.4. Example 1: canonical commutation relations.} Replace $M$
as the target space with $T^{*}M$ and consider
$J^{\infty}(T^{*}M_{\Sigma/\Sigma''})$ with $d=dim\Sigma=2$. There
arises the projection
$$
\pi:\; J^{\infty}(T^{*}M_{\Sigma/\Sigma''})\rightarrow
T^{*}M_{\Sigma/\Sigma''}. \eqno{(1.5.13)}
$$
Let $\omega^{o}$ be the canonical symplectic form on $T^{*}M$ and
$\omega=\omega^{o}\wedge d\sigma$; the latter is a $(1,2)$-form on
$T^{*}M_{\Sigma/\Sigma''}$ -- we are taking advantage of coordinates
(1.1.2-4). There arises then $\pi^{*}\omega$, the pull-back of
$\omega$ onto $J^{\infty}(T^{*}M_{\Sigma/\Sigma''})$ under (1.5.13).
Let us now compute
$\CH^{\pi^{*}\omega}_{J^{\infty}(T^{*}M_{\Sigma/\Sigma''})}$.

In this case
$$
\Omega^{adm}_{J^{\infty}(T^{*}M_{\Sigma/\Sigma''})}=\Omega^{1,0}_{J^{\infty}(T^{*}M_{\Sigma/\Sigma''})}.
\eqno{(1.5.14)}
$$
To see why, note that, $\omega^{o}$ being non-degenerate, any
section of $\Omega^{1,1}_{T^{*}M_{\Sigma/\Sigma''}}$ can be written
as $\iota_{\xi^{o}}\omega$ for some vector field $\xi^{o}$ on
$T^{*}M_{\Sigma}$. Pulling back on
$J^{\infty}(T^{*}M_{\Sigma/\Sigma''})$ one sees that likewise any
section of $\pi^{*}\Omega^{1,1}_{T^{*}M_{\Sigma/\Sigma''}}$ can be
written as $\iota_{\xi^{o}}\pi^{*}\omega$, where $\xi^{o}$ is now a
vector field on $J^{\infty}(T^{*}M_{\Sigma/\Sigma''})$, which is
locally a linear combination of $\delta/\delta x^{i}$, $x^{i}$ being
local coordinates on $T^{*}M$. Thinking of $\xi^{o}$ as a
characteristic, one can  prolong it to an evolutionary vector field
$\xi$, as in (1.2.2), and thus obtain
$$
\pi^{*}\Omega^{1,1}_{T^{*}M_{\Sigma/\Sigma''}}=\{\iota_{\xi}\pi^{*}\omega,\;
\xi\in\text{Evol}(J^{\infty}(T^{*}M_{\Sigma/\Sigma''}))\}\eqno{(1.5.15)}.
$$

Now observe that
$\Omega^{1,1}_{J^{\infty}(T^{*}M_{\Sigma/\Sigma''})}$ is generated,
as a $D_{\Sigma/\Sigma''}$-module, by
$\pi^{*}\Omega^{1,1}_{T^{*}M_{\Sigma/\Sigma''}}$. Hence (1.5.15)
holds true for the entire
$\Omega^{1,1}_{J^{\infty}(T^{*}M_{\Sigma/\Sigma''})}$ modulo
$d_{\rho/\Sigma''}$-exact terms, i.e.,
$$
\Omega^{1,1}_{J^{\infty}(T^{*}M_{\Sigma/\Sigma''})}=\{\iota_{\xi}\pi^{*}\omega+d_{\rho/\Sigma''}\beta,\;
\xi\in\text{Evol}(J^{\infty}(T^{*}M_{\Sigma/\Sigma''})),\beta\in\Omega^{0,1}_{J^{\infty}(T^{*}M_{\Sigma/\Sigma''})}
\}\eqno{(1.5.16)}
$$
and (1.5.14) follows.

\bigskip

Because of its importance, the  sheaf of Lie algebras arising in
this way will be denoted thus
$$
\CH^{can}\buildrel\text{def}\over =
\CH^{\pi^{*}\omega}_{J^{\infty}(T^{*}M_{\Sigma/\Sigma''})}=
\Omega^{1,0}_{J^{\infty}(T^{*}M_{\Sigma/\Sigma''})}/d_{\rho/\Sigma''}\left(\Omega^{0,0}_{J^{\infty}(T^{*}M_{\Sigma/\Sigma''})/}\right).
\eqno{(1.5.17)}
$$

Computationally, the gist of our discussion is as follows. The
algebra of functions on the cotangent bundle with the canonical
Poisson bracket is a Lie subalgebra of $\CH^{can}$:
$$
\pi^{\#}:\;\left(\pi^{-1}\CO_{T^{*}M_{\Sigma/\Sigma''}}d\sigma,\{.,.\}_{T^{*}M}\right)\hookrightarrow
\CH^{can}, \eqno{(1.5.18)}
$$
and the rest of the Lie algebra structure is determined by
$$
\{\overline{F},\overline{GHd\sigma}\}=\overline{\{F,Gd\sigma\}H}+\overline{G\{F,Hd\sigma\}},
\eqno{(1.5.19a)}
$$
$$
\{\overline{F},\overline{G\rho(\partial_{\sigma})Hd\sigma}\}=
\overline{\{F,Gd\sigma\}\rho(\partial_{\sigma})H}+\overline{G\rho(\partial_{\sigma})\{F,Hd\sigma\}},
\eqno{(1.5.19b)}
$$
because an evolutionary vector field is a derivation commuting with
$\rho(\partial_{\sigma})$, (1.4.4.).

To see what all of this means, let us compute some brackets.Let
$F,G$ be   functions on $M_{\Sigma/\Sigma''}$, $\xi,\eta$  vector
fields on $M_{\Sigma}$ vertical w.r.t. $M\rightarrow\Sigma$, which
we regard as fiberwise linear functions on $T^{*}M_{\Sigma}$. Then,
$$
\{\overline{Fd\sigma},\overline{Gd\sigma}\}=0, \eqno{(1.5.20a)}
$$
$$
\{\overline{\xi d\sigma},\overline{G d\sigma}\}=\overline{\xi G
d\sigma}, \eqno{(1.5.20b)}
$$
$$
\{\overline{\xi d\sigma},\overline{\eta
d\sigma}\}=\overline{[\xi,\eta] d\sigma}. \eqno{(1.5.20c)}
$$
The first instance of the bracket jet nature  manifesting itself is
as follows. If $F_{i}dx^{i}$ is a 1-form on $M$,then
$\alpha=F_{i}(x)x^{i}_{(1)}$ is a well-defined $(0,0)$-form on
$J^{\infty}(T^{*}M_{\Sigma})$. Having thus embedded $\Omega^{1}_{M}$
into $\Omega^{0,0}_{J^{\infty}(T^{*}M_{\Sigma/\Sigma''})}$, one uses
(1.5.19b) to obtain
$$
\{\overline{\xi d\sigma},\overline{\alpha
d\sigma}\}=\overline{(\text{Lie}_{\xi}\alpha )d\sigma},
\eqno{(1.5.20d)}
$$
if $\xi$ does not depend on $\sigma$ explicitly, and
$$
\{\overline{\xi d\sigma},\overline{\alpha
d\sigma}\}=\overline{(\text{Lie}_{\xi}\alpha
+\iota_{\partial_{\sigma}\xi}\alpha)d\sigma} \eqno{(1.5.20e)}
$$
in general, where $\text{Lie}_{\xi}\alpha$ is the Lie derivative of
$\alpha$ along $\xi$.

Formulas (1.5.20a--d), without functions explicitly depending on
$\tau,\sigma$, are a familiar definition of the Lie algebra
associated with the {\it Courant algebroid} on $TM\oplus T^{*}M$.
The idea that the Courant algebroid has infinite dimensional nature
apparently goes back to I.Dorfman [Dor]. It was revived recently, in
a slightly different context, by P.Bressler [Bre].

Note that identities (1.5.19a,b)  seem to incorporate the Leibnitz
identity, which they do not, because $\CH^{can}$ is not an
associative algebra. It is, however, a quotient of
$\Omega^{1,0}_{J^{\infty}(T^{*}M_{\Sigma/\Sigma''})}$, and the
latter is. In fact,
$\Omega^{1,0}_{J^{\infty}(T^{*}M_{\Sigma/\Sigma''})}$ is a sheaf of
vertex Poisson algebras, and its quotient $\CH^{can}$ is the
canonically associated to it sheaf of
 Lie algebras, see Proposition 2.7.3 below.

\bigskip

{\bf 1.5.5. Example 2: the solution space of an order 2 system in
the pseudo-Cauchy-Kovalevskaya form.} Let us place ourselves in the
situation of 1.4 and let $Sol$ satisfy (1.4.2). Then
$$
Sol_{\Sigma/\Sigma''}\iso J^{\infty}(TM_{\Sigma/\Sigma''}).
$$
 The latter does not carry any canonical 2-form, but let
us fix a diffeomorphism
$$
g: TM\iso T^{*}M, \eqno{(1.5.21)}
$$
which in practice is most often defined by a metric on $M$. There
arises then a diffeomorphism (for any base $S$)
$$
g: J^{\infty}(TM_{\Sigma/S})\iso J^{\infty}(T^{*}M_{\Sigma/S}),
$$
hence a diffeomorphism
$$
g: Sol\iso J^{\infty}(T^{*}M_{\Sigma/\Sigma''}),
$$
and a sheaf  isomorphism
$$
 g^{\#}:\Omega^{.,.}_{Sol/\Sigma''}\iso
 g^{-1}\Omega^{.,.}_{J^{\infty}(T^{*}M_{\Sigma})/\Sigma''},
\eqno{(1.5.22)}
$$
where $g^{-1}$ stands for the inverse image in the category of
sheaves.

 Let $g^{*}\omega$ be the symplectic form on $TM$
obtained by pulling back the canonical symplectic form $\omega$ on
$T^{*}M$. We have arrived at

\bigskip
{\bf 1.5.5.1. Lemma.} {\it Mapping (1.5.20) induces an isomorphism
of Lie algebra sheaves
$$
g^{\#}: \CH^{g^{*}\omega}_{Sol/\Sigma''}\iso g^{-1}\CH^{can}.
$$}

\bigskip

\bigskip\bigskip

{\bf 1.6. Calculus of variations and integrals of motion.}

\bigskip

Calculus of variations is the principal source of the brackets
discussed in 1.4-5.

{\bf 1.6.1.} An action $S$ is a global section
$$
S\in
\Omega^{d,0}_{J^{\infty}(M_{\Sigma})}/d_{\rho}\Omega^{d-1,0}_{J^{\infty}(M_{\Sigma})}
\eqno{(1.6.0)}
$$
 It can be represented by a  Lagrangian which is a
collection of sections
$$
L=\{L^{(i)}\in\Gamma(U_{i},\Omega^{d,0}_{J^{\infty}(M_{\Sigma})})\text{
s.t. } L^{(j)}-L^{(i)}\in\text{Im}d_{\rho}\text{ on }U_{j}\cap
U_{i}\}, \eqno{(1.6.1a)}
$$
determined up to a transformation
$$
L^{(i)}\mapsto L^{(i)}+d_{\rho}\beta^{(i)},\eqno{(1.6.1b)}
$$
where $\{U_{i}\}$ is an open covering of $J^{\infty}(M_{\Sigma})$.

 Choosing local coordinates one observes that
$$
\delta L^{(i)}=-d_{\rho}\gamma^{(i)}+ E_{j}^{(i)}\delta
x^{j},\eqno{(1.6.2)}
$$
for some $\gamma^{(i)}\in \Omega^{d-1,1}_{J^{\infty}(M_{\Sigma})}$,
known as a variational 1-form, and some
$E_{j}\in\Omega^{d,0}_{J^{\infty}(M_{\Sigma})}$.

Since representation (1.6.2) is unique [Di],  and transformation
(1.6.1b) leaves $E^{(i)}_{j}$ unaffected (because
$[\delta,d_{\rho}]=0$, see 1.3), associated to $S$ there arises the
sheaf of Euler-Lagrange ideals
$$
 \CJ_{L}=<D_{\Sigma }E_{1}, D_{\Sigma }E_{2},...,D_{\Sigma
 }E_{n}>\subset\CO_{J^{\infty}(M_{\Sigma})}.\eqno{(1.6.3)}
$$
 Let $Sol_{L}\subset
J^{\infty}(M_{\Sigma})$ be the corresponding zero locus. We will
assume that $\CJ_{L}$ is of pseudo Cauchy-Kovalevskaya type and
usually (1.4.2) will hold.

The variational 1-form $\gamma^{(i)}$ is not quite uniquely defined,
but due to the well-known acyclicity theorem [T,Di], locally it is
determined up to a $d_{\rho}$-exact term.  Therefore,
$$
\omega^{(i)}=\delta\gamma^{(i)}.
$$
unambiguously defines a  section of the quotient sheaf
$\Omega^{d-1,2}_{J^{\infty}(M_{\Sigma})}/d_{\rho}\Omega^{d-2,2}_{J^{\infty}(M_{\Sigma})}$
over $U_{i}$. Transformation (1.6.1b) leaves it invariant, and
definition (1.6.1a) shows that the collection
$$
\omega_{L}=\{\omega^{(i)}=\delta\gamma^{(i)}\}. \eqno{(1.6.4)}
$$
unambiguously defines a global section of the quotient sheaf
$\Omega^{d-1,2}_{J^{\infty}(M_{\Sigma})}/d_{\rho}\Omega^{d-2,2}_{J^{\infty}(M_{\Sigma})}$.
Hence on $Sol_{L}$ there arises the sheaf of Lie algebras
$\CH^{\omega_{L}}_{Sol}$, see Definition 1.5.3 and remark 1.5.3.1.
Our task now is to detect inside it a subalgebra of integrals of
motion. As we have seen already, the nature of the argument tends to
be purely local, and until further notice it will be assumed that
$$
L\in
\Gamma(J^{\infty}(M_{\Sigma}),\Omega^{d,0}_{J^{\infty}(M_{\Sigma})}).
$$

{\bf 1.6.2.} A symmetry of $L$ is an evolutionary vector field $\xi$
s.t.
$$
\text{Lie}_{\xi}L=d_{\rho}\alpha_{\xi}, \eqno{(1.6.5)}
$$
for some $\alpha_{\xi}\in \Omega^{d-1,0}_{J^{\infty}(M_{\Sigma})}$.

An integral of motion is an $F\in
\Omega^{d-1,0}_{J^{\infty}(M_{\Sigma})}$ s.t.
$$
d_{\rho}F\in\CJ_{L}.\eqno{(1.6.6)}
$$
If $\xi$ is a symmetry of $L$ with characteristic $\{Q^{j}\}$, see
(1.2.2), then the computation
$$
d_{\rho}\alpha_{\xi}=\iota_{\xi}\delta
L\buildrel(1.6.2)\over=-\iota_{\xi}d_{\rho}\gamma+E_{j}Q^{j}=d_{\rho}\iota_{\xi}\gamma+E_{j}Q^{j}
$$
shows that $\alpha_{\xi}-\iota_{\xi}\gamma$ is an integral of
motion.

{\bf 1.6.3. N\"other's Theorem.} {\it If $\xi$ is a symmetry of $L$,
then
$$
\xi\mapsto\alpha_{\xi}-\iota_{\xi}\gamma \eqno{(1.6.7)}
$$
establishes a 1-1 correspondence between symmetries of $L$ and
integrals of motion.}

Let $\CI_{L}$ be the space of integrals of motion. Thanks to
N\"other's Theorem its elements are pairs ``integral of motion,
symmetry of $L$'' to be denoted,  depending on context, as either
$(F_{\xi},\xi)$ or $(F,\xi_{F})$.

\bigskip
{\bf 1.6.4.} It is easy to derive from (1.6.5) that any symmetry of
$L$ preserves $\CJ_{L}$, see [Dickey], and thus defines a vector
field on $Sol_{L}$. Motivated by this, let us restrict all the
relevant objects to $Sol_{L}\subset J^{\infty}(M_{\Sigma})$; thus
from now on $\gamma$, $\omega_{L}$, $\CI_{L}$ will denote the
pull-backs of these on $Sol_{L}$.

If $\xi$ is a symmetry of $L$ and $G$ is an integral of motion, then
$\xi G$ is also an integral of motion because, due to (1.4.4),
$$
d_{\rho}\xi G=\xi d_{\rho}G=0.
$$
It is then natural to define a bracket
$$
\CI_{L}\times\CI_{L}\rightarrow\CI_{L};\;(F,\xi_{F}),(G,\xi_{G})\mapsto(\xi_{F}G,[\xi_{F},\xi_{G}]),
\eqno{(1.6.8)}
$$
hoping that it  makes $\CI_{L}$ into a Lie algebra. It does not, but
it agrees with the bracket on $Sol_{L}$ discussed in 1.5 for
$\omega=\omega_{L}$.

{\bf 1.6.5. Lemma.}

{\it (i) $\CI_{L}$ is a subspace of $\Gamma(Sol,\CF)$, see
definition (1.5.4)

(ii) Definition (1.6.8) makes $\CI_{L}$ into a subalgebra of
$\Gamma(Sol,\CF)$, c.f. Proposition 1.5.2.}

\bigskip

{\bf Proof.} (i) is essentially   [Di] 19.6.17 or [DF] Proposition
2.76, which assert that if $\xi$ is a symmetry of $L$ for which
(1.6.5) holds, then
$$
\text{Lie}_{\xi}\gamma=\delta \alpha_{\xi}+d_{\rho}\beta
$$
for some $\beta$. An application of $\delta$ to both sides of this
equality shows that $\xi$ is homotopy Hamiltonian, see (1.5.2). The
corresponding integral of motion
$F_{\xi}=\alpha_{\xi}-\iota_{\xi}\gamma$ is admissible because
$$
\delta
F_{\xi}=\delta\alpha_{\xi}-\delta\iota_{\xi}\gamma=\delta\alpha_{\xi}-\text{Lie}_{\xi}\gamma+\iota_{\xi}\delta\gamma=-d_{\rho}\beta
+\iota_{\xi}\omega_{L}.
$$
Hence $(F_{\xi},\xi)\in\CF$.

(ii) is obvious.

\bigskip

{\bf 1.6.6. Corollary.} {\it If $\tilde{\CI}_{L}$ is the homomorphic
image of $\CI_{L}$ in $\Gamma(Sol,\CH^{\omega_{L}}_{Sol_{L}})$, see
Definition 1.5.3, then
$\tilde{\CI}_{L}\subset\Gamma(Sol,\CH^{\omega_{L}}_{\Sigma})$ is a
Lie subalgebra.}

\bigskip
It is this Lie algebra that is often referred to as the algebra of
integrals of motion or current algebra.

{\bf 1.6.7.} Let us now drop the requirement that $L$ be globally
defined. The exposition above has to be altered a little. An
evolutionary vector field is a symmetry of $L$, see (1.6.1), if it
is of each $L^{(i)}$:
$$
\text{Lie}_{\xi}L^{(i)}=d_{\rho}\alpha_{\xi}^{(i)}.
$$
There may arise discrepancies
$\alpha_{\xi}^{(i)}-\alpha_{\xi}^{(j)}$ on double intersections
$U_{i}\cap U_{j}$, but (1.6.1a)and (1.4.4) ensure that they are
$d_{\rho}$-exact. Therefore, while the collection
$$
\{F_{\xi}^{(i)}=\alpha_{\xi}^{(i)}-\iota_{\xi}\gamma^{(i)}\}
$$
does not define a global section of
$\Omega^{d-1,0}_{J^{\infty}(M_{\Sigma})}$, taken modulo $d_{\rho}$
it defines a global section of
$\Omega^{d-1,0}_{J^{\infty}(M_{\Sigma})}/d_{\rho}\Omega^{d-2,0}_{J^{\infty}(M_{\Sigma})}$.
Therefore, while $\CI_{L}$ does not make sense if $L$ is not
globally defined, its homomorphic image $\tilde{\CI}_{L}$ inside
$\Gamma(Sol,\CH^{\omega_{L}}_{Sol})$ does. Hence

\bigskip

{\bf 1.6.8. Corollary.} {\it For any Lagrangian (1.6.1a,b),
$\tilde{\CI}_{L}\subset\Gamma(Sol,\CH^{\omega_{L}}_{Sol})$ is a Lie
subalgebra.}

\bigskip

Along with $\CH^{\omega_{L}}_{Sol}$, there is its relative version,
$\CH^{\omega_{L}}_{Sol/\Sigma''}$ and the Lie algebra sheaf morphism
$$
\CH^{\omega_{L}}_{Sol}\rightarrow \CH^{\omega_{L}}_{Sol/\Sigma''}
$$
defined in (1.5.11), which is neither surjection nor injection,
generally speaking.

\bigskip

{\bf 1.6.9. Lemma.} {\it If $Sol_{L}$ satisfies (1.4.2), then the
composition
$$
\tilde{\CI}_{L}\hookrightarrow
\Gamma(Sol,\CH^{\omega_{L}}_{Sol})\rightarrow
\Gamma(Sol,\CH^{\omega_{L}}_{Sol/\Sigma''})
$$
is an injection.}

\bigskip

{\bf Proof.} Assume that $\tilde{F}\in \tilde{\CI}_{L}$ is
annihilated by the composite map. This means that if $F\in \CI_{L}$
is a representative of $\tilde{F}$, then $F=F^{o}\wedge d\tau$, and
$d_{\rho/\Sigma''}F^{o}=0$. Due to the Takens acyclicity theorem [T]
(applicable thanks to (1.4.2)), $F^{o}=d_{\rho/\Sigma''}G$ for some
$G$. Therefore, $F=\pm d_{\rho}(Gd\tau)$ and $\tilde{F}=0$, as
desired.

\bigskip

Now we would like to explain that for an important class of
Lagrangians, the sheaf $\CH^{\omega_{L}}_{Sol/\Sigma''}$ is
isomorphic to the canonical $\CH^{can}$, see (1.5.14), and exhibit
some concrete Lie algebras of integrals of motion.

{\bf 1.6.10.}{\it Order 1 Lagrangians and the Legendre transform.}
Let us assume that the Lagrangian $L$ depends only on 1-jets of the
coordinates $x^{j}$. If we let
$$
L=\tilde{L}d\sigma^{0}\wedge\cdots\wedge d\sigma^{d-1},
$$
then (1.6.2) becomes
$$
\aligned
 \delta L = &-d_{\rho}\left((-1)^{p+1}
\left(\frac{\partial
\tilde{L}}{\partial(\partial_{\sigma^{p}}x^{j})}\right) \delta
x^{j}\wedge d\sigma^{0}\wedge\cdots\wedge
\widehat{d\sigma^{p}}\wedge\cdots\wedge d\sigma^{d-1}\right)
\\
&+ \left(\frac{\partial \tilde{L}}{\partial
x^{j}}-\partial_{\sigma^{p}}\left(\frac{\partial
\tilde{L}}{\partial(\partial_{\sigma^{p}}x^{j})}\right)\right)\delta
x^{j}\wedge d\sigma^{0}\wedge\cdots\wedge d\sigma^{d-1},
\endaligned
\eqno{(1.6.9)}
$$
where $\widehat{ }$ means that the term is omitted and summation
w.r.t repeated indices is assumed.

Assume now that on $\Sigma$ there is a distinguished coordinate, say
$\tau=\sigma^{0}$, such that $L$ is a convex function of jets of
coordinates in the $\tau$-direction. It follows then that $Sol_{L}$
satisfies (1.4.2). Applying (1.3.4) to $\gamma$ we obtain
$$
\gamma':=\gamma|_{d\tau=0}= \left(\frac{\partial
\tilde{L}}{\partial(\partial_{\tau}x^{j})}\right) \delta x^{j}\wedge
d\sigma^{1}\wedge\cdots\wedge d\sigma^{d-1}.
$$
Note that, as a function of $\partial_{\tau}x^{j}$, $\tilde{L}$ is
canonically a function on the tangent space $TM$. It follows that
$\gamma'$ is unambiguously  a 1-form on $TM$.

The convexity of $L$ implies that the {\it Legendre transform}
$$
d_{TM}\tilde{L}: TM\rightarrow T^{*}M\eqno{(1.6.10)}
$$
is a diffeomorphism.  A moment's thought shows that $\gamma'$ is the
pull-back of the canonical 1-form on $T^{*}M$ w.r.t.
$d_{TM}\tilde{L}$, which places us in the situation of Lemma
1.5.5.1. In a coordinate form, we have: if $x^{j}$ are coordinates
on $M$, $x_{j}=\partial/\partial x^{j}$ are fiberwise linear
functions on $T^{*}M$, then
$$
(d_{TM}\tilde{L})^{\#}(x^{j})=x^{j},\;(d_{TM}\tilde{L})^{\#}(x_{j})=\frac{\partial
\tilde{L}}{\partial(\partial_{\tau}x^{j})},
$$
and
$$
\aligned  \gamma'&=(d_{TM}\tilde{L})^{\#}\left(x_{j}\delta
x^{j}\wedge
d\sigma^{1}\wedge\cdots\wedge d\sigma^{d-1}\right),\\
\omega_{L}'=(d_{TM}\tilde{L})^{\#}\left(\delta \gamma'\right)&=
(d_{TM}\tilde{L})^{\#}\left(\delta x_{j}\wedge\delta x^{j}\wedge
d\sigma^{1}\wedge\cdots\wedge d\sigma^{d-1}\right),
\endaligned
\eqno{(1.6.11)}
$$
are the pull-backs of the canonical degenerate symplectic form.
Hence  Lemmas 1.5.5.1 and 1.6.9 specialized to the present situation
read as follows.

\bigskip

{\bf 1.6.10.1. Lemma.} {\it If $L$ depends only on the 1-jets of
coordinates and is convex, then in the case where $d=2$, then there
are  the following  Lie algebra (sheaf) morphisms
$$
\CH^{\omega_{L}}_{Sol/\Sigma''} \iso
(d_{TM}\tilde{L})^{-1}\CH^{can},\;
\tilde{\CI}_{L}\hookrightarrow\Gamma(M,
\CH^{\omega_{L}}_{Sol/\Sigma''}) \iso\Gamma(M,\CH^{can}).
$$}

\bigskip
{\bf 1.6.10.2.} This lemma explains the universality of $\CH^{can}$.
One can argue, therefore, that the Lie algebra content of the
``theory'' is independent of the Lagrangian. What captures the
properties of an individual Lagrangian is the subalgebra of
integrals of motion. For example, if $L$ is independent of $\tau$,
the intrinsic time, then $\rho(\partial_{\tau})$ is a symmetry of
$L$, and (1.6.7) produces the corresponding integral of motion as
follows:

since
$$
\rho(\partial_{\tau}) L= d_{\rho}(\tilde{L}
d\sigma^{1}\wedge\cdots\wedge d\sigma^{d-1}),
$$
 the corresponding integral of motion, upon restriction to the
 fibers of $Sol_{L}\rightarrow\Sigma''$,
 becomes
$$
 H_{\rho(\partial_{\tau})}=
 \alpha_{\rho(\partial_{\tau})}-\iota_{\rho(\partial_{\tau})}\gamma'
 =\left(\tilde{L}-\left(\frac{\partial
\tilde{L}}{\partial(\partial_{\tau}x^{j})}\right) \partial_{\tau}
x^{j}\right)d\sigma^{1}\wedge\cdots\wedge d\sigma^{d-1},
\eqno{(1.6.12)}
$$
which is the familiar energy function, of course.
\bigskip

{\bf 1.6.11.} {\it Bosonic string, left/right movers, and a rudiment
of generalized geometry.}

Let $M$ be a Riemannian manifold with metric $(.,.)$, $\Sigma$ be
2-dimensional with coordinates $\tau$ and $\sigma$. By definition, a
point in $ J^{1}(M_{\Sigma})$ is a triple $(t,x,\partial x)$, where
$t\in\Sigma$, $x\in M$, and $\partial x$ is a linear map
$$
\partial x: T_{t}\Sigma\rightarrow T_{x}M,\; \xi\mapsto\partial_{\xi}x.
$$
This makes sense out of the symbol
$(\partial_{\xi}x,\partial_{\eta}x)$ as a function on
$J^{1}(M_{\Sigma})$. The following
$$
L=\frac{1}{2}\left(\left(\partial_{\sigma}-\partial_{\tau}\right)x,
\left(\partial_{\sigma}+\partial_{\tau}\right)x\right)d\sigma\wedge
d\tau \eqno{(1.6.13)}
$$
is then a well-defined Lagrangian, the celebrated {\it
$\sigma$-model Lagrangian}. In terms of local coordinates
$x^{1},...,x^{n}$ s.t. $(.,.)=g_{ij}dx^{i}dx^{j}$ it looks as
follows:
$$
L=\frac{1}{2}\left(g_{ij}\partial_{\sigma}x^{i}\partial_{\sigma}x^{j}
-g_{ij}\partial_{\tau}x^{i}\partial_{\tau}x^{j}\right)d\sigma\wedge
d\tau.
$$

A direct computation shows that
$$
\aligned
 \delta L=&-d_{\rho}\left(\left(\partial_{\tau}x,\delta
 x\right)d\sigma-\left(\partial_{\sigma}x,\delta
 x\right)d\tau\right)\\
 &+\left(\nabla_{\partial_{\tau}x}\partial_{\tau}x-\nabla_{\partial_{\sigma}x}\partial_{\sigma}x\right)
 d\sigma\wedge
d\tau,
 \endaligned
 \eqno{(1.6.14)}
 $$
 where $\nabla_{\partial_{\bullet}x}$ is the value of the Levi-Civita connection
 on $\partial_{\bullet}x$. It is clear that $L$ satisfies
  all the conditions of Lemma 1.6.10.1.

 The Lagrangian being independent of $\tau$ and $\sigma$,
 associated to $\rho(\partial_{\tau})$ and $\rho(\partial_{\tau})$ there
 arise two integrals of motion, energy and  momentum, and any linear
 combination thereof.
 But much more is true. In fact, any vector field of the type
 $$
 \text{either }
 \xi^{-}=\frac{1}{2}f(\sigma-\tau)\rho(\partial_{\sigma}-\partial_{\tau})
 \text{ or }
 \xi^{+}=\frac{1}{2}f(\sigma+\tau)\rho(\partial_{\sigma}+\partial_{\tau})
 \eqno{(1.6.15)}
 $$
 is a symmetry of $L$. Indeed, precisely because
 $(\partial_{\sigma}\pm\partial_{\tau})(\sigma\mp\tau)=0$, one has
 $$
 \aligned
 \xi^{-}L&=d_{\rho}\left(\frac{1}{4}f(\sigma-\tau)\left(\left(\partial_{\sigma}-\partial_{\tau}\right)x,
\left(\partial_{\sigma}+\partial_{\tau}\right)x\right)\left(d\sigma+d\tau\right)\right),\\
\xi^{+}L&=-d_{\rho}\left(\frac{1}{4}f(\sigma+\tau)\left(\left(\partial_{\sigma}-\partial_{\tau}\right)x,
\left(\partial_{\sigma}+\partial_{\tau}\right)x\right)\left(d\sigma-d\tau\right)\right)
\endaligned
\eqno{(1.6.16)}
$$
Using (1.6.7) one obtains the corresponding integrals of motion,
inside $\CH^{\omega_{L}}_{Sol_{L}/\Sigma''}$,
$$
\aligned
F_{\xi^{-}}&=\frac{1}{4}f(\sigma-\tau)\left(\left(\partial_{\sigma}-\partial_{\tau}\right)x,
\left(\partial_{\sigma}-\partial_{\tau}\right)x\right)d\sigma,\\
F_{\xi^{+}}&=-\frac{1}{4}f(\sigma+\tau)\left(\left(\partial_{\sigma}+\partial_{\tau}\right)x,
\left(\partial_{\sigma}+\partial_{\tau}\right)x\right)d\sigma.
\endaligned
\eqno{(1.6.17)}
$$
Upon  Legendre transform (1.6.10), which in terms of local
coordinates  is this
$$
x_{i}=g_{i\alpha}\partial_{\tau}x^{\alpha},\;\partial_{\tau}x^{i}=g^{i\alpha}x_{\alpha},
$$
formulas (1.6.17) become
$$
\aligned
F_{\xi^{-}}&=f(\sigma-\tau)\left(\frac{1}{4}g^{ij}x_{i}x_{j}
+\frac{1}{4}g_{ij}\partial_{\sigma}x^{i}\partial_{\sigma}x^{j}
-\frac{1}{2}x_{j}\partial_{\sigma}x^{j}\right)
d\sigma,\\
 F_{\xi^{+}}&=f(\sigma+\tau)\left(-\frac{1}{4}g^{ij}x_{i}x_{j}
-\frac{1}{4}g_{ij}\partial_{\sigma}x^{i}\partial_{\sigma}x^{j}
-\frac{1}{2}x_{j}\partial_{\sigma}x^{j}\right) d\sigma,
\endaligned
\eqno{(1.6.18)}
$$
and this computes the image of $F_{\xi^{\pm}}$ under the composite
map of Lemma 1.6.10.1. Let
$$
\Vir^{\pm}=\text{span}\{\overline{F_{\xi^{\pm}}}\}\subset
\Gamma(M,\CH^{can}).\eqno{(1.6.19)}
$$
All of this means that the space of global sections of the sheaf of
Lie algebras $\CH^{can}$ contains 2 commuting copies of the Lie
algebra of vector fields on $\Sigma$. In the case where
$\Sigma=S^{1}\times\Sigma''$, each is the centerless Virasoro
algebra, hence the notation. In view of canonical commutation
relations discussed in 1.5.5, formulas (1.6.18) are 2 bozonizations
of the Virasoro algebra -- in the quasiclassical limit.

This prompts the following definitions:

{\bf 1.6.11.1. Definition.}

{\it (i) Denote by $\CH^{\omega_{L},+}_{Sol_{L}/\Sigma''}$ the
centralizer of $\Vir^{-}$ in $\CH^{\omega_{L}}_{Sol_{L}/\Sigma''}$
and call it the right moving algebra.

(ii) Denote by $\CH^{\omega_{L},-}_{Sol_{L}/\Sigma''}$ the
centralizer of $\Vir^{+}$ in $\CH^{\omega_{L}}_{Sol_{L}/\Sigma''}$
and call it the left moving algebra.}

We will present a computation of left/right moving algebra in the
context of the WZW model in sect. 2.9.2. Let us also note that
 $\CH^{\omega_{L},-}_{Sol_{L}/\Sigma''}$ contains yet
 another Virasoro algebra -- the sum of the first two, which upon
 restriction to $\{\tau=0\}$ becomes
 $$
 \aligned
\Vir^{o}&=\text{span}\{\overline{F_{\xi^{+}}}+\overline{F_{\xi^{-}}}\}\\
&=\text{span}\{\overline{f(\sigma)x_{j}\partial_{\sigma}x^{j}d\sigma}\}.
\endaligned
\eqno{(1.6.20)}
$$
Bosonization (1.6.20) is much simpler than (1.6.18) and was
thoroughly investigated in [MSV,GMS1], but the corresponding
Virasoro algebra is neither right nor left moving.

{\bf 1.6.11.2.} {\it Generalized geometry interpretation.}

Formulas (1.6.18) admit a nice, Lagrangian free, interpretation in
the spirit of  Hitchin's ``generalized geometry'', [G]. The idea of
generalized geometry is that the tangent bundle of a manifold must
be consistently replaced with the direct sum of the tangent and
cotangent bundles. From this point of view, a metric on $M$ is a
reduction of the structure group of $TM\oplus T^{*}M$ from $SO(n,n)$
to $SO(n,0)\times SO(0,n)$. Letting $\{e_{i}\}$, $\{e^{j}\}$ be a
pair of relatively dual bases of the $SO(n,0)$-subbundle and letting
$\{f_{i}\}$, $\{f^{j}\}$ the same for the $SO(0,n)$-subbundle, one
can form 2 invariantly defined tensors, $e^{i}e_{i}$ and
$f^{i}f_{i}$.  Noticing that $x_{i}$, in (1.6.18), is naturally
identified with $\partial_{x^{i}}$, $\partial_{\sigma}x^{j}$ with
$dx^{i}$, one concludes that $\Vir^{+}$ is generated by $e^{i}e_{i}$
and $\Vir^{-}$ by $-f^{i}f_{i}$

To talk about these and other issues coherently, one must change
gears and introduce vertex Poisson algebras.

\bigskip\bigskip

\centerline{{\bf 2. Vertex Poisson Algebras}}

\bigskip

Our presentation of this well-known topic, see e.g. [FB-Z], will be
a little different in the following respects. First of all, we will
fix
 $B$, an associative commutative $\BC$-algebra, to be
  the ground ring for all linear algebra constructions of
this section.

Second of all, we will let $\fg=\text{Der}B$ and demand that all the
structures be $\fg$-equivariant. These assumptions are intended to
handle  functions of $\tau$ and $\sigma$ should they appear.
Therefore, two examples to be kept in mind are these:
$$
B=C^{\infty}(\Sigma),\fg=\CT_{\Sigma}(\Sigma)\text{ or }
B=\BC,\fg=0.\eqno{(2.1)}
$$
The case at hand, where $M_{\Sigma}=M\times\Sigma$, is rather
special, and we could have avoided including $B$ and $\fg$ as part
of data (which is customary in works on vertex algebras), but we
decided against it. That the natural setting for what follows is
equivariant was  pointed out by Beilinson and Drinfeld [BD,3.9].

{\bf 2.1. Definition.} An equivariant vertex Poisson algebra is a
collection $(V,T,_{(n)},\fg;n\geq -1)$, where $V$ is a $B$-module,
$$
T:V\rightarrow V
$$
 is  a $B$-linear map, and
 $$
 _{(n)}:V\otimes V\rightarrow V,\; a_{(n)}b=0\text{ if }n>>0
 $$
  is a family of $B$-bilinear multiplications, such that the
following axioms hold:

{\bf I.} The triple $(V,T,_{(-1)})$ is a commutative associative
algebra with derivation $T$.

{\bf II.} The collection $(V,T,_{(n)}; n\geq 0)$ is a vertex Lie
algebra, i.e., the following holds:

\indent\indent{\bf II.1} skew-commutativity
$$
a_{(n)}b=(-1)^{n+1}\sum_{j=0}^{\infty}\frac{)(-1)^{j}}{j!}T^{j}(b_{(n+j)}a),
$$
\indent\indent{\bf II.2.} Jacobi identity
$$
a_{(m)}b_{(n)}c-b_{(n)}a_{(m)}c=\sum_{j=0}^{\infty}\binom
mj(a_{(j)}b)_{(n+m-j)}c,
$$
\indent\indent{\bf II.3.} properties of $T$:
$$
(Ta)_{(n)}=[T,a_{(n)}]=-na_{(n-1)}.
$$
\indent{\bf III.} Leibnitz identity: for any $n$, $a_{(n)}$ is a
derivation of $_{(-1)}$.

{\bf IV.}$\fg$-equivariance: $V$ is a $\fg$-module, and the maps
$_{(n)}$ and $T$ are $\fg$-module morphisms.

\bigskip

In addition, we will always be assuming that a vertex Poisson
algebra $(V,T,_{(n)};n\geq -1)$ is $\BZ_{+}$-graded, i.e.,
$$
V=\bigoplus_{n=0}^{\infty}V_{n},\; T(V_{n})\subset
V_{n+1},\;\fg(V_{n})\subset V_{n},\; V_{m (j)}V_{(n)}\subset
V_{(m+n-j-1)}.\eqno{(2.1.1)}
$$
We will unburden the notation by letting $V$ stand for
$(V,T,_{(n)},\fg ;n\geq -1)$ when this does not lead to confusion
and by suppressing $_{(-1)}$ so that $ab$ stands for $a_{(-1)}$. We
will also  tend to drop the adjective ``equivariant'' whenever it
seems appropriate.

Note that if $m=n=0$, then {\bf II.2} becomes
$$
a_{(0)}b_{(0)}c-b_{(0)}a_{(0)}c=\left(a_{(0)}b\right)_{(0)}c,
\eqno{(2.1.2)}
$$
which is the usual Jacobi identity for $(V,_{(0)})$.
Anticommutativity fails, but {\bf II.1} ensures that it holds up to
$T(...)$. This almost proves the following  important

\bigskip
{\bf 2.2. Lemma.} {\it If $V$ is a vertex Poisson algebra, then
$T(V)\subset V$ is a 2-sided ideal w.r.t. $_{(0)}$, and
$(V/T(V),_{(0)})$ is a Lie algebra.}

\bigskip

{\bf 2.3. Tensor products}

\bigskip
 The simplest example of a vertex Poisson
algebra is a commutative associative algebra $V$ with derivation
$T$. Defining $a_{(-1)}$ to be  multiplication by $a$ and letting
$a_{(n)}=0$ if $n\geq 0$. makes $V$ into a vertex Poisson algebra
with $T=\partial$.

If $(V_{1},_{(n)} T_{1})$ and $(V_{2},_{(n)}, T_{2})$ are two vertex
Poisson algebras, then $V_{1}\otimes V_{2}$ carries at least two
vertex Poisson algebra structures. First of all, one can simply
regard $V_{1}\otimes V_{2}$ as an extension of scalars whereby
$V_{1}\otimes V_{2}$ becomes a vertex Poisson algebra over $V_{1}$
with derivation $T_{(2)}$ and multiplications coming from $V_{2}$.

Second of all, one can define $T=T_{1}+T_{2}$ and
$$
(a\otimes b)_{(n)}=\left\{\aligned a_{(-1)}b_{(-1)}&\text{ if }
n=-1\\
\sum_{i=0}^{\infty}\frac{1}{i!}\left(\left(T_{1}^{i}a\right)_{(-1)}b_{(n+i)}+a_{(n+i)}\left(T_{1}^{i}b\right)_{(-1)}\right)
&\text{ if } n\geq 0
\endaligned\right.
\eqno{(2.3.1)}
$$
If, in addition, $V_{1}$ is of the type we started with, i.e.,
$(V_{1})_{(n)}(V_{1})=0$ if $n\geq 0$, then (2.3.1) is simplified as
follows
$$
(a\otimes b)_{(n)}=\left\{\aligned a_{(-1)}b_{(-1)}&\text{ if }
n=-1\\
\sum_{i=0}^{\infty}\frac{1}{i!}\left(\left(T_{1}^{i}a\right)_{(-1)}b_{(n+i)}\right)
&\text{ if } n\geq 0
\endaligned\right.
\eqno{(2.3.2)}
$$
In a sense, the second version is a twist of the first by derivation
$T_{1}\in\text{Der}(V_{1})$. In the context of equivariant vertex
Poisson algebras this can be generalized as follows.

If $(V,_{(n)}, T)$ is an equivariant vertex Poisson algebra over $B$
and $\xi\in\fg$, then letting
$$
a_{(n)_{\xi}}=\sum_{i=0}^{\infty}\frac{1}{i!}\left(\xi^{i}a\right)_{(n+i)}
\eqno{(2.3.3)}
$$
defines a vertex Poisson algebra $(V,_{(n)_{\xi}},T+\xi)$. We will
refer to this construction as the $\xi$-twist. Note that the
$\xi$-twist reduces the constants from $B$ to the algebra of
$\xi$-invariants,$B^{\xi}$.

\bigskip

{\bf 2.4. From Poisson vertex algebras to Courant algebroids}

\bigskip
 The
Poisson vertex algebra structure on $V=\oplus_{n=0}^{\infty}V_{n}$
defines on the subspace $V_{0}+V_{1}$ the following operations:
$$
_{(-1)}: V_{0}\otimes V_{0}\rightarrow V_{0},\eqno{(2.4.1a)}
$$
$$
_{(-1)}: V_{0}\otimes V_{1}\rightarrow V_{1},\;  V_{1}\otimes
V_{0}\rightarrow V_{1},\eqno{(2.4.1b)}
$$
$$
_{(0)}: V_{1}\otimes V_{0}\rightarrow V_{0},\;  V_{0}\otimes
V_{1}\rightarrow V_{0},\eqno{(2.4.1c)}
$$
$$
_{(0)}: V_{1}\otimes V_{1}\rightarrow V_{1},\eqno{(2.4.1d)}
$$
$$
_{(1)}: V_{1}\otimes V_{1}\rightarrow V_{0},\eqno{(2.4.1e)}
$$
$$
T:V_{0}\rightarrow V_{1}, \eqno{(2.4.1f)}
$$
 all the other operations
either not preserving the subspace $V_{0}+V_{1}$ or being zero due
to condition (2.1.1).

Vertex Poisson algebra axioms imply that (2.4.1a-f) satisfy certain
conditions; e.g., (2.3.1a) is such that $(V_{0},_{(-1)})$ is an
associative commutative $B$-algebra, and (2.4.1b) is such that
$V_{1}$ is a $V_{0}$-module due. In [GMS1], these conditions were
written down explicitly and made into an axiomatic definition of a
vertex algebroid -- in a more complicated, quantum, situation. It is
a nice observation due to Bressler [Bre] that under some
non-degeneracy assumptions a quasiclassical limit of a vertex
algebroid is an exact Courant $V_{0}$-algebroid; e.g. (2.4.1d)
defines the Dorfman barcket [Dor,G] on $V_{1}$. Therefore, the
assignment $V\mapsto \left(V_{0}\oplus V_{1},\; T, _{(-1)},
_{(0)},_{(1)}\right)$ defines a functor from a subcategory of vertex
Poisson algebras to the category of exact Courant
$V_{0}$-algebroids. This functor is actually an equivalence of
categories, and a classification of exact Courant algebroids
furnishes that of a subclass of vertex Poisson algebras. For the
future use, and for the reader's convenience -- after all the
present situation is somewhat different --  let us now reproduce the
essence of this argument.

\bigskip

{\bf 2.4.1.} We have seen already that the pair $(V_{0},_{(-1)})$ is
an associative commutative $B$-algebra. Let $A=V_{0}$. The entire
$V$, hence $V_{1}$, is  an $A$-module and $A_{(n)}A=0$ if $n\geq 0$.

\bigskip

By virtue of Axiom I, the map
$$
T: A\rightarrow AT(A)\subset V_{1}\text{ is a $B$-derivation, }
\eqno{(2.4.2)}
$$
 i.e., $T(ab)=aT(b)+bT(a)$ and $T(B)=0$. Therefore, $AT(A)$ is a
quotient of the module of relative K\"ahler differentials,
$\Omega_{A/B}$.

It is easy to see that
$$
(AT(A))_{(n)}(AT(A))=0,\;n\geq 0. \eqno{(2.4.3)}.
$$

\bigskip

{\it Assumption 1.} {\it Let $(A;T: A\rightarrow AT(A))$ be
isomorphic to $(A;d:A\rightarrow \Omega_{A/B})$.}

\bigskip

There arises an exact sequence of $A$-modules
$$
0\rightarrow \Omega_{A/B}\rightarrow V_{1}\rightarrow
V_{1}/\Omega_{A/B}\rightarrow 0.
 \eqno{(2.4.4)}
$$
Let $\CT=V_{1}/\Omega_{A/B}$. It is an $A$-module and a Lie algebra
w.r.t. the operation $_{(0)}$, by virtue of Lemma 2.2. Furthermore,
the map
$$
_{(0)}:\CT\otimes A\rightarrow A \eqno{(2.4.5)}
$$
gives $A$ a $\CT_{A/B}$-module structure (see (2.1.2)) compatible
with the $A$-module structure in that $(a\xi)_{(0)}b=a(\xi_{(0)}b$.

 For each
$\tau\in\CT$, $\tau_{(0)}\in\text{End}(A)$ is actually a
$B$-derivation of $A$, and this defines a Lie algebra homomorphism
over $A$
$$
\CT\rightarrow\text{Der}_{B}(A), \eqno{(2.4.6)}
$$
All of this can be summarized by saying that $\CT$ is an
$A$-algebroid Lie.

\bigskip
{\it Assumption 2.} {\it  Morphism (2.4.6) is an isomorphism.}

\bigskip

Consider the map
$$
_{(0)}: \CT\otimes\Omega_{A/B}\rightarrow\Omega_{A/B} \eqno{(2.4.7)}
$$
arising by virtue of (2.4.3). It is defined by the Lie derivative:
$$
\xi_{(0)}\omega=\text{Lie}_{\xi}\omega,\eqno{(2.4.8)}
$$
cf. (1.5.20d). (Indeed,
$\xi_{(0)}(aTb)=(\xi_{(0)}a)Tb+a(T\xi_{(0)}b)$.)

Next, again thanks to (2.3.3), there arises the map
$$
_{(1)}:\CT\otimes\Omega_{A/B}\rightarrow A. \eqno{(2.4.9)}
$$
It is the natural pairing of vector fields and forms:
$$
\xi_{(1)}\omega=\iota_{\xi}\omega. \eqno{(2.4.10)}
$$
(Indeed,
$\xi_{(1)}(aTb)=(\xi_{(1)}a)Tb+a(\xi_{(1)}Tb)=a(\xi_{(0)}b)$, where
 axioms II.3 and III are used.)

 This determines all of  (2.3.1a-f) that makes sense on the graded
 object $A\oplus(\CT\oplus\Omega_{A/B})$. To continue our analysis we need to
 make the following

 \bigskip
 {\it Assumption 3. Let sequence (2.4.4) be splitting.}

 \bigskip

 Let us fix a splitting
 $$
 s:\CT\rightarrow V_{1}.\eqno{(2.4.11)}
 $$
 Then there arise the following two maps:
 $$
 _{(1)_{s}}:\CT\otimes\CT\rightarrow A,\eqno{(2.4.12)}
 $$
 $$
 _{(0)_{s}}:\CT\otimes\CT\rightarrow \Omega_{A/B},\eqno{(2.4.13)}
 $$
 where (2.3.12) is the restriction of $_{(1)}$ to $s(\CT)$, and
 (2.3.13) is the composition of the restriction of $_{(0)}$ to
 $s(\CT)$ with the projection $V_{1}\rightarrow\Omega_{A/B}=V_{1}/s(\CT)$.
 These two maps determine all of (2.4.1a-f).

 The map $_{(1)_{s}}$ is, in fact, a symmetric $A$-bilinear form on $\CT$. By
 varying the splitting $s$ it can killed. Indeed, letting
 $h(.,.)=_{(1)_{s}}$, we obtain, for any $\xi\in\CT$, an $A$-linear
form $h(\xi,.)\in\Omega_{A/B}$. Replacing $s$ with $s_{h}$ defined
to be
$$
s_{h}(\xi)=s(\xi)-\frac{1}{2}h(\xi,.)
$$
we get $_{(1)_{s_{h}}}=0$.

Therefore, we can, and usually will,  assume that
$$
V_{1}=\CT\oplus\Omega_{A/B}\eqno{(2.4.14)}
$$
and
$$
_{(1)}:\left(\CT\oplus\Omega_{A/B}\right)\otimes\left(\CT\oplus\Omega_{A/B}\right)\rightarrow
A\eqno{(2.4.15)}
$$
is the canonical pairing
$(\xi+\omega)_{(1)}(\xi'+\omega')=\iota_{\xi}\omega'+\iota_{\xi'}\omega$,
cf. (2.4.10).

Sometimes the following version of (2.4.15) will be used: let $h$ be
a symmetric $A$-bilinear form on $\CT$ and define
$$
\aligned
_{(1)}:\left(\CT\oplus\Omega_{A/B}\right)\otimes\left(\CT\oplus\Omega_{A/B}\right)\rightarrow
A\\
(\xi+\omega)_{(1)}(\xi'+\omega')=\iota_{\xi}\omega'+\iota_{\xi'}\omega+h(\xi,\xi').
\endaligned
\eqno{(2.4.15_{h})}
$$

\bigskip
{\bf 2.4.2.}  Therefore, all  moduli, if any, come from
$_{(0)_{s}}$. A short computation shows that it is $A$-linear.
Furthermore, axiom IV implies that
$$
 _{(0)_{s}}\in\text{Hom}_{\fg}\left(\CT\otimes\CT,\Omega_{A/B}\right).\eqno{(2.4.16)}
 $$
Hence $_{(0)_{s}}$ can be considered as an $A$-trilinear
$\fg$-invariant functional on $\CT$, and as such it will be denoted
by $H$:
$$
_{(0)_{s}}\approx H\in\left(\Omega_{A/B}^{\otimes 3}\right)^{\fg}.
\eqno{(2.4.17)}
$$
 Skew-commutativity
II.1 implies that it is anti-commutative in the first 2 variables:
$H(\xi,\eta,.)=-H(\eta,\xi,.)$.

 Jacobi identity II.2 applied to
$[\xi_{(1)},\eta_{(0)}](\zeta)$, $\xi,\eta,\zeta\in s(\CT_{A})$,
shows that, in fact, $H(.,.,.)$, is totally anti-commutative, hence
belongs to $\left(\Omega_{A/B}^{3}\right)^{\fg}$.

 Jacobi identity II.2 applied to $[\xi_{(0)},\eta_{(0)}](\zeta)$,
$\xi,\eta,\zeta\in s(\CT_{A})$, shows that $H$ is closed, i.e.,
$$
H\in \left(\Omega_{A/B}^{3,cl}\right)^{\fg}.\eqno{(2.4.18)}
$$
Conditions (2.4.14-15 or $15_{h}$) do not determine the splitting
$s$; they are respected by the shearing transformation
$$
\CT\ni\xi\mapsto\xi+\iota_{\xi}\alpha\text{ for a fixed }\alpha\in
\left(\Omega_{A/B}^{2}\right)^{\fg} \eqno{(2.4.19)}.
$$
The effect of this transformation on $H$ is
$$
H\mapsto H+d_{DR}\alpha.\eqno{(2.4.20)}
$$

\bigskip

{\bf 2.4.3.} Checking the various properties of maps (2.4.1a-f)
derived in 2.3.1 against the definition of an exact Courant
$A$-algebroid [LWX] (especially in the form proposed in [Bre]) shows
the following. If Assumptions 1--3 hold, then the equivariant
Poisson vertex algebra structure on $V$ defines an equivariant exact
Courant $A$-algebroid structure on $\CT\oplus\Omega_{A}$ such that
$$
_{(0)}:\left(\CT\oplus\Omega_{A}\right)\otimes\left(\CT\oplus\Omega_{A}\right)
\rightarrow \CT\oplus\Omega_{A}
$$
is the Dorfman [Dor,G] bracket,
$$
_{(1)}:\left(\CT\oplus\Omega_{A}\right)\otimes\left(\CT\oplus\Omega_{A}\right)
\rightarrow A
$$
is the symmetric pairing, and (2.4.6) is the anchor.

The discussion in 2.3.2 practically proves (see [Bre, GMS1] for a
complete analysis) that the category of exact equivariant Courant
$A$-algebroids is an
$\left(\Omega_{A/B}^{3,cl}\right)^{\fg}$-torsor. Indeed, if $\CC$ is
one such algebroid and $H\in
\left(\Omega_{A/B}^{3,cl}\right)^{\fg}$, then the $H$-twisted
Courant algebroid $\CC\plus H$ is defined by replacing $_{(0)_{s}}$
with $_{(0)_{s}}+H$; a moment's thought shows that this is
independent of the choice of $s$.  In fact, a ``canonical'' Courant
algebroid $\CC_{0}$ can be chosen by letting the only ``unknown''
operation $_{(0)_{s}}$ be zero:
$$
\CC_{0}\text{ s.t. (2.3.14,15) hold and $_{(0)_{s}}=0$}.
\eqno{(2.4.21)}
$$
This identifies the category of equivariant exact Courant
$A$-algebroids with $\left(\Omega_{A/B}^{3,cl}\right)^{\fg}$ s.t.
$$
\left(\Omega_{A/B}^{3,cl}\right)^{\fg}\ni H\mapsto
\CC_{H}=\CC_{0}\plus H. \eqno{(2.4.22)}
$$

Shearing transformation (2.4.20) implies the following description
of morphisms
$$
\text{Mor}(\CC,\CC\plus H)=\{\alpha\in
\left(\Omega_{A/B}^{2}\right)^{\fg}\text{ s.t. }d_{DR}\alpha=H\}
\eqno{(2.4.23)}
$$
and automorphisms
$$
\text{Aut}(\CC)=\left(\Omega_{A/B}^{2,cl}\right)^{\fg}.\eqno{(2.4.24)}
$$
In particular, the set of isomorphism classes of exact Courant
$A$-algebroids is identified with the $\fg$-invariant de Rham
cohomology group,
$$
\left(\Omega_{A/B}^{3,cl}\right)^{\fg}/d_{DR}\left(\Omega_{A/B}^{2}\right)^{\fg}.
\eqno{(2.4.25)}
$$

\bigskip

{\bf 2.5. Symbols of vertex differential operators}

\bigskip
 Let $\Sigma$ be
an open subset of a $\BR^d$, $U$ of $\BR^{n}$ and
$U_{\Sigma}=U\times\Sigma$. Define $B=C^{\infty}(\Sigma)$,
$\fg=\CT_{\BR^{d}}(\Sigma)$. Identify $\fg$ with the subalgebra of
horizontal vector fields on $U_{\Sigma}$, thereby making
$C^{\infty}(U_{\Sigma})$ into a $\fg$-module. These are the
prerequisites to the definition of an equivariant vertex Poisson
algebra over $B$.

{\bf 2.5.1. Definition.} {\it  Call $V$ an algebra of symbols of
vertex differential operators, SVDO for short, if

(i) $V_{0}=C^{\infty}(U_{\Sigma})$,  $V_{1}$ is an equivariant exact
Courant $C^{\infty}(U_{\Sigma})$-algebroid over
$B=C^{\infty}(\Sigma)$,

 (ii) $V$ is generated as an associative
commutative algebra with derivation $T$ by $V_{0}\oplus V_{1}$.}

\bigskip

We have obtained, therefore, a functor, say $\CF$, from the category
of SVDOs to the category of equivariant exact Courant
$C^{\infty}(U_{\Sigma})$-algebroids:
$$
\CF:\{\text{SVDOs}\}\rightarrow \{\text{Courant algebroids}\}.
\eqno{(2.5.1)}
$$

\bigskip

{\bf 2.5.2. Theorem.} ([GMS1, Bre]). {\it This functor is an
equivalence of categories.}

\bigskip

To be precise, [GMS1, Bre] only construct $\CF^{*}$, the left
adjoint to $\CF$, but a simple representation-theoretic argument
shows that the ``vertex envelope'', $\CF^{*}(\CC)$, is simple.

\bigskip

{\bf 2.6. A sheaf-theoretic version}

\bigskip
 All of this can be spread over
manifolds. The geometric prerequisite is a fiber bundle
$$
\pi: M_{\Sigma}\rightarrow\Sigma\eqno{(2.6.1a)}
$$
with a flat connection
$$
\nabla: \CT_{\Sigma}\rightarrow\CT_{M_{\Sigma}}\eqno{(2.6.1b)}
$$

A sheaf of SVDOs, $\CV$, over  $M_{\Sigma}$ is a sheaf of vector
spaces s.t. the space of sections $\CV(U)$ is an SVDO for each open
$U\subset M_{\Sigma}$ such that
$$
\CV(U)_{0}=\CO_{M_{\Sigma}}(U_{\Sigma}),
B(U_{\Sigma})=\pi^{*}\CO_{\Sigma}(U),\eqno{(2.6.2)}
$$
and the equivariant structure is defined by $\nabla$.

The condition that $\CV(U)_{0}=\CO_{M_{\Sigma}}(U)$ implies that
$\CV$ is automatically a sheaf of $\CO_{M_{\Sigma}}(U)$-modules. It
follows from (2.4.6) that the next homogeneous component, $\CV_{1}$
is an extension of vertical vector fields by relative 1-forms:
$$
0\rightarrow\Omega_{M_{\Sigma}/\Sigma}\rightarrow\CV_{1}\rightarrow\CT_{M_{\Sigma}/\Sigma}\rightarrow
0. \eqno{(2.6.3)}
$$
As to the existence of such sheaves, they are plentiful locally: for
any sufficiently small open $U\subset M_{\Sigma}$, the category of
such sheaves over $U$ is an
$\left(\Omega^{3,cl}_{M_{\Sigma}/\Sigma}\right)^{\nabla}(U)$-torsor,
as follows from (2.4.22) and Theorem 2.5.2. If $\CV_{U}$ is one such
sheaf and
$H\in\left(\Omega^{3,cl}_{M_{\Sigma}/\Sigma}\right)^{\nabla}(U)$,
then
$$
\text{Mor}(\CV_{U},\CV_{U}\plus
H)=\{\alpha\in\left(\Omega^{2}_{M_{\Sigma}/\Sigma}\right)^{\nabla}(U)\text{
s.t. }d_{DR}\alpha=H\}, \eqno{(2.6.4)}
$$
cf. (2.4.23).

Technically, (2.6.4) means that there is a gerbe, i.e. a sheaf of
categories, of SVDOs bound by a sheaf complex
$$
0\rightarrow
\left(\Omega^{2}_{M_{\Sigma}/\Sigma}\right)^{\nabla}(U)\buildrel
d_{DR}\over\rightarrow\left(\Omega^{3,cl}_{M_{\Sigma}/\Sigma}\right)^{\nabla}(U)\rightarrow
0.
$$
so that the categories over sufficiently small $U$ are equivalent to
that of SVDOs.

{\it A priori} there may be no single sheaf of SVDOs on the entire
$M$; an obstruction to its existence is a certain canonical
characteristic class lying in
$$
H^{2}(M_{\Sigma},
\left(\Omega^{2}_{M_{\Sigma}/\Sigma}\right)^{\nabla}
\rightarrow\left(\Omega^{3,cl}_{M_{\Sigma}/\Sigma}\right)^{\nabla}).
$$

At this point let us return to the concrete situation of interest to
us where $M_{\Sigma}=M\times\Sigma$ and $\nabla$ is the horizontal
connection. If so,  the above discussion is simplified in that the
sheaves $\left(\Omega^{\bullet}_{M_{\Sigma}/\Sigma}\right)^{\nabla}$
can be replaced with $\Omega^{\bullet}_{M}$. For example, the
obstruction becomes a class lying in
$$
H^{2}(M, \Omega^{2}_{M} \rightarrow\Omega^{3,cl}_{M}).
$$

This class vanishes; the obstruction (equal to the 1st Pontryagin
class) computed in [GMS1], see also [Bre], is a purely quantum
phenomenon, and in any case, an example of such sheaf will be
exhibited shortly.

Furthermore,  (2.4.23-25) imply that the set of isomorphism classes
of such sheaves is an $H^{1}(M,
\Omega^{2}_{M}\rightarrow\Omega^{3,cl}_{M})$-torsor, and the group
of automorphisms of any such sheaf is isomorphic to $H^{0}(M,
\Omega^{2}_{M}\rightarrow\Omega^{3,cl}_{M})\iso
H^{0}(M,\Omega_{M}^{2,cl})$.

Note that since the sequence
$$
0\rightarrow\Omega^{2,cl}_{M}\hookrightarrow \Omega^{2}_{M}
\buildrel d_{DR}\over \rightarrow \Omega^{3,cl}_{M}\rightarrow
0\eqno{(2.6.5)}
$$
is exact, we obtain  isomorphisms
$$
\aligned H^{1}(M, \Omega^{2}_{M}\rightarrow\Omega^{3,cl}_{M})&\iso
H^{1}(M, \Omega^{2,cl}_{M}),\\
 H^{0}(M, \Omega^{2}_{M}\rightarrow\Omega^{3,cl}_{M})&\iso H^{0}(M,
\Omega^{2,cl}_{M}). \endaligned \eqno{(2.6.6)}
$$
The long exact cohomology sequence associated with (2.6.5) implies,
in addition, that
$$
H^{1}(M, \Omega^{2,cl}_{M})\iso H^{0}(M, \Omega^{3}_{M})/d H^{0}(M,
\Omega^{2}_{M})\iso H^{3}(M, \BR),\eqno{(2.6.7)}
$$
where the last isomorphism is the de Rham theorem. This proves

\bigskip

{\bf 2.6.1. Proposition.} {\it a) The set of isomorphism classes of
sheaves of SVDOs on $M_{\Sigma}$ is identified with either of the
isomorphic groups $H^{1}(M, \Omega^{2,cl}_{M})$ and $ H^{3}(M,
\BR)$.

 b) If $\CV$ is a sheaf of SVDOs over $M$, then
$$
\text{Aut}\CV\iso H^{0}(M, \Omega^{2,cl}_{M}).
$$}

{\bf 2.6.2.} Here is an explicit construction of identifications a)
and b) of Proposition 2.5.1. The presentation of the set of
isomorphism classes as $H^{1}(M, \Omega^{2,cl}_{M})$ emphasizes that
locally all such sheaves are isomorphic (this is an immediate
consequence of (2.4.25)). Indeed, let $\{U_{i}\}$ be a covering by
balls. Let $\CV_{i}$ be the restriction $\CV$ to $U_{i}$. Then there
arise canonical identifications,
$$
\phi_{ij}:\CV_{i}|_{U_{i}\cap U_{j}}\iso \CV_{j}|_{U_{i}\cap U_{j}},
\eqno{(2.6.8)}
$$
to be thought of as gluing functions. Let now
$\alpha_{ij}\in\Omega^{2,cl}_{M}(U_{i}\cap U_{j})$ be a \v Cech
cocycle representing  $\alpha\in H^{1}(M, \Omega^{2,cl}_{M})$.
Regarding $\alpha_{ij}$ as an automorphism of $\CV_{j}|_{U_{i}\cap
U_{j}}$, define
$$
\hat{\phi}_{ij}\buildrel \text{def}\over
=\phi_{ij}\plus\alpha_{ij}:\CV_{i}|_{U_{i}\cap U_{j}}\iso
\CV_{j}|_{U_{i}\cap U_{j}}, \eqno{(2.6.9)}
$$
to be the composition of $\phi_{ij}$ and the shear by$\alpha_{ij}$
defined in (2.4.19). The \v Cech cocycle condition satisfied by
$\{\alpha_{ij}$ implies that
$\hat{\phi}_{ik}\circ\hat{\phi}_{kj}\circ\hat{\phi}_{ji}=\text{id}$
on the triple intersection $U_{i}\cap U_{j}\cap U_{k}$ for any
$i,j,k$. Thus $\hat{\phi}_{ij}$ are gluing functions of a new sheaf
of SVDOs, to be denoted $\CV\plus\alpha$.

Contrary to this, the presentation of the set of isomorphism classes
as $H^{3}(M, \BR)$ has nothing to do with gluing functions or even
the $\CO_{M}$-module structure. Indeed, for an element of $H^{3}(M,
\BR)$, pick a global closed 3-form $H$ representing it. The sheaf
$\CV\plus H$ is different from $\CV$ only in that the operation
$$
_{(0)}:\;\CV_{0}\otimes\CV_{0}\rightarrow\CV_{0}
$$
is replaced with $_{(0)}+ H$ as explained in 2.4.3 ( and the sheaf
$\CV\plus d\beta$, $\beta$ a global 2-form, is canonically
isomorphic to $\CV$.)

The relation of one point of view to another is as follows. For
example, given $\CV\plus H$, find a collection
$\beta=\{\beta_{i}\in\Omega^{2}_{M}(U_{i})\}$ so that
$d\beta_{i}=H|_{U_{i}}$. Then $d_{\check{C}}(\beta)$ is de
Rham-closed and hence is a \v Cech 1-cocycle with coefficients in
$\Omega^{2,cl}_{M}$. The map
$$
H^{0}(M,\Omega^{3,cl}_{M})\ni H\mapsto \beta\mapsto \text{class of
}d_{\check{C}}(\beta)\in H^{1}(M,\Omega^{2,cl}_{M}) \eqno{(2.6.10)}
$$
descends to  the inverse of (2.6.7).

Now,  $(\CV\plus H)|_{U_{i}}=\CV_{i}$ as vector spaces but not as
SVDOs; to obtain an SVDO isomorphism, the shear by $\beta_{i}$ is
needed:
$$
\beta_{i}:\; \CV_{i}\rightarrow (\CV\plus H)|_{U_{i}}
\eqno{(2.6.11)}.
$$
The effect of this transformation on the gluing functions is as
follows:
$$
\phi_{ij}\mapsto\phi_{ij}\plus d_{\check{C}}\beta, \eqno{(2.6.12)}
$$
cf. (2.4.20), and this delivers the desired isomorphism
$$
\CV\plus H\iso \CV\plus (\text{ class of }d_{\check{C}}\beta).
\eqno{(2.6.13)}
$$

\bigskip
{\bf 2.7. A natural sheaf of SVDOs.}

\bigskip
 Let us attach to any smooth $M$ a sheaf of SVDOs, which
  depends on $M$ functorially. In order to do so, let us place ourselves in
the situation where $T^{*}M_{\Sigma}=T^{*}M\times\Sigma$, $\Sigma$
satisfies (1.1.2-4) and carries, in particular, a distinguished
coordinate system, $\sigma$ and $\tau$.

Taking advantage of (1.1.7), we note that the operator of the jet
connection, (1.1.5), splits in the vertical and horizontal
components, e.g.,
$$
\rho(\partial_{\sigma})=\partial^{v}_{\sigma}+\partial_{\sigma}^{h},\eqno{(2.7.1)}
$$
where the latter stands for the operator of differentiation w.r.t.
$\sigma$ ``appearing explicitly''.

Let
$$
\pi: J^{\infty}(T^{*}M_{\Sigma/\Sigma''})\rightarrow M\eqno{(2.7.2)}
$$
be the natural projection. There arises the direct image of the
structure sheaf $\pi_{*}\CO_{J^{\infty}(T^{*}M_{\Sigma/\Sigma''})}$
which we will take the liberty to denote also by
$\CO_{J^{\infty}(T^{*}M_{\Sigma/\Sigma''})}$ because this is
unlikely to cause confusion. Thus, for example, if $U\subset M$ is
open, then $\CO_{J^{\infty}(T^{*}M_{\Sigma/\Sigma''})}(U)$ will
stand for the space of sections over $U$.

Being a structure sheaf,
$\CO_{J^{\infty}(T^{*}M_{\Sigma/\Sigma''})}$ carries a canonical
multiplication. Let us define a grading
$$
\aligned
\CO_{J^{\infty}(T^{*}M_{\Sigma/\Sigma''})}=&\bigoplus_{i=0}^{\infty}
\CO_{J^{\infty}(T^{*}M_{\Sigma/\Sigma''})}^{i}\text{ s.t.}\\
 \CO_{J^{\infty}(T^{*}M_{\Sigma/\Sigma''})}^{i}\cdot
\CO_{J^{\infty}(T^{*}M_{\Sigma/\Sigma''})}^{j}&\subset\CO_{J^{\infty}(T^{*}M_{\Sigma/\Sigma''})}^{i+j}
\endaligned
\eqno{(2.7.3)}
$$
by requiring that the pull-back of functions on $M$ have degree 0,
the pull-back of fiberwise linear functions on $T^{*}M$ have degree
1, and the operator $\partial_{\sigma}^{v}$, defined in (2.7.1),
have degree 1, i.e., that
$\partial_{\sigma}^{v}(\CO_{J^{\infty}(T^{*}M_{\Sigma/\Sigma''})}^{i})\subset
\CO_{J^{\infty}(T^{*}M_{\Sigma/\Sigma''})}^{i+1}$. Thus, for
example,
$$
\CO_{J^{\infty}(T^{*}M_{\Sigma/\Sigma''})}^{0}=\CO_{M_{\Sigma}},\;
\CO_{J^{\infty}(T^{*}M_{\Sigma/\Sigma''})}^{1}=\CT_{M_{\Sigma}/\Sigma}\oplus\Omega_{M_{\Sigma}/\Sigma},
\eqno{(2.7.4)}
$$
cf. (2.6.3), where $\CT_{M_{\Sigma}/\Sigma}$ is realized inside
$\CO_{J^{\infty}(T^{*}M_{\Sigma/\Sigma''})}$ as the pull-back of
fiberwise linear functions on $T^{*}M$, and
$\Omega_{M_{\Sigma}/\Sigma}$ is realized as
$\CO_{M}\partial_{\sigma}^{v}\CO_{M}$, cf. sect. 2.4.1, Assumption
1.

\bigskip
{\bf 2.7.1. Proposition.} {\it The sheaf
$\CO_{J^{\infty}(T^{*}M_{\Sigma/\Sigma''})}$ carries a unique
structure of a sheaf of SVDOs over $B=\CO_{\Sigma}$ such that
$_{(-1)}$ is the canonical multiplication, $T=\partial_{\sigma}^{v}$
(which furnishes (2.4.1a,b,f) in this case), and (2.4.1c-e) take the
following form: if $\xi,\xi'\in\CT_{M}$,
$\omega,\omega'\in\Omega_{M}$, then
$$
\aligned
 _{(0)}: \left(\CT_{M_{\Sigma}/\Sigma}\oplus\Omega_{M_{\Sigma}/\Sigma}\right)\otimes
\CO_{M_{\Sigma}}\rightarrow
\CO_{M_{\Sigma}},\;\CO_{M_{\Sigma}}\otimes
\left(\CT_{M_{\Sigma}/\Sigma}\oplus\Omega_{M_{\Sigma}/\Sigma}\right)\rightarrow
\CO_{M_{\Sigma}}\\
\left(\xi+\omega\right)_{(0)}F=-F_{(0)}\left(\xi+\omega\right)=\xi
F,
\endaligned
\eqno{(2.7.5)}
$$
$$
\aligned
 _{(0)}: \left(\CT_{M_{\Sigma}/\Sigma}\oplus\Omega_{M_{\Sigma}/\Sigma}\right)\otimes
\left(\CT_{M_{\Sigma}/\Sigma}\oplus\Omega_{M_{\Sigma}/\Sigma}\right)\rightarrow
\CT_{M_{\Sigma}/\Sigma}\oplus\Omega_{M_{\Sigma}/\Sigma},\\
\left(\xi+\omega\right)_{(0)}\left(\xi'+\omega'\right)=[\xi,\xi']+\text{Lie}_{\xi}\omega'
-\text{Lie}_{\xi'}\omega+\partial_{\sigma}^{v}\left(\iota_{\xi'}\omega\right)
\endaligned
\eqno{(2.7.6)}
$$
$$
\aligned
 _{(1)}: \left(\CT_{M_{\Sigma}/\Sigma}\oplus\Omega_{M_{\Sigma}/\Sigma}\right)\otimes
\left(\CT_{M_{\Sigma}/\Sigma}\oplus\Omega_{M_{\Sigma}/\Sigma}\right)\rightarrow
\CO_{M_{\Sigma}},\\
\left(\xi+\omega\right)_{(1)}\left(\xi'+\omega'\right)=\iota_{\xi}\omega'+\iota_{\xi'}\omega
\endaligned
\eqno{(2.7.7)}
$$}

\bigskip

Note that (2.7.5--7) restricted to some $U\subset M$ are nothing but
the definition of the canonical Courant
$C^{\infty}(U_{\Sigma})$-algebroid $\CC_{0}$ of (2.4.21); therefore
$\CO_{J^{\infty}(T^{*}M_{\Sigma/\Sigma''})}(U)$ is nothing but
$\CF^{*}(\CC_{0})$, where $\CF$ is  equivalence of categories
(2.5.1).

The vertex Poisson algebra structure of Proposition 2.7.1 is not
quite what we need. Being $\CT_{\Sigma}$-equivariant, it is subject
to the $\xi$-twist, see (2.3.3), for any $\xi\in
H^{0}(\Sigma,\CT_{\Sigma})$.

\bigskip

{\bf 2.7.2. Definition.} {\it Let
$\widehat{\CO}_{J^{\infty}(T^{*}M_{\Sigma/\Sigma''})}$ denote the
sheaf $\CO_{J^{\infty}(T^{*}M_{\Sigma/\Sigma''})}$ with the vertex
Poisson algebra structure defined in Proposition 2.7.1 and let
$\CO_{J^{\infty}(T^{*}M_{\Sigma/\Sigma''})}$ denote the latter's
$\partial_{\sigma}^{h}$-twist, see (2.7.1).}

\bigskip

Note that in the case of
$\CO_{J^{\infty}(T^{*}M_{\Sigma/\Sigma''})}$, the derivation $T$
becomes
$$
T=\rho(\partial_{\sigma}). \eqno{(2.7.8)}
$$
In particular, (2.7.6) is changed as follows
$$
\aligned
 _{(0)}: \left(\CT_{M_{\Sigma}/\Sigma}\oplus\Omega_{M_{\Sigma}/\Sigma}\right)\otimes
\left(\CT_{M_{\Sigma}/\Sigma}\oplus\Omega_{M_{\Sigma}/\Sigma}\right)\rightarrow
\CT_{M_{\Sigma}/\Sigma}\oplus\Omega_{M_{\Sigma}/\Sigma},\\
\left(\xi+\omega\right)_{(0)}\left(\xi'+\omega'\right)=[\xi,\xi']+\text{Lie}_{\xi}\omega'
-\text{Lie}_{\xi'}\omega+\rho(\partial_{\sigma})\left(\iota_{\xi'}\omega\right),
\endaligned
\eqno{(2.7.9)}
$$
and the operations on $\CO_{J^{\infty}(T^{*}M_{\Sigma/\Sigma''})}$
are no longer linear over $\CO_{\Sigma}$, only over
$\CO_{\Sigma''}$.

Let us now relate $\CO_{J^{\infty}(T^{*}M_{\Sigma/\Sigma''})}$ to
the canonical Lie algebra sheaf $\CH^{can}$ defined in (1.5.17).

Lemma 2.2 associates with
$\CO_{J^{\infty}(T^{*}M_{\Sigma/\Sigma''})}$ the sheaf of Lie
algebras
$$
\text{Lie}(\CO_{J^{\infty}(T^{*}M_{\Sigma/\Sigma''})})=
\CO_{J^{\infty}(T^{*}M_{\Sigma/\Sigma''})}/\rho(\partial_{\sigma})\CO_{J^{\infty}(T^{*}M_{\Sigma/\Sigma''})}.
$$

\bigskip

{\bf 2.7.3. Proposition.} {\it The Lie algebra  sheaves $\CH^{can}$
and} $\text{Lie}(\CO_{J^{\infty}(T^{*}M_{\Sigma/\Sigma''})})$ {\it
are canonically isomorphic.}

\bigskip

{\it Proof.} The sheaf isomorphism
$$
\CO_{J^{\infty}(T^{*}M_{\Sigma/\Sigma''})}\rightarrow
\Omega^{1,0}_{J^{\infty}(T^{*}M_{\Sigma/\Sigma''})},\; F\mapsto
Fd\sigma
$$
descends to
$$
\CO_{J^{\infty}(T^{*}M_{\Sigma/\Sigma''})}/\rho(\partial_{\sigma})\CO_{J^{\infty}(T^{*}M_{\Sigma/\Sigma''})}
\rightarrow \Omega^{1,0}_{J^{\infty}(T^{*}M_{\Sigma/\Sigma''})}/
d_{\rho/\Sigma} \Omega^{0,0}_{J^{\infty}(T^{*}M_{\Sigma/\Sigma''})}.
\eqno{(2.7.10)}
$$
(It is at this point that we needed the
$\partial_{\sigma}^{h}$-twist, see (2.7.8).) Map (2.7.10) respects
all defining relations (1.5.18-19a,b): (1.5.18) is (part of)
(2.7.5,6), (1.5.19a) is sect.2.1, Axiom III, and (1.5.19b) is
sect.2.1, Axiom II.3 (another point where the
$\partial_{\sigma}^{h}$-twist is necessary).

\bigskip

{\bf 2.7.4. Terminology.} We have obtained two families of sheaves
of vertex Poisson algebras. First, those provided by the combination
of Propositions 2.6.1a) and 2.7.3. They can be realized as either
$\widehat{\CO}_{J^{\infty}(T^{*}M_{\Sigma/\Sigma''})}\plus H$, where
$H\in H^{0}(M,\Omega^{3,cl}_{M})$ represents a 3-dimensional
cohomology class, or
$\widehat{\CO}_{J^{\infty}(T^{*}M_{\Sigma/\Sigma''})}\plus
\left(\{\alpha_{ij}\}\right)$, where $\left(\{\alpha_{ij}\}\right)$
is a cocycle representing an element of
$H^{1}(M,\Omega^{2,cl}_{M})$.

Second, their $\partial_{\sigma}^{h}$-twisted versions, to be
denoted by $\CO_{J^{\infty}(T^{*}M_{\Sigma/\Sigma''})}\plus H$ and
$\CO_{J^{\infty}(T^{*}M_{\Sigma/\Sigma''})}\plus
\left(\{\alpha_{ij}\}\right)$. As Proposition 2.7.3 indicates, it is
the latter that will be of importance. Note, however, that these
choices have arisen only because we have included functions of
$\tau$ and $\sigma$. In fact, both
$\widehat{\CO}_{J^{\infty}(T^{*}M_{\Sigma/\Sigma''})}\plus
\left(\{\alpha_{ij}\}\right)$ and
$\CO_{J^{\infty}(T^{*}M_{\Sigma/\Sigma''})}\plus
\left(\{\alpha_{ij}\}\right)$ induce the same
   vertex Poisson algebra structure on the fiber at any point
  $(\sigma,\tau)\in\Sigma$. For this reason sheaves such as
  $\CO_{J^{\infty}(T^{*}M_{\Sigma/\Sigma''})}\plus
\left(\{\alpha_{ij}\}\right)$, where $\left(\{\alpha_{ij}\}\right)$
will also be referred to as sheaves of SVDOs.

\bigskip

{\bf 2.8.} {\bf The Lagrangian interpretation}

\bigskip
 Let us place
ourselves in the situation of 1.6.10 and assume that the Lagrangian
$L\in
H^{0}(J^{\infty}(M_{\Sigma}),\Omega^{2,0}_{J^{\infty}(M_{\Sigma})})$
is of order 1, globally defined, and convex. Recall that associated
to $L$ there are the diffiety $Sol_{L}$ carrying the canonical
variational 2-form $\omega_{L}$, the Lie algebra of local
functionals $\CH_{Sol_{L}}$ and the subalgebra of integrals of
motion $\tilde{\CI}_{L}\subset \CH_{Sol_{L}}$.

Being a diffiety, $Sol_{L}$ is a fiber bundle with flat connection
over $\Sigma$. One can also consider the version of this connection
realtive to $\Sigma\rightarrow\Sigma''$. There arises then the
relative version $\CH_{Sol_{L}/\Sigma''}$.

If $L$ is as above required, then the corresponding Legendre
transform $d_{TM}\tilde{L}:TM\rightarrow T^{*}M$ identifies
$$
d_{TM}L: Sol_{L}\rightarrow J^{\infty}(T^{*}M_{\Sigma/\Sigma''})
\eqno{(2.8.1)}
$$
as manifolds with connection over $\Sigma\rightarrow\Sigma''$, cf.
(1.5.20,1.6.10), and, by virtue of Lemma 1.6.10.1, gives rise to
maps of Lie algebra sheaves
$$
\CH^{\omega_{L}}_{Sol/\Sigma''} \iso \CH^{can} \eqno{(2.8.2)}.
$$
Combined with Proposition 2.7.3 this gives
$$
\CH^{\omega_{L}}_{Sol/\Sigma''} \iso
\text{Lie}(\CO_{J^{\infty}(T^{*}M_{\Sigma/\Sigma''})}).
\eqno{(2.8.3)}
$$

In this sense the universal sheaf of SVDO's
$\CO_{J^{\infty}(T^{*}M_{\Sigma/\Sigma''})}$ governs the theory
associated to $L$.

In order to interpret  similarly all the other, twisted, sheaves
 of SVDO's provided by Proposition 2.6.1a),  one needs to consider Lagrangians
(1.6.1a,b) that do not  glue in a  global section of
$\Omega^{2,0}_{J^{\infty}(M_{\Sigma})}$. One possibility to
construct such a Lagrangian is to add what a physicist might call a
{\it Wess-Zumino term} or an $H$-{\it flux}, cf. [GHR, W1].

Fix a global closed 3-form $H$ on $M$ and let $\{U_{i}\}$ be an open
covering of $M$ fine enough that there exist a collection of
$2$-forms
$$\{\beta^{(i)}\in\Omega^{2}_{M}(U_{i})\text{ s.t. }
d\beta_{i}=H\text{ on }U_{i}\}. \eqno{(2.8.4)}
$$
Define
$$
L^{H}=\{L^{(i)}=L+\beta^{(i)}(\rho(\partial_{\tau}),\rho(\partial_{\sigma}))\}.\eqno{(2.8.5)}
$$
It follows from (2.8.4) that on double intersection
$\beta^{(i)}-\beta^{(j)}$ are closed and, provided $\{U_{i}\}$ is
fine enough, are exact, i.e., there is a collection of 1-forms,
$\{\alpha^{(ij)}\}$ such that
$\beta^{(j)}-\beta^{(i)}=d\alpha^{(ij)}$. Then a quick computation
shows that
$$
L^{(j)}-L^{(i)}=d_{\rho}\left(\left(
\iota_{\rho(\partial_{\tau})}\alpha^{(ij)}\right)d\tau+
\left(\iota_{\rho(\partial_{\sigma})}\alpha^{(ij)}\right)d\sigma\right).
$$
Therefore, collection (2.8.5) is a new Lagrangian in the sense of
(1.6.1a,b).

The  $L^{H}$ is a collection of locally defined Lagrangians, which
are still order 1 and convex, hence $Sol_{L^{H}}$ can still be
identified with the universal
$J^{\infty}(T^{*}M_{\Sigma/\Sigma''})$. One way to define such
identification is to use $L$, as in (2.8.1):
$$
d_{TM}L:\;(Sol_{L^{H}})\iso
J^{\infty}(T^{*}M_{\Sigma/\Sigma''}),\eqno{(2.8.6)}
$$
but the  obvious counterpart of (2.8.3) fails in this case. Instead,
(2.8.6) delivers an isomorphism of the twisted sheaf
$$
 \text{Lie}(\CO_{J^{\infty}(T^{*}M_{\Sigma/\Sigma''})}\plus H)\iso
\CH^{\omega_{L^{H}}}_{\Sigma/\Sigma''}.\eqno{(2.8.7)}
$$
(Remember that the $H$-twist does nothing to the  sheaf structure
and only affects the $_{(n)}$-multiplications.) This attaches the
twisted sheaf $\CO_{J^{\infty}(T^{*}M_{\Sigma/\Sigma''})}\plus H$ to
the Lagrangian $L^{H}$.

To see how the twist comes about note that the Legendre transform
$d_{TM}L$ used in (2.8.6) does not respect the canonical variational
2-form $\omega_{L^{H}}$. This can be straightened out locally.
According to (2.8.1), one way to proceed is to choose, over $U_{i}$,
the mapping to be $d_{TM}L^{(i)}$. Since $L^{(i)}=L+\beta^{(i)}$,
$$
d_{TM}L^{(i)}(\xi)=d_{TM}L(\xi)+\frac{1}{2}\iota_{\xi}\beta^{(i)},
\eqno{(2.8.8)}
$$
as follows, e.g., from local formulas (1.6.11). But mappings (2.8.8)
are incompatible on double intersections $U_{i}\cap U_{j}$, the
obstruction being the \v Cech cocycle
$$
d_{\check{C}}\{\beta^{(i)}\}=\{\beta^{(j)}-\beta^{(i)}\}\in
Z^{1}_{\text{\v Cech}}(M,\Omega^{2,cl}_{M}). \eqno{(2.8.9)}
$$
In order to restore the compatibility, let us introduce the twisted
sheaf
$\CO_{J^{\infty}(T^{*}M_{\Sigma/\Sigma''})}\plus\left(d_{\check{C}}\{\beta^{(i)}\}\right)$
obtained by twisting the gluing functions of
$\CO_{J^{\infty}(T^{*}M_{\Sigma/\Sigma''})}$ over $U_{i}\cap U_{j}$
by the 2-form $\beta^{(j)}-\beta^{(i)}$, as we did in (2.6.9). Then
the collection of mappings
$$
\{(d_{TM}L^{(i)})^{*}:\;
\;\CO_{J^{\infty}(T^{*}M_{\Sigma/\Sigma''})}(U_{i})\iso
\CO_{Sol_{L^{H}}}(U_{i})\eqno{(2.8.10a)}
$$
delivers a map of the twisted sheaf
$$
\left(\CO_{J^{\infty}(T^{*}M_{\Sigma/\Sigma''})}\plus\left(d_{\check{C}}\{\beta^{(i)}\}\right)\right)
\rightarrow \CO_{Sol_{L^{H}}},\eqno{(2.8.10b)}
$$
so that the arising
$$
 \text{Lie}\left(\CO_{J^{\infty}(T^{*}M_{\Sigma/\Sigma''})}\plus
\left(d_{\check{C}}\{\beta^{(i)}\}\right)\right)\iso
\CH^{\omega_{L^{H}}}_{\Sigma/\Sigma''}.\eqno{(2.8.11)}
$$
is a Lie algebra sheaf isomorphism.
 The twisted sheaf of SVDOs
$\CO_{J^{\infty}(T^{*}M_{\Sigma/\Sigma''})}\plus
\left(d_{\check{C}}\{\beta^{(i)}\}\right)$ appearing in (2.8.11) is
defined in sect. 2.5.2, formula (2.5.8). As is also explained in
some detail in 2.5.2, this sheaf is the same as
$\CO_{J^{\infty}(T^{*}M_{\Sigma/\Sigma''})}\plus H$; hence (2.8.11)
is equivalent to (2.8.7).

Incidentally, the classification of automorphisms of SVDO's,
Proposition 2.6.1b) is also accurately reflected in the Lagrangian
approach. Given a globally defined Lagrangian and a closed 2-form
$\beta$, {\it a B-field}, let
$L^{\beta}=L+\beta(\rho(\partial_{\tau}),\rho(\partial_{\sigma}))$,
cf. (2.7.5). This does nothing to the corresponding equations of
motion; hence $Sol_{L}=Sol_{L^{\beta}}$, literally, but there arise
two competing Legendre transforms, $d_{TM}L$ and $d_{TM}L^{\beta}$.
A moment's thought shows that the latter is the composition of the
former with the {\it B-field transform},
$\xi\mapsto\xi+\iota_{\xi}\beta$, and this provides the Lagrangian
realization of the automorphism of the SVDO
$\CO_{J^{\infty}(T^{*}M_{\Sigma/\Sigma''})}$ associated to $\beta$
in Proposition 2.6.1b.

The subalgebras of integrals of motion $\hat{\CI}_{L}\hookrightarrow
\Gamma(Sol_{L},\CH^{\omega_{L}}_{\Sigma/\Sigma''})$, arising by
virtue of Lemma 1.6.9, also tend to come from vertex Poisson
subalgebras of $\CO_{J^{\infty}(T^{*}M_{\Sigma/\Sigma''})}$. For
example, the three Virasoro algebras, left, right, and
``half-twisted'', see (1.6.19,20), are the $\text{Lie}$-functor
applied to the three subalgebras of
$\Gamma(M,\CO_{J^{\infty}(T^{*}M_{\Sigma/\Sigma''})})$ generated by
$$
\aligned &\frac{1}{4}g^{ij}x_{i}x_{j}
+\frac{1}{4}g_{ij}\partial_{\sigma}x^{i}\partial_{\sigma}x^{j}
-\frac{1}{2}x_{j}\partial_{\sigma}x^{j},\\
 -&\frac{1}{4}g^{ij}x_{i}x_{j}
-\frac{1}{4}g_{ij}\partial_{\sigma}x^{i}\partial_{\sigma}x^{j}
-\frac{1}{2}x_{j}\partial_{\sigma}x^{j},\\
-&x_{j}\partial_{\sigma}x^{j},
\endaligned
\eqno{(2.8.12)}
$$
respectively. The global nature of these local formulas formulas was
unraveled  in 1.6.11.2.

\bigskip\bigskip

{\bf 2.9. An example: WZW model.}

\bigskip

Let us see how all of this plays out in the case where the target
manifold is a real Lie group $G$, either compact and simple or
$\text{GL}(n,\BR)$.

{\bf 2.9.1.} {\it Classification.} Let $\fg=\text{Lie}G$ be the
corresponding Lie algebra. Fix an invariant bilinear form
 $g\in S^{2}(\fg^{*})^{\fg}$ and an invariant trilinear form
 $$
 H(x,y,z)=g([x,y],z).\eqno{(2.9.1)}
 $$
 The left translates of these generate the invariant metric and 3-form
 (resp.)
 on $G$, which we will take the liberty of denoting by the same
 letters
 $$
 g\in H^{0}(G,\CT^{\otimes 2}_{G}),\; H\in H^{0}(G,\Omega^{3}_{G}).
 \eqno{(2.9.2)}
 $$
 Note that the latter is closed:
 $$
 H\in H^{0}(G,\Omega^{3,cl}_{G}).
 \eqno{(2.9.3)}
 $$

It is well known that
$$
H^{3}(G,\BR)=\BR\cdot(\text{class of }H).\eqno{(2.9.4)}
$$
Therefore, Proposition 2.6.1a) implies that the set of isomorphism
classes of SVDO's on $G_{\Sigma}$ form a 1-parameter family:
$$
SD_{G,k}\buildrel\text{def}\over =
\CO_{J^{\infty}(G_{\Sigma/\Sigma''})}\plus
\frac{-k}{2}H.\eqno{(2.9.5)}
$$
As was explained in sect. 2.6.2, the structure of $SD_{G,k}$ is
determined by the following:

there is a canonical splitting
$$
\left(SD_{G,k}\right)_{1}=\CT_{G_{\Sigma}/\Sigma}\oplus\Omega_{G_{\Sigma}/\Sigma},
\eqno{(2.9.6)}
$$
and the vertex Poisson algebra structure makes
$\left(D^{poiss}_{G,k}\right)_{1}$ into the Courant algebroid that
satisfies
$$
(2.3.14,15)\text{ hold true, and
}_{(0)_{s}}=-\frac{k}{2}H,\eqno{(2.9.7)}
$$
cf. the definition, (2.4.21), of the canonical Courant algebroid
$\CC_{0}$.

Induced by the action on the left and on the right, there are the
corresponding Lie algebra $\fg=\text{Lie}G$  embeddings in the space
of global vector fields
$$
j^{0}_{l}:\fg\hookrightarrow
H^{0}(G_{\Sigma},\CT_{G_{\Sigma}/\Sigma}),\;
j^{0}_{r}:\fg\hookrightarrow
H^{0}(G_{\Sigma},\CT_{G_{\Sigma}/\Sigma}) \text{ s.t.}
[j^{0}_{l}(\fg),j^{0}_{r}(\fg)]=0.\eqno{(2.9.8)}
$$
These embeddings respect the SVDO structure on
$\left(SD_{G,0}\right)_{1}$ in that
$$
j^{0}_{l}\left([x,y]\right)=\left(j^{0}_{l}(x)_{(0)}j^{0}_{l}(y)\right),\;
\left(j^{0}_{l}(x)_{(n)}j^{0}_{l}(y)\right)\text{ if
}n>0,\eqno{(2.9.9a)}
$$
$$
j^{0}_{r}\left([x,y]\right)=\left(j^{0}_{r}(x)_{(0)}j^{0}_{r}(y)\right),\;
\left(j^{0}_{r}(x)_{(n)}j^{0}_{r}(y)\right)\text{ if
}n>0,\eqno{(2.9.9b)}
$$
and
$$
\left(j^{0}_{l}(x)\right)_{(n)}\left(j^{0}_{l}(x)\right)=0\text{ if
}n\geq 0,\eqno{(2.9.9c)}
$$
as follows from either (2.7.6) or  (2.4.6).

 Technically, (2.9.9a--c)
mean the following. Associated to $\fg$ there is a $\BZ_{+}$-graded
vertex Poisson algebra, $V(\fg)_{k}$, see e.g. [FB-Z]. It is the
universal vertex Poisson algebra generated by
$$
\left(V(\fg)_{k}\right)_{0}=\BR,\; \left(V(\fg)_{k}\right)_{1}=\fg,
\eqno{(2.9.10)}
$$
subject to relations
$$
x_{(n)}y=\left\{\aligned kg(x,y)&\text{ if }n=1\\
[x,y]&\text{ if }n=0\\
0&\text{ if }n> 1.
\endaligned\right.
\eqno{(2.9.11)}
$$
By definition, (2.9.9a--c) imply that maps (2.9.8) can be extended
to vertex Poisson algebra maps
$$
 j^{0}_{l}:V(\fg)_{0}\hookrightarrow H^{0}(G_{\Sigma},SD_{G,0}),\;
j^{0}_{r}:V(\fg)_{0}\hookrightarrow
H^{0}(G_{\Sigma},SD_{G,0})\eqno{(2.9.12)}
$$
such that
$$
\left(j^{0}_{l}(V(\fg)_{0})\right)_{(n)}\left(j^{0}_{r}(V(\fg)_{0})\right)=0
\text{ if }n\geq 0. \eqno{(2.9.13)}
$$
In order to carry this  over to $k\neq 0$, the maps $j^{0}_{l/r}$
must be deformed. Let
$$
j^{k}_{l}:\fg\rightarrow \left(D^{poiss}_{G,k}\right)_{1},\;
j^{k}_{l}(x)=j^{0}_{l}(x)+\frac{k}{2}g\left(j^{0}_{l}(x),.\right),\eqno{(2.9.14)}
$$
$$
j^{k}_{r}:\fg\rightarrow \left(D^{poiss}_{G,k}\right)_{1},\;
j^{k}_{r}(x)=j^{0}_{r}(x)-\frac{k}{2}g\left(j^{0}_{r}(x),.\right),\eqno{(2.9.15)}
$$

\bigskip

{\bf 2.9.1.1. Theorem.} [FP,F,AG,GMS2] {\it Maps (2.9.14,15) extend
to vertex Poisson algebra embeddings
$$
V(\fg)_{k}\buildrel j^{k}_{l}\over\hookrightarrow
H^{0}(G_{\Sigma},SD_{G,k})\buildrel j^{k}_{r}\over\hookleftarrow
V(\fg)_{-k}\eqno{(2.9.16)}
$$
such that
$$
\left(j^{k}_{l}(V(\fg)_{k})\right)_{(n)}\left(j^{k}_{r}(V(\fg)_{-k})\right)=0
\text{ if }n\geq 0. \eqno{(2.9.17)}
$$}

\bigskip

{\bf 2.9.1.2. Remark.} This appealing result has a long and somewhat
unhappy history. A version of it first appeared in [FP] (in a more
complicated, quantum, situation) but apparently had been known even
earlier to E.Frenkel, [F] -- all of this before the introduction of
sheaves of vertex algebras -- and then was thoroughly forgotten.
Arkhipov and Gaitsgory [AG] gave a proof in the language of chiral
algebras. Our presentation is close to [GMS2].

\bigskip

The algebra $V(\fg)_{k}$ has a well-known family of modules,
$V_{\lambda,k}$, induced from $V_{\lambda}$, the simple finite
dimensional $\fg$-module with highest weight $\lambda$, see e.g.
[FBZ]. According  to Theorem 2.9.1.1, $H^{0}(G_{\Sigma},SD_{G,k})$
is a $V(\fg)_{k}\otimes V(\fg)_{-k}$-module, see sect. 2.3 for the
definition of the tensor product of vertex Poisson algebras.

\bigskip

{\bf 2.9.1.3. Proposition.} {\it If $k\neq 0$, then there is an
isomorphism of $V(\fg)_{k}\otimes V(\fg)_{-k}$-modules
$$
H^{0}(G_{\Sigma},SD_{G,k})\iso
C^{\infty}(\Sigma)\hat{\otimes}\left(\oplus_{\lambda}V_{\lambda,k}\otimes
V_{\lambda^{*},-k}\right),\eqno{(2.9.18)}
$$
where $\lambda^{*}$ stands for the highest weight of the
$\fg$-module dual to $V_{\lambda}$.}

\bigskip

\begin{sloppypar}
{\it Sketch of Proof. } The validity of  decomposition (2.9.18) for
the subspace $H^{0}(G_{\Sigma},(SD_{G,k})_{0})$ is the content of
the Peter-Weyl theorem. It is not hard to deduce from (2.9.14,15)
that any $\BR$-basis, $\CB$, of $j^{k}_{l}(\fg)\oplus
j^{k}_{r}(\fg)$ is a basis of $H^{0}(G_{\Sigma},(SD_{G,k})_{1})$
over functions. Hence, the entire
$H^{0}(G_{\Sigma},D^{poiss}_{G,k})$ is the space of differential
polynomials in $\CB$ over functions. Decomposition (2.9.18) follows
at once from the induced nature of modules $V_{\lambda,k}$.
\end{sloppypar}

\bigskip

{\bf 2.9.1.4. Remark.} A proof -- in the quantum case -- of (2.9.18)
for a generic $k$ first appeared in [FS]. Our proof goes through in
the quantum case  for all but one value of $k$.

\bigskip

Decomposition (2.9.21) is tantalizingly similar to space of states
of WZW model to which  $SD_{G,k}$ is intimately related.

\bigskip

{\bf 2.9.2. WZW.} Consider the standard $\sigma$-model Lagrangian
with target $G$:
$$
L_{\kappa}=\frac{\kappa}{2}g\left(\left(\partial_{\tau}-\partial_{\sigma}\right)x,
\left(\partial_{\tau}+\partial_{\sigma}\right)x\right)d\tau\wedge
d\sigma, \eqno{(2.9.19)}
$$
cf. (1.6.13), where $g(.,.)$ is the invariant metric (2.9.2) and
$\kappa$ is an arbitrary constant.

Next use the 3-form $H$ of (2.9.2) to obtain $L_{\kappa}^{-k/2H}$ as
explained in (2.8.4-5). The WZW Lagrangian [W1] is
$$
L_{WZW}=L_{k/2}^{-\frac{k}{2}H}.\eqno{(2.9.20)}
$$
As follows from (2.8.7) and normalization (2.9.5), the sheaf
$SD_{G,k}$ governs the theory associated to $L_{\kappa}^{-k/2H}$ for
any $\kappa$. It is clear why the $H$-twist of (2.9.19) is needed --
the pleasing decomposition (2.9.18) is valid only if $k\neq 0$.

Let us now explain the choice of $\kappa$ made in (2.9.20). Recall
that Lagrangian (2.9.19) is conformally invartiant, i.e., the
corresponding algebra of integrals of motion contains two Virasoro
subalgebras $\Vir^{\pm}$, see (1.6.17). It is easy to see that the
twisted version, $L_{\kappa}^{-k/2H}$, is also, and $\Vir^{\pm}$ are
still the corresponding integrals of motion. By virtue of (2.8.7)
the Legendre transform delivers the embeddings
$$
\Vir^{\pm}\hookrightarrow\Gamma(G,\text{Lie}\left(SD_{G,k}\right)).\eqno{(2.9.21)}
$$
On the other hand,  each $V(\fg)_{k}$ carries its own Virasoro
element -- well-known fact. By virtue of Theorem 2.9.1.1, there
arise then two more Virasoro subalgebras
$$
\Vir^{l}\hookrightarrow\Gamma(G,\text{Lie}\left(SD_{G,k}\right))\hookleftarrow\Vir^{r}.\eqno{(2.9.22)}
$$
\bigskip

{\bf 2.9.2.1. Lemma.} {\it Upon taking the images of (2.9.21-22)
$$
\Vir^{+}=\Vir^{l},\; \Vir^{-}=\Vir^{r}\eqno{(2.9.23)}
$$
if and only if $\kappa=k/2$.}

\bigskip

This allows to compute the left/right moving subalgebra, see
Definition 1.6.11.1.

\bigskip

{\bf 2.9.2.2. Corollary.} {\it The right moving subalgebra of WZW is
$\text{Lie}(C^{\infty}(\Sigma)\otimes V(\fg)_{k})$ and the left
moving is $\text{Lie}(C^{\infty}(\Sigma)\otimes V(\fg)_{-k})$.}

\bigskip

The Lie-functor appearing in 2.9.2.1-2 only obscures the matter, of
course. Armed with the notion of a vertex Poisson algebra we can
easily refine both Definition 1.6.11.1 and  2.9.2.1-2. The Lie
algebra $\Vir$ itself is the Lie-functor applied to a certain vertex
Poisson algebra, $\text{Vir}$. Embeddings (2.9.21-22) are engendered
by vertex Poisson algebra embeddings of 4 copies of $\text{Vir}$:
$$
\text{Vir}^{\pm}\hookrightarrow \Gamma(G,SD_{G,k}),\eqno{(2.9.24)}
$$
$$
\text{Vir}^{l}\hookrightarrow \Gamma(G,SD_{G,k}) \hookleftarrow
\text{Vir}^{r}.\eqno{(2.9.25)}
$$
Lemma 2.9.2.1 can be refined as follows: upon taking the images of
(2.9.24,25)
$$
\text{Vir}^{+}=\text{Vir}^{l},\; \text{Vir}^{-}=\text{Vir}^{r}
\text{ iff } \kappa=\frac{k}{2}. \eqno{(2.9.26)}
$$
Definition 1.6.11.1 can be similarly refined:

\bigskip

{\bf 2.9.2.3. Definition.} Let the left/right moving subalgebras of
$SD_{G,k}$ be
$$
SD^{+}_{G,k}=\{v\in SD_{G,k}\text{ s.t. }
v_{(n)}\text{Vir}^{-}=0\text{ if }n\geq 0\}\eqno{(2.9.27)}
$$
$$
SD^{-}_{G,k}=\{v\in SD_{G,k}\text{ s.t. }
v_{(n)}\text{Vir}^{+}=0\text{ if }n\geq 0\}\eqno{(2.9.28)}
$$
\bigskip

The refined form of Corollary 2.9.2.2 is this:
$$
SD^{+}_{G,k}=C^{\infty}(\Sigma)\otimes V(\fg)_{k},\;
SD^{-}_{G,k}=C^{\infty}(\Sigma)\otimes V(\fg)_{-k}. \eqno{(2.9.29)}
$$
\bigskip

{\bf 2.9.2.4. Proofs.}  We will prove (2.9.26) and (2.9.29) from
which Lemma 2.9.2.1 and Corollary 2.9.2.2 follow instantaneously.
Proving (2.9.26) amounts to painstakingly translating from sect.
2.9.1 to sect. 2.9.2, the Legendre transform being the main tool.

To facilitate  bookkeeping, we will assume that
$G=\text{GL}(n,\BR)$; an extension via a faithful representation to
compact Lie groups is immediate. Let then $x^{ij}$ be coordinates,
$\partial_{ij}=\partial/\partial x^{ij}$, and $\{E_{ij}\}$ the
standard basis of $gl(n,\BR)$.

The invariant metric is
$$
g=x_{t\alpha}dx^{\alpha j}x_{j\beta}dx^{\beta t},\eqno{(2.9.30)}
$$
where $x_{t\alpha}$ are defined so that $x_{t\alpha} x^{\alpha
j}=\delta_{t}^{j}$, and the summation w.r.t. repeated indices is
always assumed.

Embeddings (2.9.8.) take on the form
$$
j_{l}^{0}(E_{ij})=x^{\alpha i}\partial_{\alpha j},\eqno{(2.9.31)}
$$
$$
j_{r}^{0}(E_{ij})=-x^{j\alpha }\partial_{i\alpha}.\eqno{(2.9.32)}
$$
By virtue of (2.9.30), definitions (2.9.14,15) read
$$
j_{l}^{k}(E_{ij})=x^{\alpha i}\partial_{\alpha
j}+\frac{k}{2}x_{j\gamma}\partial_{\sigma}x^{\gamma
i},\eqno{(2.9.33)}
$$
$$
j_{r}^{k}(E_{ij})=-x^{j\alpha
}\partial_{i\alpha}+\frac{k}{2}x_{\gamma i}\partial_{\sigma}x^{j
\gamma}.\eqno{(2.9.34)}
$$

Finally, the elements that generate the two corresponding Virasoro
vertex Poisson algebras inside $SD_{G,k}$, cf. (2.9.25), are
$$
\text{Vir}^{l}=<\frac{1}{k}j_{l}^{k}(E_{ij})j_{l}^{k}(E_{ji})>,\;
\text{Vir}^{r}=<\frac{1}{k}j_{r}^{k}(E_{ij})j_{r}^{k}(E_{ji})>.
\eqno{(2.9.35)}
$$
To recapitulate  all of this in terms intrinsic to the Lagrangian
$L_{\kappa}^{-k/2H}$, one needs to use the twisted version of the
Legendre transform, see (2.8.6), i.e., apply (1.6.10-11) not to
$L_{\kappa}^{-k/2H}$ but to $L_{\kappa}^{0}$. This amounts to
letting
$$
\partial_{ij}=\frac{\partial
L_{\kappa}^{0}}{\partial(\partial_{\tau}x^{ij})};
$$
thus
$$
\partial_{ij}=\kappa x_{\alpha
i}\partial_{\tau}x^{\beta\alpha}x_{j\beta}.\eqno{(2.9.36)}
$$
Plugging this in (2.9.33-34) gives
$$
j_{l}^{k}(E_{ij})=\left(\kappa\partial_{\tau}x^{\alpha
i}+\frac{k}{2}\partial_{\sigma}x^{\alpha
i}\right)x_{j\alpha},\eqno{(2.9.37)}
$$
$$
j_{r}^{k}(E_{ij})=\left(-\kappa\partial_{\tau}x^{j\alpha}
+\frac{k}{2}\partial_{\sigma}x^{j\alpha}\right)x_{\alpha
i}.\eqno{(2.9.38)}
$$
It is pleasing to notice that precisely when $\kappa=k/2$, the
latter formulas become the WZW currents, see [W1], (15) or [GW],
(2.3),
$$
j_{l}^{k}(E_{ij})=k\partial_{+}x^{\alpha
i}x_{j\alpha},\eqno{(2.9.39)}
$$
$$
j_{r}^{k}(E_{ij})=k\partial_{-}x^{j\alpha}x_{\alpha
i},\eqno{(2.9.40)}
$$
where $\partial_{\pm}=(\partial_{\sigma}\pm\partial_{\tau})/2$.

Now to the Virasoro subalgebras. Plugging (2.9.37,38) in (2.9.35)
one finds similarly that precisely when $\kappa=k/2$ the
corresponding Virasoro elements are
$$
\text{Vir}^{l}=<k b(\partial_{+}x,\partial_{+}x)>,\eqno{(2.9.41)}
$$
$$
\text{Vir}^{r}=<k b(\partial_{-}x,\partial_{-}x)>,\eqno{(2.9.42)}
$$
i.e., defined by the familiar, see  (1.6.17), formulas for
$\text{Vir}^{\pm}$. This concludes our proof of (2.9.26).

\bigskip

Now to (2.9.29).  Having it our disposal (2.9.26), we infer from
Theorem 2.9.1.1  that
$$
C^{\infty}(\Sigma)\otimes V(\fg)_{k}\subset SD^{+}_{G,k},
C^{\infty}(\Sigma)\otimes V(\fg)_{-k}\subset
SD^{-}_{G,k}.\eqno{(2.9.43)}
$$
To prove the reverse inclusions, let
$$
L^{l}=\frac{1}{k}j_{l}^{k}(E_{ij})j_{l}^{k}(E_{ji}),
L^{r}=\frac{1}{k}j_{r}^{k}(E_{ij})j_{r}^{k}(E_{ji}).
$$
It follows easily from the definition of the modules $V_{\lambda,k}$
that
$$
\text{Ker}L_{(0)}^{l}=V_{0,-k}\buildrel\text{def}\over =
V(\fg)_{-k},  \text{Ker}L_{(0)}^{r}=V_{0,k}\buildrel\text{def}\over
= V(\fg)_{k}.\eqno{(2.9.44)}
$$
By definition then
$$
C^{\infty}(\Sigma)\otimes V(\fg)_{k}\supset SD^{+}_{G,k},
C^{\infty}(\Sigma)\otimes V(\fg)_{-k}\supset
SD^{-}_{G,k},\eqno{(2.9.45)}
$$
which concludes the proof of (2.9.29).

\bigskip\bigskip

\begin{center}
{\bf 3. Supersymmetric Analogues}
\end{center}

\bigskip\bigskip

{\bf 3.1. Bits of supergeometry.}

\bigskip

All of the geometric background of sect. 1 allows for more or less
straightforward super-generalization. We will explain this very
briefly, and in less generality, because our exposition will be more
example-oriented. Such sources  as [DM,L,M1] provide an introduction
to supermathematics.

\bigskip

{\bf 3.1.1. Super world-sheet.} The world-sheet  is now  a
$2|2$-dimensional real $C^{\infty}$-manifold either with a fixed
coordinate system
$$
(u,v,\theta^{+},\theta^{-}): \hat{\Sigma}\hookrightarrow\BR^{2|2}
\eqno{(3.1.1a)}
$$
or a fixed \'etale coordinate system
$$
(u,v,\theta^{+},\theta^{-}): \BR^{2|2}\rightarrow\hat{\Sigma},
\eqno{(3.1.1b)}
$$
where $(u,v)$ are even and $\theta^{\pm}$ are odd.

We have the underlying even manifold
$$
\Sigma=\{\theta^{+}=\theta^{-}=0\}\hookrightarrow\hat{\Sigma}\eqno{(3.1.2)}
$$
and the bundle
$$
\hat{\Sigma}\buildrel (u,v)\over\rightarrow\Sigma.\eqno{(3.1.3)}
$$
The time-fibration will be defined to be the composition
$$
\hat{\Sigma}\buildrel
(u,v)\over\rightarrow\Sigma\buildrel\tau\over\rightarrow\Sigma''\subset\BR
\eqno{(3.1.4)}
$$
for some fibration $\tau$, where $\Sigma$ is an even manifold
underlying $\hat{\Sigma}$.

The Lie algebra of vector fields on $\hat{\Sigma}$ contains two
remarkable elements
$$
D_{+}=\frac{\partial}{\partial\theta^{+}}-\theta^{+}\frac{\partial}{\partial
u},\;
D_{-}=\frac{\partial}{\partial\theta^{-}}-\theta^{-}\frac{\partial}{\partial
v}.\eqno{(3.1.5)}
$$
The following relations hold true
$$
\aligned &[D_{+},D_{+}]=-2\frac{\partial}{\partial u},
[D_{-},D_{-}]=-2\frac{\partial}{\partial v}, [D_{+},D_{-}]=0\\
&[\frac{\partial}{\partial v},D_{\pm}]= [\frac{\partial}{\partial
u},D_{\pm}]=0.
\endaligned
\eqno{(3.1.6)}
$$

\bigskip

{\bf 3.1.2. Super-jets.} Let $M$ be a $C^{\infty}$-supermanifold
with underlying even manifold $M^{even}$. Define
$$
M_{\hat{\Sigma}}=M\times\hat{\Sigma}.\eqno{(3.1.7)}
$$
It is fibered over $\hat{\Sigma}$:
$$
M_{\hat{\Sigma}}\rightarrow\hat{\Sigma}.\eqno{(3.1.8)}
$$
The manifold of  $\infty$-jets of sections of this bundle,
$J^{\infty}(M_{\hat{\Sigma}})$, is defined in a straightforward
manner as follows (cf. [BD p.80]).

\bigskip

{\bf 3.1.2.1. Definition.} {\it $J^{\infty}(M_{\hat{\Sigma}})$ is a
supermanifold with underlying even manifold
$J^{\infty}(M^{even}_{\Sigma})$ and the structure sheaf
$\CO_{J^{\infty}(M_{\hat{\Sigma}})}$ defined to be the symmetric
algebra on
$$
D_{\hat{\Sigma}}\otimes_{\CO_{\hat{\Sigma}}} \CO_{M_{\Sigma}}
$$
modulo the relations
$$\aligned
&1\otimes f\cdot 1\otimes g=1\otimes fg,\; 1\otimes 1=1,\\
 &\xi\otimes
fg=(\xi\otimes f)\cdot(1\otimes
g)+(-1)^{\tilde{\xi}\tilde{f}}(1\otimes f)\cdot(\xi\otimes g).
\endaligned
\eqno{(3.1.9)}
$$
for any $\xi\in\CT_{\hat{\Sigma}}$, $f,g\in \CO_{M_{\Sigma}}$, where
$\tilde{ }$ stands for the parity.}

\bigskip

There arises a fiber bundle
$$
J^{\infty}(M_{\hat{\Sigma}})\rightarrow\hat{\Sigma}\eqno{(3.1.10)}
$$
with connection
$$
\rho:\CT_{\hat{\Sigma}}\rightarrow\CT_{J^{\infty}(M_{\hat{\Sigma}})}
\text{ s.t. }\rho(\eta)(\xi\otimes f)=(\eta\xi)\otimes
f\eqno{(3.1.11)}
$$
in complete analogy with (1.1.5).

The relative versions, such as
$J^{\infty}(M_{\hat{\Sigma}/\Sigma''})$, are immediate.

Note that connection (3.1.11) is constant in the direction of
$(\theta^{+},\theta^{-})$, i.e.,  if we let
$J^{\infty}(M_{\hat{\Sigma}})^{o}=\{\theta^{+}=\theta^{-}=0\}\hookrightarrow
J^{\infty}(M_{\hat{\Sigma}})$, then there is a diffeomorphism
$$
\left(
J^{\infty}(M_{\hat{\Sigma}}),\rho(\partial_{\theta^{\pm}})\right)\rightarrow
\left(J^{\infty}(M_{\hat{\Sigma}})^{o}\times
\BR^{0|2},\rho^{o}(\partial_{\theta^{\pm}})=\partial_{\theta^{\pm}}\right)
\eqno{(3.1.12)}
$$
of $\BR^{0|2}$-manifolds with connection.

Indeed, given a local coordinate system $X^{i}$ on $M$, the
collection
$$
\{X^{i}_{(m),(\epsilon)},u,v,\theta^{+},\theta^{-};
(m)\in\BZ_{+}^{2},(\epsilon)\in\BZ_{2}^{2}\}\eqno{(3.1.13a)}
$$
constitutes a local coordinate system on
$J^{\infty}(M_{\hat{\Sigma}})$, where
$$
X^{i}_{(m_{1},m_{2}),(\epsilon_{1},\epsilon_{2})}=
(\partial_{u}^{m_{1}}\partial_{v}^{m_{2}}\partial_{\theta^{+}}^{\epsilon_{1}}\partial_{\theta^{-}}^{\epsilon_{2}})\otimes
X.\eqno{(3.1.13b)}
$$
Letting
$$
\aligned
&\tilde{F}^{i}=(\partial_{\theta^{-}}\partial_{\theta^{+}})\otimes X^{i},\\
&\psi^{i}_{+}=(\partial_{\theta^{+}})\otimes X^{i}-\theta^{-}\tilde{F}^{i},\\
&\psi^{i}_{-}=(\partial_{\theta^{-}})\otimes X^{i}+\theta^{+}\tilde{F}^{i},\\
&x^{i}=X^{i}-\theta^{+}\psi^{i}_{+}-\theta^{-}\psi^{i}_{-}-\theta^{+}\theta^{-}\tilde{F}^{i},
\endaligned
\eqno{(3.1.14)}
$$
we obtain another local coordinate system
$$
\{x^{i}_{(m)},\psi^{i}_{\pm,(m)},
\tilde{F}^{i}_{(m)};u,v,\theta^{+},\theta^{-};
(m)\in\BZ_{+}^{2},(\epsilon)\in\BZ_{2}^{2}\}\eqno{(3.1.15)}
$$
such that
$$
\partial_{\theta^{\pm}}x^{i}_{(m)}=\partial_{\theta^{\pm}}\psi^{i}_{\pm,(m)}=
\partial_{\theta^{\pm}}\tilde{F}^{i}_{(m)}=0,\eqno{(3.1.16)}
$$
and (3.1.12) follows.

Note that  change of variables  (3.1.14) is nothing but the formal
Taylor series expansion at $J^{\infty}(M_{\hat{\Sigma}})^{o}$:
$$
X^{i}=x^{i}+\theta^{+}\psi^{i}_{+}+\theta^{-}\psi^{i}_{-}+\theta^{+}\theta^{-}\tilde{F}^{i}.\eqno{(3.1.17)}
$$
Along $M$, $\{x^{i}\}$ are coordinates and
$$
\psi_{\pm}^{i}\text{ transform as (even or odd)
}dx^{i}.\eqno{(3.1.18)}
$$

\bigskip

{\bf 3.1.3. Differential equations.} The definition and discussion
of a submanifold $Sol\subset J^{\infty}(M_{\hat{\Sigma}})$ as the
zero locus of a differential ideal $\CJ$ is quite parallel to sect.
1.4.  Since our exposition is strongly focused on one particular
example, that of the (2,2)-supersymmetric $\sigma$-model, we will
restrict ourselves to the case where $\CJ$ is locally generated by
4n functions, $E^{i}_{\alpha}$, $1\leq i\leq n$, $1\leq \alpha\leq
4$, such that (cf. (1.4.1))
$$
\aligned
&E^{i}_{1}=\tilde{F}^{i}+\cdots\\
&E^{i}_{2}=\partial_{\tau}\psi^{i}_{-}+\cdots\\
&E^{i}_{3}=\partial_{\tau}\psi^{i}_{+}+\cdots\\
&E^{i}_{4}=\partial_{\tau}^{2}x^{i}+\cdots,
\endaligned
\eqno{(3.1.19)}
$$
where  the omitted terms are independent of $\tilde{F}^{\cdot}$, of
non-zero order jets of $\psi_{\pm}^{\cdot}$ in the direction of
$\tau$, and of order $>1$ jets of $x^{\cdot}$ also in the direction
of $\tau$. ($\tau$ is time-function (3.1.4) tacitly assumed to have
been included in a coordinate system.)

Letting
$$
Sol^{o}=\{\theta^{+}=\theta^{-}\}\hookrightarrow Sol,\eqno{(3.1.20)}
$$
one obtains a diffeomorphism of $\BR^{0|2}$-manifolds with
connection
$$
(Sol,\rho(\partial_{\theta^{\pm}}))\iso
(Sol^{o}\times\BR^{0|2},\rho^{o}(\partial_{\theta^{\pm}})=\partial_{\theta^{\pm}}),\eqno{(3.1.21)}
$$
by restricting (3.1.12).

If (3.1.19) holds, then (3.1.18) implies  a diffeomorphism of
$D_{\Sigma/\Sigma''}$-manifolds
$$
Sol_{\Sigma/\Sigma''}^{o}\iso J^{\infty}\left(T\left(\Pi
TM\right)_{\Sigma/\Sigma''}\right),\eqno{(3.1.22)}
$$
where $\Pi$ is the familiar parity change functor. Similarly,
$$
Sol_{\hat{\Sigma}/\Sigma''}\iso J^{\infty}\left(T\left(\Pi
TM\right)_{\Sigma/\Sigma''}\right)\times\BR^{0|2}\eqno{(3.1.23)}
$$
as $D_{\hat{\Sigma}/\Sigma''}$-manifolds. Both (3.1.22,23) are
analogous to (1.4.2).

\bigskip

{\bf 3.2. Homotopy presymplectic structure.}

\bigskip
The right framework for super-generalization of 1.5 is provided by
integral, rather than differential, forms on $\hat{\Sigma}$
[L,M1,DM].

\bigskip
{\bf 3.2.1.} Recall that the sheaf of integral forms is defined to
be
$$
I^{*}_{\hat{\Sigma}}=\bigoplus_{i=-\infty}^{4}I^{i}_{\hat{\Sigma}}\text{
s.t.
}I^{4-i}_{\hat{\Sigma}}=\Lambda^{i}\CT_{\hat{\Sigma}}\otimes_{\CO_{\hat{\Sigma}}}
\text{Ber}(\Omega_{\hat{\Sigma}}),\eqno{(3.2.1)}
$$
where $\text{Ber}(\Omega_{\hat{\Sigma}})$ is the Berezinian of
$\Omega_{\hat{\Sigma}}$.

By definition, $I^{*}_{\hat{\Sigma}}$ is a locally free
$\Lambda^{*}\CT_{\hat{\Sigma}}$-module defined by
$$
\CT_{\hat{\Sigma}}\rightarrow
End_{\CO_{\hat{\Sigma}}}(I^{*}_{\hat{\Sigma}}),\;
\xi\mapsto\iota_{\xi}\text{ where
}\iota_{\xi}\beta\buildrel\text{def}\over=\xi\wedge\beta.
\eqno{(3.2.2)}
$$
Next, $I^{*}_{\hat{\Sigma}}$ carries a unique structure of a module
over the Clifford algebra,
$Cl(\CT_{\hat{\Sigma}}\oplus\Omega_{\hat{\Sigma}})$, such that
$$
\Omega_{\hat{\Sigma}}^{i}\otimes_{\CO_{\hat{\Sigma}}}I^{j}_{\hat{\Sigma}}\rightarrow
I^{i+j}_{\hat{\Sigma}}
,\;\alpha\otimes\beta\mapsto\alpha\wedge\beta,\;
[\iota_{\xi},\alpha\wedge]=\alpha(\xi). \eqno{(3.2.3)}
$$
The Berezinian, $\text{Ber}(\Omega_{\hat{\Sigma}})$, carries the Lie
derivative operation
$$
\CT_{\hat{\Sigma}}\otimes_{\BR}\text{Ber}(\Omega_{\hat{\Sigma}})\rightarrow
\text{Ber}(\Omega_{\hat{\Sigma}}),\;
\xi\otimes\beta\mapsto\text{Lie}_{\xi}\beta,\eqno{(3.2.4)}
$$
which is naturally extended to
$$
\CT_{\hat{\Sigma}}\otimes_{\BR}I^{i}_{\hat{\Sigma}}\rightarrow
I^{i}_{\hat{\Sigma}},\;
\xi\otimes\beta\mapsto\text{Lie}_{\xi}\beta.\eqno{(3.2.5)}
$$
The sheaf of integral forms is a complex with differential
$$
d:I^{i}_{\hat{\Sigma}}\rightarrow
I^{i+1}_{\hat{\Sigma}}\eqno{(3.2.6a)}
$$
determined by
$$
[d,\iota_{\xi}]=\text{Lie}_{\xi},\;\xi\in\CT_{\hat{\Sigma}}.\eqno{(3.2.6b)}
$$
Many other differential-geometric identities, such as
$$
[\text{Lie}_{\xi},\iota_{\eta}]=\iota_{[\xi,\eta]},\;
[\text{Lie}_{\xi},\beta\wedge ]=(\text{Lie}_{\xi}\beta)\wedge,\;
\xi,\eta
\in\CT_{\hat{\Sigma}},\beta\in\Omega_{\hat{\Sigma}},\eqno{(3.2.7)}
$$
keep on holding true.

Since our particular $\hat{\Sigma}$ carries a fixed (\'etale)
coordinate system $(u,v,\theta^{+},\theta^{-}$, there is an integral
form $[d\theta^{+}d\theta^{-}]$ such that $du\wedge dv\wedge
[d\theta^{+}d\theta^{-}]$ trivilaizes the Berezinian
$\text{Ber}(\Omega_{\hat{\Sigma}})$. Letting
$$
[d\theta^{\pm}]=\iota_{\partial_{\theta^{\mp}}}[d\theta^{+}d\theta^{-}]=[d\theta^{\pm}]
$$
one discovers a part of $I^{*}_{\hat{\Sigma}}$ pleasingly -- and
deceptively -- similar to the de Rham complex; e.g.,
$$
\text{Lie}_{\partial_{\theta^{\pm}}}[d\theta^{+}d\theta^{-}]=0,\;
d([d\theta^{+}d\theta^{-}])=d([d\theta^{\pm}])=0,d(\theta^{\pm}[d\theta^{\mp}]=[d\theta^{+}d\theta^{-}].
\eqno{(3.2.8)}
$$
Once a projection
$$
\hat{\Sigma}\rightarrow\Sigma
$$
is given, integration over fibers delivers a morphism
$$
I_{\hat{\Sigma}}^{4}\rightarrow\Omega^{2}_{\Sigma},\eqno{(3.2.9)}
$$
which, in the case where the projection is (3.1.3), means that
$$
f(u,v,\theta^{+},\theta^{-})du\wedge dv\wedge
[d\theta^{+}d\theta^{-}]\mapsto
\partial_{\theta^{-}}\partial_{\theta^{+}}f(u,v,\theta^{+},\theta^{-})du\wedge
dv, \eqno{(3.2.10)}
$$
cf. (3.2.8). This is often referred to as integrating out
$\theta^{+}$ and $\theta^{-}$.

\bigskip

{\bf 3.2.2.} Back to super-presymplectic forms. Let $\CM$ be either
$Sol$ or $Sol^{o}$ or any version of an $\infty$-jet space
considered in 3.1.2, cf. 1.5. Suppose $\CM$ is fibered over
$\hat{\Sigma}$ and let
$$
\tilde{\Omega}^{*,*}_{\CM}=\Omega^{*}_{\CM/\hat{\Sigma}}\otimes_{\CO_{\hat{\Sigma}}}
I^{*}_{\hat{\Sigma}}.\eqno{(3.2.11)}
$$
If we wish to work in a relative situation determined by $\tau$,
(3.1.4), then we  write
$$
\tilde{\Omega}^{*,*}_{\CM/\Sigma''}=\Omega^{*}_{\CM/\hat{\Sigma}}\otimes_{\CO_{\hat{\Sigma}}}
I^{*}_{\hat{\Sigma}/\Sigma''}.\eqno{(3.2.12)}
$$
In any case, we get a bi-complex with an obvious vertical
differential
$$
\delta: \tilde{\Omega}^{*,i}_{\CM/S}\rightarrow
\tilde{\Omega}^{*,i+1}_{\CM/S}\eqno{(3.2.13)}
$$
 and a horizontal differential
 $$
 d_{\rho/S}: \tilde{\Omega}^{i,*}_{\CM/S}\rightarrow
\tilde{\Omega}^{i+1,*}_{\CM/S}, \eqno{(3.2.14)}
 $$
 which owes its existence to connection (3.1.11) and is defined in
 exactly the same way as its counterpart in sect. 3.1; here and
 elsewhere $S$ is either $\Sigma''$ or a point.

 With $\tilde{\Omega}^{*,*}_{\CM/\Sigma''}$ taken for a replacement
 of $\Omega^{*,*}_{\CM/\Sigma''}$, the discussion of 1.5.1-3 carries
 over to the super-case practically word for word. For example, cf. (1.5.1), a
 homotopy pre-symplectic form is $\omega\in H^{0}(\CM,\tilde{\Omega}_{\CM/S}^{3,2}/d_{\rho}\tilde{\Omega}_{\CM/S}^{2,2})$
 such that,
$$
\delta\omega\in H^{0}\left(\CM,
d_{\rho/S}\left(\tilde{\Omega}_{\CM/S}^{2,3}\right)\right),
\eqno{(3.2.15)}
$$
cf. (1.5.1) and Remark 1.5.3.1. The outcome is the Lie superalgebra
sheaf over $\CM$
$$
\tilde{\CH}^{\omega}_{\CM/S}.\eqno{(3.2.16)}
$$
Here is an operation that does not have an adequate purely even
analogue. In all our examples,
$$
\CM\iso\CM^{o}\times\BR^{0|2},
$$
in a way respecting the connection, cf. (3.1.12, 21, 23).
Integrating over fibers one obtains, as in (3.2.9), a Lie algebra
sheaf morphism, an isomorphism in fact,
$$
\tilde{\CH}^{\omega}_{\CM/S}\iso
\tilde{\CH}^{\omega}_{\CM^{o}/S}.\eqno{(3.2.17)}
$$
As a practical matter, (3.2.17) amounts to carrying out a Taylor
expansion as in (3.1.17) and then extracting the coefficient of
$\theta^{+}\theta^{-}[d\theta^{+}d\theta^{-}]$ as in (3.2.10).

Similar in spirit is a morphism
$$
\tilde{\CH}^{\omega}_{\CM}\rightarrow
\tilde{\CH}^{\omega}_{\CM/\Sigma''}\eqno{(3.2.18)}
$$
that relates the relative and absolute versions and amounts to
letting $d\tau=0$, cf. (1.5.11).

\bigskip

{\bf 3.2.3.} {\it Example: canonical commutation relations.}

Let $M$ be an n-dimensional purely even $C^{\infty}$-manifold. The
$2n|2n$-dimensional supermanifold $T^{*}(\Pi TM)$ carries a
well-known closed 2-form $\omega^{o}$. If we let $\{x^{i}\}$ be
coordinates on $M$, then $\{x^{i}, x_{i}=\partial_{x^{i}}\}$ along
with their superpartners $\{\phi^{i},\phi_{i}\}$ form a system of
local coordinates on $T^{*}(\Pi TM)$, and
$$
\omega^{o}=\delta x_{i}\wedge\delta x^{i}+\delta
\phi_{i}\wedge\delta \phi^{i}.\eqno{(3.2.19)}
$$
Now use the projection
$$
\pi: J^{\infty}(T^{*}(\Pi TM)_{\Sigma/\Sigma''})\rightarrow
T^{*}(\Pi TM)
$$
to introduce $\pi^{*}\omega^{o}$, a closed 2-form on
$J^{\infty}(T^{*}(\Pi TM)_{\Sigma/\Sigma''})$.

 A suitable analogue of $\CH^{can}$, see 1.5.4, is provided by
 fixing a suitable $\sigma$ so that $\sigma$, $\tau$ is a coordinate
 system on $\Sigma$, see (3.1.4), letting
$$
\omega=\pi^{*}\omega^{o}\wedge d\sigma, \eqno{(3.2.20)}
$$
and defining
$$
\tilde{\CH}^{can}=\CH^{\omega}_{J^{\infty}(T^{*}(\Pi
TM)_{\Sigma/\Sigma''})}. \eqno{(3.2.21)}
$$
The rest of the discussion in 1.5.4 carries over to the present
situation practically word for word; we will not dwell upon this any
longer.

\bigskip

{\bf 3.2.4.} {\it Legendre transform?} In practice, the manifold
$J^{\infty}(T(\Pi TM)_{\Sigma/\Sigma''})$ may be more important than
$J^{\infty}(T^{*}(\Pi TM)_{\Sigma/\Sigma''})$ because of (3.1.22).
The possibility to apply $\tilde{\CH}^{can}$ then rests on the
existence of the diffeomorphism, cf. (1.5.19),
$$
g: J^{\infty}(T(\Pi TM)_{\Sigma/\Sigma''})\rightarrow
J^{\infty}(T^{*}(\Pi TM)_{\Sigma/\Sigma''}),\eqno{(3.2.22)}
$$
because given (3.2.21) there arises at once a Lie algebra sheaf
ismorphism, cf. Lemma 1.5.5.1,
$$
g^{\#}:\CH^{g^{*}\omega}\iso g^{-1}\tilde{\CH}^{can}.
\eqno{(3.2.23)}
$$
Isomorphism (3.2.22), however, is a more subtle matter in the
present situation than the usual Legendre transform. While the
purely even manifolds underlying both the manifolds in (3.2.22) are
the familiar $J^{\infty}(TM_{\Sigma/\Sigma''})$ and
$J^{\infty}(T^{*}M_{\Sigma/\Sigma''})$, and they are easy to
identify via a metric, the structure sheaves are more substantially
different. The essence of this difference is that while
$$
\Omega_{M}\hookrightarrow \CO_{T(\Pi TM)}\eqno{(3.2.24)}
$$
 as a direct summand, its $T^{*}(\Pi TM)$-counterpart, $\CT_{M}$,
appears via the extension
$$
0\rightarrow
\text{End}\Omega_{M}\rightarrow\CA_{\Omega_{M}}\rightarrow\CT_{M}\rightarrow
0, \eqno{(3.2.25)}
$$
where $\CA_{\Omega_{M}}$ is the Atiyah algebra, i.e., the algebra of
order 1 differential operators acting on the sections of
$\Omega_{M}$. One way to construct (3.2.22) seems to be this:  split
(3.2.25) by means of a connection
$$
0\rightarrow
\text{End}\Omega_{M}\rightarrow\CA_{\Omega_{M}}\buildrel
{\buildrel\nabla\over\leftarrow}\over\rightarrow\CT_{M}\rightarrow
0, \eqno{(3.2.26)}
$$
and then identify
$$
\Omega_{M}\iso\CT_{M}\eqno{(3.2.27)}
$$
by means of a metric. This is exactly what the Lagrangian of a
(1,1)-supersymmetric $\sigma$-model allows to do.

\bigskip

{\bf 3.3. Calculus of variations.}

\bigskip

{\bf 3.3.1.} The discussion of sect. 1.6 carries over in a
straightforward manner. Here are a few highlights.

An action is
$$
S\in H^{0}\left(J^{\infty}\left(M_{\hat{\Sigma}}\right),
\tilde{\Omega}^{4,0}_{J^{\infty}\left(M_{\hat{\Sigma}}\right)}/
d_{\rho}\tilde{\Omega}^{3,0}_{J^{\infty}\left(M_{\hat{\Sigma}}\right)}\right).
\eqno{(3.3.1)}
$$
It is represented by a collection of Lagrangians
$$
\{L^{(j)}\in
\tilde{\Omega}^{4,0}_{J^{\infty}\left(M_{\hat{\Sigma}}\right)}(U_{j})\}
\eqno{(3.3.2)}
$$
determined up to $d_{\rho}$-exact terms and equal to each other on
intersections $U_{i}\cap U_{j}$ up to $d_{\rho}$-exact terms, cf.
(1.6.0-1ab).

An analogue of (1.6.2) is immediate, the outcome is the
supermanifold $Sol_{L}$ over $\hat{\Sigma}$ with a variational
1-form $\gamma_{L}$ and 2-form $\omega_{L}=\delta\gamma_{L}$.

The definition of a symmetry of $L$ is also an obvious modification
of (1.6.5). N\"other's Theorem establishes a bijection between
symmetries and integrals of motion as follows
$$
\xi\leftrightarrow
\alpha_{\xi}+(-1)^{\tilde{\xi}+1}\iota_{\xi}\gamma_{L},
\eqno{(3.3.3)}
$$
cf. (1.6.7); the change of sign occurs when swapping $\iota_{\xi}$
and $d_{\rho}$ as in the last formula of sect.1.6.2.

Thus there arise the Lie algebra sheaf $\CH^{\omega_{L}}_{Sol_{L}}$,
containing the algebra of first integrals $\tilde{I}_{L}$, its
relative version, $\CH^{\omega_{L}}_{Sol_{L}/\Sigma''}$, and
morphisms
$$
 \CH^{\omega_{L}}_{Sol_{L}}\rightarrow
\CH^{\omega_{L}}_{Sol_{L}/\Sigma''},\; \tilde{I}_{L}\hookrightarrow
\Gamma(Sol_{L},\CH^{\omega_{L}}_{Sol_{L}})\rightarrow
\Gamma(Sol_{L},\CH^{\omega_{L}}_{Sol_{L}/\Sigma''})\eqno{(3.3.4)}
$$
whose composition is an injection provided (3.1.19) holds.

A familiar novelty is that in all of this $\theta^{\pm}$ can be
integrated out. The result is this: the action
$$
S\in
H^{0}\left(\left(J^{\infty}\left(M_{\hat{\Sigma}}\right)\right)^{o},
\Omega^{2,0}_{\left(J^{\infty}\left(M_{\hat{\Sigma}}\right)\right)^{o}}/
d_{\rho}\Omega^{1,0}_{\left(J^{\infty}\left(M_{\hat{\Sigma}}\right)\right)^{o}}\right),
\eqno{(3.3.5)}
$$
 (cf. (3.1.12), the Lagrangian
$$
\{L^{(j)}\in
\Omega^{2,0}_{\left(J^{\infty}\left(M_{\hat{\Sigma}}\right)\right)^{o}}(U_{j})\}
\eqno{(3.3.6)},
$$
and, since nothing is gained or lost, the integrated version of
(3.3.4) as follows
$$
 \CH^{\omega_{L}}_{Sol_{L}^{o}}\rightarrow
\CH^{\omega_{L}}_{Sol_{L}^{o}/\Sigma''},\;
\tilde{I}_{L}\hookrightarrow
\Gamma(Sol_{L}^{o},\CH^{\omega_{L}}_{Sol_{L}^{o}})\rightarrow
\Gamma(Sol_{L}^{o},\CH^{\omega_{L}}_{Sol_{L}^{o}/\Sigma''}),\eqno{(3.3.7)}
$$
where $Sol_{L}^{o}$ is as defined in (3.1.20)

In view of (3.1.22), this means $J^{\infty}(T(\Pi
TM)_{\Sigma/\Sigma''})$ equipped with $\omega_{L}$ and an embedding
$$
\tilde{I}_{L}\hookrightarrow\Gamma(M,
\CH^{\omega_{L}}_{J^{\infty}(T(\Pi
TM)_{\Sigma/\Sigma''})}).\eqno{(3.3.8)}
$$
We will now exhibit an example where
$$
\CH^{\omega_{L}}_{J^{\infty}(T(\Pi
TM)_{\Sigma/\Sigma''})}\iso\tilde{\CH}^{can}.\eqno{(3.3.9)}
$$

\bigskip

{\bf 3.4. An example: (1,1)-supersymmetric $\sigma$-model.}

\bigskip

{\bf 3.4.1.} Let $M$ be an $n$-dimensional purely even Riemannian
manifold with metric $(.,.)$. Analogously to sect. 1.6.11, we
observe that a point in $J^{1}(M_{\hat{\Sigma}})$ is a triple
$(\hat{t}, X,
\partial X)$, a point in $\hat{\Sigma}$, a point in $M$, and a map
$$
\partial X: T_{\hat{t}}\hat{\Sigma}\rightarrow T_{X}M,\;
\xi\mapsto \partial_{\xi}X,\eqno{(3.4.1)}
$$
cf. (3.1.13a,b).  Hence for fixed vector fields $\xi$, $\eta$,
$(\partial_{\xi}X,\partial_{\eta}X)$ is a global section of
$\CO_{J^{1}(M_{\hat{\Sigma}})}$. (Of course, to be precise, we
should have used $B$-points.) We will unburden the notation by
letting $\xi X$ stand for $\partial_{\xi}X$. Here is a coordinate
expression for this function
$$
g_{ij}(X)\xi X^{i}\eta X^{j}.\eqno{(3.4.1)}
$$

The (1,1)-supersymmetric $\sigma$-model Lagrangian is defined to be
$$
L=(D_{+}X,D_{-}X) du\wedge dv\wedge[d\theta^{+}d\theta^{-}],
\eqno{(3.4.2)}
$$
where the vector fields $D_{\pm}$ are from (3.1.5), cf. (1.6.13).
Integrating out $\theta^{+}$ and $\theta^{-}$ gives (an exercise in
differential geometry, see e.g. [QFS, p. 666])
$$
\aligned
L_{11}=&(-(\partial_{u}x,\partial_{v}x)+(\nabla_{\partial_{v}x}\psi_{+},\psi_{+})+
(\nabla_{\partial_{u}x}\psi_{-},\psi_{-})+\\
&(R(\psi_{+},\psi_{-})\psi_{+},\psi_{-})-(F,F))du\wedge dv,
\endaligned
\eqno{(3.4.3)}
$$
In this formula
$$
\partial_{u}x=\partial_{u}X|_{\theta^{+}=\theta^{-}=0},
\partial_{v}x=\partial_{v}X|_{\theta^{+}=\theta^{-}=0},
\psi_{\pm}=\partial_{\theta^{\pm}}X|_{\theta^{+}=\theta^{-}=0}
\eqno{(3.4.4a)}
$$
and coincide with their namesakes from (3.1.14), $\nabla$ is the
Levi-Civita connection associated to the metric $(.,.)$,   $R$ is
the curvature tensor, and
$$
F=\nabla_{D_{+}X}D_{-}X|_{\theta^{+}=\theta^{-}=0}.
\eqno{(3.4.4b)}
$$

In fact, with a little extra effort the entire Taylor series
expansion of $L$ in $\theta^{+}$, $\theta^{-}$, cf. (3.1.17), can be
computed to the effect that
$$
\aligned
L=(&(\psi_{+},\psi_{-})-\\
&\theta^{+}((\partial_{u}x,\psi_{-})+(\psi_{+},F))+\theta_{-}((\partial_{v}x,\psi_{+})-(\psi_{-},F))+\\
&\theta^{+}\theta^{-}L_{11})du\wedge
dv\wedge[d\theta^{+}d\theta^{-}].
\endaligned
\eqno{(3.4.5)}
$$

To see better what all of this means, let us write down the first
three terms of  (3.4.3) in local coordinates (3.1.14,17); the result
is
$$
\aligned
 L_{11}=&(-g_{ij}(x)\partial_{u}x^{i}\partial_{v}x^{j}+
\partial_{v}x^{\alpha}\psi_{+}^{i}\psi_{+}^{j}\Gamma^{s}_{\alpha
i}g_{sj}(x)+
\partial_{u}x^{\alpha}\psi_{-}^{i}\psi_{-}^{j}\Gamma^{s}_{\alpha i}g_{sj}(x)
+\\
&g_{ij}(x)\psi_{+}^{j}\delta\psi_{+}^{i}+
g_{ij}(x)\psi_{-}^{j}\delta\psi_{-}^{i}\cdots)du\wedge dv.
\endaligned
\eqno{(3.4.6)}
$$

Computation of $\delta L_{11}$, cf. (1.6.2), yields the
Euler-Lagrange equations and a variational 1-form. The former are as
follows, see also [QFS, p.666],
$$
\aligned F&=0\\
\nabla_{\partial_{u}x}\psi_{-}&=-R(\psi_{+},\psi_{-})\psi_{+}\\
\nabla_{\partial_{v}x}\psi_{+}&=-R(\psi_{+},\psi_{-})\psi_{-}\\
\nabla_{\partial_{u}x}\partial_{v}x&=
\frac{1}{2}\left(R(\psi_{-},\psi_{-})\partial_{v}x+R(\psi_{+},\psi_{+})\partial_{u}x\right)
-(\nabla_{\psi_{+}}R)(\psi_{-},\psi_{-})\psi_{+}
\endaligned
\eqno{(3.4.7)}
$$
The latter is
$$
\aligned
 \gamma_{L}=&\left(-g_{ij}(x)\partial_{v}x^{j}\delta x^{i}+
 \psi_{-}^{i}\psi_{-}^{j}\Gamma^{s}_{\alpha i}g_{sj}(x)\delta
x^{\alpha}\right)dv+\\
&\left(g_{ij}(x)\partial_{u}x^{i}\delta
x^{j}-\psi_{+}^{i}\psi_{+}^{j}\Gamma^{s}_{\alpha i}g_{sj}(x)\delta
x^{\alpha}\right)du+\\
&-g_{ij}(x)\psi_{+}^{j}\delta\psi_{+}^{i}du+g_{ij}(x)\psi_{-}^{j}\delta\psi_{-}^{i}dv.
\endaligned
\eqno{(3.4.8)}
$$
This, unlike more challenging (3.4.7), is a straightforward
consequence of (3.4.6). Note that we have computed after projection
(3.2.9), i.e., with $\theta^{\pm}$ integrated out, (3.2.10); nothing
is gained or lost because $\gamma_{L}$ matters only modulo
$d_{\rho}$-exact terms.

\bigskip

{\bf 3.4.2.} {\it (1,1)-supersymmetry.}

In addition to $D_{\pm}$, (3.1.5), $\hat{\Sigma}$ carries another
pair of distinguished vector fields,
$$
\Xi_{+}=\frac{\partial}{\partial\theta^{+}}+\theta^{+}\frac{\partial}{\partial
u},\;
\Xi_{-}=\frac{\partial}{\partial\theta^{-}}+\theta^{-}\frac{\partial}{\partial
v}.\eqno{(3.4.9)}
$$
They enjoy similar properties
$$
\aligned &[\Xi_{+},\Xi_{+}]=2\frac{\partial}{\partial u},
[\Xi_{-},\Xi_{-}]=2\frac{\partial}{\partial v}, [\Xi_{+},\Xi_{-}]=0\\
&[\frac{\partial}{\partial v},\Xi_{\pm}]= [\frac{\partial}{\partial
u},\xi_{\pm}]=0,
\endaligned
\eqno{(3.4.10)}
$$
and
$$
[\Xi_{\bullet},D_{\bullet}]=0. \eqno{(3.4.11)}
$$
Relations (3.4.10) imply that
$$
\aligned
&\CN1^{+}\buildrel\text{def}\over=\text{span}\{f(u)\Xi_{+}\}\subset
\CT_{\hat{\Sigma}}(\hat{\Sigma})\\
&\CN1^{-}\buildrel\text{def}\over=\text{span}\{f(v)\Xi_{-}\}\subset
\CT_{\hat{\Sigma}}(\hat{\Sigma})
\endaligned
\eqno{(3.4.12)}
$$
are two commuting copies of the N=1- supersymmetric superalgebra Lie
realized in vector fields on $\hat{\Sigma}$; note that each contains
a copy, $\Vir^{\pm}$, of the algebra of vector fields on $\Sigma$.

In fact, both are subalgebras of the algebra of symmetries of $L$:
$$
\CN1^{\pm}\hookrightarrow\tilde{I}_{L}.\eqno{(3.4.13)}
$$
Indeed, using (3.4.5), one computes easily that
$$
\text{Lie}_{\rho f(u)\Xi_{+}}L_{11}=d_{\rho}(f(u)L_{01})dv,\;
\text{Lie}_{\rho f(v)\Xi_{-}}L_{11}=d_{\rho}(f(v)L_{10})du.
\eqno{(3.4.14)}
$$
It is then rather straightforward, and pleasing, to use (3.3.3,
3.4.8) in order to compute the corresponding integrals of motion
$$
\aligned
 Q^{+}_{f}\buildrel\text{def}\over=Q_{\rho
f(u)\Xi_{+}}=2f(u)g_{ij}(x)\psi_{+}^{i}\partial_{u}x^{j}du-g_{ij}(x)F^{i}\psi_{-}^{j}dv\\
Q^{-}_{f}\buildrel\text{def}\over=Q_{\rho
f(v)\Xi_{-}}=-2f(v)g_{ij}(x)\psi_{-}^{i}\partial_{v}x^{j}dv+g_{ij}(x)F^{i}\psi_{+}^{j},
\endaligned
\eqno{(3.4.15)}
$$
or, taking advantage of the Euler-Lagrange equation $F=0$, see
(3.4.7), simply
$$
 Q^{+}_{f}=2f(u)g_{ij}(x)\psi_{+}^{i}\partial_{u}x^{j}du,\;
Q^{-}_{f}=-2f(v)g_{ij}(x)\psi_{-}^{i}\partial_{v}x^{j}dv.
\eqno{(3.4.16)}
$$
This furnishes the embeddings
$$
\Vir^{\pm}\hookrightarrow
\CN1^{\pm}\hookrightarrow\Gamma(Sol_{L}^{o},\CH^{\omega_{L}}_{Sol_{L}^{o}}),
\eqno{(3.4.17a)}
$$ and the definitions (cf. 1.6.11.1) of right/left moving
subalgebras
$$
\CH^{\omega_{L},\pm}_{Sol_{L}^{o}}=\{F\in\CH^{\omega_{L}}_{Sol_{L}^{0}}:\;
[F,\Vir^{\mp}]=0\}.\eqno{(3.4.17b)}
$$
Next, we will see that all of this unfolds within the canonical Lie
algebra sheaf $\tilde{\CH}^{can}$ of sect. 3.2.3.

\bigskip

{\bf 3.4.3. Proposition.} {\it There is a diffeomorphism
$$
\hat{g}:\; Sol_{L}^{o}\iso J^{\infty}(T^{*}(\Pi
TM)_{\Sigma/\Sigma''})
$$
of $D_{\Sigma/\Sigma''}$-manifolds, which delivers the Lie algebra
sheaf isomorphism
$$
g^{\#}:\CH^{\omega_{L},\pm}_{Sol_{L}^{o}/\Sigma''}\iso
\hat{g}^{-1}\tilde{\CH}^{can},
$$
cf. Lemmas 1.5.5.1 and 1.6.10.1. }

\bigskip

{\bf 3.4.4.} {\it Proof. Super-Legendre transform.}

In order to proceed, we need to make sure that $Sol_{L}$ as defined
by (3.4.7) satisfies Cauchy-Kovalevskaya condition (3.1.19).
Apparently neither $u$ nor $v$ can play the role of time, but the
change variables as follows
$$
u=\sigma+\tau,\; v=\sigma-\tau \eqno{(3.4.18)}
$$
so that
$$
\partial_{u}=\frac{1}{2}(\partial_{\sigma}+\partial_{\tau}),\;
\partial_{v}=\frac{1}{2}(\partial_{\sigma}-\partial_{\tau}) \eqno{(3.4.19)}
$$
does the job. Therefore, cf. (3.1.22),
$$
Sol_{L}^{o}\iso J^{\infty}(T(\Pi TM)_{\Sigma/\Sigma''}),
\eqno{(3.4.20)}
$$
and our task is to find
$$
p: J^{\infty}(T(\Pi TM)_{\Sigma/\Sigma''})\iso J^{\infty}(T^{*}(\Pi
TM)_{\Sigma/\Sigma''})\eqno{(3.4.21)}
$$
that identifies $\omega_{L}$ on the L.H.S. with the pull-back of the
canonical $\omega^{o}$, (3.2.19), on the R.H.S
$$
p^{*}\omega|_{d\tau=0} =\omega_{L}|_{d\tau=0}+
d_{\rho/\Sigma''}(\cdots) \eqno{(3.4.22)}
$$
modulo $d_{\rho/\Sigma''}$-exact terms. Variational 1-form
$\gamma_{L}$, computed in (3.4.8), is not well suited for this
purpose. In addition to (3.4.18), let us introduce variables
$$
\rho^{j}=\frac{1}{2}(\psi_{+}^{j}+\psi_{-}^{j}),
\phi^{j}=\frac{1}{2}(\psi_{-}^{j}-\psi_{-}^{j}).\eqno{(3.4.23)}
$$
(Since $\psi_{\pm}^{j}$ are sections of 2 copies of the bundle of
1-forms, see (3.4.4a) and (3.1.17,18), this change of variables
makes sense globally.) Plugging these variables in Lagrangian
(3.4.3) gives
$$
L_{11}=(\frac{1}{2}(\partial_{\tau}x,\partial_{\tau}x)+
2(\nabla_{\partial_{\tau}x}\rho,\phi)+2(\nabla_{\partial_{\tau}x}\phi,\rho)+\cdots
d\tau\wedge )d\sigma
$$
where $\cdots$ stand for the terms not containing $\partial_{\tau}$.
Since
$$
(\nabla_{\partial_{\tau}x}\phi,\rho)=-(\phi,\nabla_{\partial_{\tau}x}\rho)+
\partial_{\tau}(\phi,\rho)=(\nabla_{\partial_{\tau}x}\rho,\phi,)+
\partial_{\tau}(\phi,\rho), \eqno{(3.4.24)}
$$
we obtain
$$
\tilde{L}_{11}=(\frac{1}{2}(\partial_{\tau}x,\partial_{\tau}x)+
4(\nabla_{\partial_{\tau}x}\rho,\phi)+\cdots )d\tau\wedge d\sigma=
L_{11}\; mod\; d_{\rho}(...) \eqno{(3.4.25)}
$$
 It is immediate to derive from
(3.4.25)
$$
\gamma_{L}'=(g_{ij}\partial_{\tau}x^{j}+4\Gamma_{i
\alpha}^{s}\rho^{\alpha}\phi^{j}g_{sj}(x))\delta x^{i}\wedge
d\sigma+ 4g_{ij}(x)\delta \rho^{i}\phi^{j}\wedge
d\sigma,\eqno{(3.4.26)}
$$
which is equal to $\gamma_{L}|_{d\tau=0}$ modulo $d_{\rho}$-exact
terms.

If we let
$$
\rho_{i}=\frac{1}{4}g_{ij}(x)\phi^{j},\eqno{(3.4.27)}
$$
then
$$
\gamma_{L}'=(g_{ij}\partial_{\tau}x^{j}+\Gamma_{i
\alpha}^{s}\rho^{\alpha}\rho_{s})\delta x^{i}\wedge d\sigma+ \delta
\rho^{i}\rho_{i}\wedge d\sigma,\eqno{(3.4.28)}
$$
The substitution
$$
x_{i}\mapsto g_{ij}\partial_{\tau}x^{j}+\Gamma_{i
\alpha}^{s}\rho^{\alpha}\rho_{s}, \phi^{i}\mapsto\rho^{i},
\phi_{i}\mapsto\rho_{i}\eqno{(3.4.29)}
$$
makes sense as a globally defined map
$$
d_{TM}\tilde{L}_{11}: T(\Pi TM)\rightarrow T^{*}(\Pi
TM).\eqno{(3.4.30)}
$$
It is a super-analogue of the Legendre transform, (1.6.10), which
was envisaged in sect.3.2.4; indeed, if $x_{i}=\partial_{x^{i}}$,
then the first of assignments (3.4.29) is exactly splitting
(3.2.26). The $D_{\Sigma/\Sigma''}$-manifold property allows to
extend this map unambiguously to the jet-spaces, and it is clear
that such map identifies $\delta \gamma'_{L}$ with $\omega$ from
(3.2.19).

\bigskip

{\bf 3.4.5.} Therefore, $\tilde{\CH}^{can}$ is to the
(1,1)-supersymmetric $\sigma$-model what $\CH^{can}$ is to the
ordinary $\sigma$-model. In particular, writing integrals of motion
(3.4.15)  in terms of the new variables introduced in sect. 3.4.4
provides a free field realization of $\CN1^{\pm}$. The result, which
we will discuss in the context of the K\"ahler geometry, see the
next section, is presumably the quasiclassical limit of the formulas
obtained in [B-ZHS].

\bigskip

{\bf 3.4.6.} {\it The K\"ahler case: (2,2)-supersymmetry and the
Witten Lie algebra.}

It is an exciting discovery going back to [Zu, A-GF] that in the
K\"ahler case the supersymmetry algebra becomes twice as large.

{\bf 3.4.6.1.} Let  then  $M$ be a complex manifold and $(.,.)$  a
K\"ahler metric on it. To handle this case, we will change the
notation somewhat: the natural vector bundles, such as $TM$, will be
assumed to be complexified, and decompositions, such as
$TM=T^{10}M\oplus T^{01}M$, will arise. What has been treated as a
vector field, e.g. $\partial_{\tau}x$, $\partial_{\tau}\psi_{+}$,
will become a section of $T^{10}M$, and $\partial_{\tau}\bar{x}$,
$\partial_{\tau}\bar{\psi}_{+}$ will stand for the complex conjugate
sections. We will also let, sloppily but customarily,
$$
\partial_{\tau} \overline {x^{j}}=\partial_{\tau} x^{\bar{j}},
\overline{\psi_{\pm}^{j}}=\psi_{\pm}^{\bar{j}}. \eqno{(3.4.31)}
$$
The defining property of the K\"ahler metric
$$
\nabla(T^{10})\subset T^{10},\; \nabla(T^{01})\subset
T^{01}\eqno{(3.4.32)}
$$
is crucial for what follows.

Computing as in 3.4.4 (and using (3.4.32)) one obtains
$$
L_{11}=(-(\partial_{u}x,\partial_{u}\bar{x})
-(\partial_{v}x,\partial_{u}\bar{x})
+(\nabla_{\partial_{v}}\psi_{+},\bar{\psi}_{+})
+(\nabla_{\partial_{u}}\psi_{-},\bar{\psi}_{-})+\cdots)du\wedge dv\;
mod\; d_{\rho},\eqno{(3.4.33)}
$$
where the terms not containing $\partial_{u}$, $\partial_{v}$ are
omitted.

Property (3.4.32) implies (and (3.4.33) supports) that w.r.t. the
grading on $\CO_{J^{\infty}(M_{\hat{\Sigma}})^{o}}$ defined by
$$
\psi_{\pm}\mapsto 1,\; \bar{\psi}_{\pm}\mapsto -1 \eqno{(3.4.34)}
$$
$L_{11}$ is homogeneous of degree 0. Therefore, any homogeneous
component of a symmetry of $L_{11}$ is also a symmetry. Integrals of
motion (3.4.15) afford decomposition
$$
Q^{+}_{f}= Q^{++}_{f}+Q^{+-}_{f}, Q^{-}_{f}=
Q^{-+}_{f}+Q^{--}_{f}\eqno{(3.4.35)}
$$
into the sum of degree $\pm 1$ components, which implies that the
entire quadruple
$$
\{Q^{++}_{f},Q^{+-}_{f},  Q^{-+}_{f},Q^{--}_{f}\}\subset
\tilde{\CI}_{L},\eqno{(3.4.36)}
$$
and this extends (3.4.17a) to an embedding of a pair of
N=2-superconformal Lie algebras
$$
\Vir^{\pm}\hookrightarrow \CN 1^{\pm}\hookrightarrow
\CN2^{\pm}\hookrightarrow\tilde{\CI}_{L}\hookrightarrow\tilde{\CH}^{can}.
\eqno{(3.4.37)}
$$

In particular, (and it follows from the consideration of the degree)
$$
[Q^{++}_{f},Q^{++}_{g}]=[Q^{--}_{f},Q^{--}_{g}]=0.\eqno{(3.4.38)}
$$
Witten has used these relations [W2,W3] to define what in the
present context becomes {\it Witten Lie algebra sheaves}:
$$
\CW^{\pm}\buildrel\text{def}\over= \frac{\{X\in\tilde{\CH}^{can}:\;
[Q^{\mp,\mp}_{1},X]=0\}} {\{[Q^{\mp,\mp}_{1},X]\text{ all
}X\in\tilde{\CH}^{can}\}}. \eqno{(3.4.39)}
$$
(There are, of course, two more versions of these sheaves.)

\bigskip

{\bf 3.4.6.2.} {\it Some formulas.} For the purpose of writing
embeddings such as (3.4.37) explicitly, rewrite (3.4.33) using
$\sigma$ and $\tau$ which were defined in (3.4.18)
$$
\aligned
L_{11}=&(g_{i\bar{j}}\partial_{\tau}x^{i}\partial_{\tau}x^{j}
-2\partial_{\tau}x^{i}\Gamma_{ij}^{s}g_{s\bar{t}}\psi_{+}^{j}\psi_{+}^{\bar{t}}
+2\partial_{\tau}x^{\bar{i}}\Gamma_{\bar{i}\bar{j}}^{\bar{s}}g_{\bar{s}\bar{t}}
\psi_{-}^{\bar{j}}\psi_{-}^{t})d\tau\wedge d\sigma\\
&+(2g_{i\bar{j}}\psi_{+}^{\bar{j}}\partial_{\tau}\psi_{+}^{i}
-2g_{\bar{i}j}\psi_{-}^{j}\partial_{\tau}\psi_{-}^{\bar{i}})d\tau\wedge
d\sigma\cdots
\endaligned
\eqno{(3.4.40)}
$$
where the terms not containing $\partial_{\tau}$ are omitted. It
follows that, cf. (3.4.26),
$$
\aligned \gamma_{L}'=&(g_{i\bar{j}}\partial_{\tau}x^{\bar{j}}
-2\Gamma_{ij}^{s}g_{s\bar{t}}\psi_{+}^{j}\psi_{+}^{\bar{t}})\delta
x^{i}\wedge d\sigma + (g_{i\bar{j}}\partial_{\tau}x^{i}
+2\Gamma_{\bar{j}\bar{i}}^{\bar{s}}g_{\bar{s}t}\psi_{-}^{\bar{i}}\psi_{-}^{t})\delta
x^{\bar{j}}\wedge d\sigma\\
&(-2g_{i\bar{j}}\psi_{+}^{\bar{j}}\delta\psi_{+}^{i}+2g_{\bar{i}j}\psi_{-}^{j}\delta\psi_{-}^{\bar{i}})
\wedge d\sigma.
\endaligned
\eqno{(3.4.41)}
$$
If we let
$$
\psi=\psi_{+},\bar{\psi}=\bar{\psi}_{-},\eqno{(3.4.42a)}
$$
and
$$
\psi_{i}=-2g_{i\bar{j}}\psi_{+}^{\bar{j}},
\psi_{\bar{i}}=2g_{\bar{i}j}\psi_{-}^{j},\eqno{(3.4.42b)}
$$
then
$$
\aligned \gamma_{L}'=&(g_{i\bar{j}}\partial_{\tau}x^{\bar{j}}
+\Gamma_{ij}^{s}\psi^{j}\psi_{s})\delta x^{i}\wedge d\sigma +
(g_{i\bar{j}}\partial_{\tau}x^{i}
+\Gamma_{\bar{j}\bar{i}}^{\bar{s}}\psi^{\bar{i}}\psi_{\bar{s}})\delta
x^{\bar{j}}\wedge d\sigma\\
&(\psi_{i}\delta\psi^{i}+\psi_{\bar{i}}\delta\psi^{\bar{i}}) \wedge
d\sigma.
\endaligned
\eqno{(3.4.43)}
$$
Therefore, the coordinate form of the super-Legendre transform
(3.4.30) is

$$
\aligned
 x_{i}&\mapsto g_{i\bar{j}}\partial_{\tau}x^{\bar{j}}
+\Gamma_{ij}^{s}\psi^{j}\psi_{s}\\
x^{\bullet}&\mapsto x^{\bullet}\\
 x_{\bar{j}}&\mapsto
g_{i\bar{j}}\partial_{\tau}x^{i}
+\Gamma_{\bar{j}\bar{i}}^{\bar{s}}\psi^{\bar{i}}\psi_{\bar{s}}\\
\phi^{\bullet}&\mapsto\psi^{\bullet},\;
\phi_{\bullet}\mapsto\psi_{\bullet}
\endaligned
\eqno{(3.4.44)}
$$
Plugging these in (3.4.17) and extracting homogeneous components as
in (3.4.35) one obtains, upon letting $d\tau=0$,
$$
\aligned
Q_{f}^{--}&=f(\sigma-\tau)(-x_{\bar{j}}\phi^{\bar{j}}+g_{i\bar{j}}\partial_{\sigma}x^{i}\phi^{\bar{j}})d\sigma,\\
Q_{f}^{-+}&=2f(\sigma-\tau)(\partial_{\sigma}x^{\bar{j}}\phi_{\bar{j}}-g^{\bar{j}i}x_{i}\phi_{\bar{j}}+
g^{\bar{j}i}\Gamma_{i\alpha}^{s}\phi^{\alpha}\phi_{s}\phi_{\bar{j}})d\sigma.
\endaligned
\eqno{(3.4.45^{-})}
$$
$$
\aligned
Q_{f}^{++}&=f(\sigma+\tau)(x_{j}\phi^{j}+g_{\bar{j}i}\partial_{\sigma}x^{\bar{j}}\phi^{i})d\sigma,\\
Q_{f}^{+-}&=-2f(\sigma+\tau)(\partial_{\sigma}x^{i}\phi_{i}+g^{i\bar{j}}x_{\bar{j}}\phi_{i}-
g^{i\bar{j}}\Gamma_{\bar{j}\bar{\alpha}}^{\bar{s}}\phi^{\bar{\alpha}}\phi_{\bar{s}}\phi_{i})d\sigma.
\endaligned
\eqno{(3.4.45^{+})}
$$

One may wish at this point to use these formulas to compute Witten's
Lie  algebra sheaf (3.4.39). Two things transpire immediately:
first, the role played by $f$ in all of this is rather superficial
and, second, if one removes from the first of $(3.4.45^{-})$ the
annoying $g_{i\bar{j}}\partial_{\sigma}x^{i}\phi^{j}$ ( and
$g_{i\bar{j}}\partial_{\sigma}x^{i}\phi^{j}$ from the first of
$(3.4.45^{+})$ resp.), then it becomes exactly the $\bar{\partial}$-
($\partial$- resp.) differential; and so, perhaps, $\CW^{\pm}$
should be of completely holomorphic (antiholomorphic resp.) nature.
This is all true, but the language suited to analysis of such issues
is that of vertex Poisson algebras.

\bigskip

{\bf 3.5. Vertex Poisson algebra interpretation. Witten's models}

\bigskip

The sheaf $\tilde{\CH}^{can}$ is the tip of an iceberg. It is, just
as its purely even counterpart $\CH^{can}$ was, sect. 1.5.4, a Lie
algebra sheaf attached to a certain sheaf of vertex Poisson
superalgebras

\bigskip

{\bf 3.5.1.} The notion of a super-SVDO is quite analogous to the
one we discussed in sect. 2. It is a $\BZ_{+}$-graded vertex Poisson
superalgebra $V=V_{0}\oplus V_{1}\oplus\cdots$ such that
$$
V_{0}=C^{\infty}(\Pi TU_{\Sigma}),\; U\subset\BR^{n},\eqno{(3.5.1)}
$$
and, non-canonically,
$$
V_{1}=\CT_{U_{\Sigma}/\Sigma}(U_{\Sigma})+\Omega_{U_{\Sigma}/\Sigma}(U_{\Sigma})+
\Pi\left(\CT_{U_{\Sigma}/\Sigma}(U_{\Sigma})+\Omega_{U_{\Sigma}/\Sigma}(U_{\Sigma})\right).
\eqno{(3.5.2)}
$$
Classification of such algebras [GMS3], under some obvious
non-degeneracy assumptions, is obtained in a way similar to sect.
2.4.3, 2.5. They form an $\Omega^{3,cl}(U)$-torsor, i.e., given a
super-SVDO $V$ and a closed 3-form $H\in \Omega^{3,cl}(U)$, an
operation
$$
(V, H)\mapsto V\plus H\eqno{(3.5.3)}
$$
is defined, where $V\plus H$ is a super-SVDO different from $V$ only
in that the operation
$$
_{(0)}:\CT_{U_{\Sigma}/\Sigma}\otimes
\CT_{U_{\Sigma}/\Sigma}\rightarrow V_{1}
$$
is replaced with
$$
_{(0)_{H}}=_{(0)}+H,\eqno{(3.5.4)}
$$
cf. (2.4.16-18). (This involves only even components of $V_{1}$.)
One has
$$
Mor(V,V\plus H)=\{\alpha\in\Omega^{2}(U)\text{ s.t.
}d\alpha=H\}.\eqno{(3.5.5)}
$$
In particular,
$$
\Omega^{2,cl}(U)\iso\text{Aut}(V), \eqno{(3.5.6a)}
$$
where the automorphism corresponding to $\alpha$ is the one
determined by the shear
$$
\CT_{U_{\Sigma}/\Sigma}(U_{\Sigma})\ni\xi\mapsto\xi+\iota_{\xi}\alpha,\eqno{(3.5.6b)}
$$
cf (2.4.19).

All of this can be spread over manifolds. There is a distinguished
such sheaf of super-SVDOs, the vertex Poisson de Rham complex [MVS],
$\Omega^{poiss}_{M}$. As an $\CO_{M}$-module,
$$
\Omega^{poiss}_{M}=\pi_{*}\CO_{J^{\infty}(T^{*}(\Pi
TM)_{\Sigma/\Sigma''})},\eqno{(3.5.7)}
$$
where $\pi$ is the projection $J^{\infty}(T^{*}(\Pi
TM)_{\Sigma/\Sigma''})\rightarrow M_{\Sigma}$,
 The operations are determined by the requirement that they all be
 of classical origin -- as in Proposition 2.7.1. Here are some
 examples written down in local coordinates, cf. (3.2.19):
 $$
 \aligned
 &(x_{i})_{(0)}f(x)=\partial_{x^{i}}f(x),
 (x_{\bar{i}})_{(0)}f(x)=\partial_{x^{\bar{i}}}f(x)\;
 (\phi_{i})_{(0)}\phi^{j}=\delta_{i}^{j},(\phi_{\bar{i}})_{(0)}\phi^{\bar{j}}=\delta_{\bar{i}}^{\bar{j}}\\
 &\xi_{(0)}\eta=[\xi,\eta],\;\xi_{(o)}\alpha=\text{Lie}_{\xi}\alpha,
 \endaligned
 \eqno{(3.5.8a)}
 $$
 where $\xi=f^{i}(x)x_{i}+f^{\bar{i}}(x)x_{\bar{i}}$,
$\eta=g^{i}(x)x_{i}+g^{\bar{i}}(x)x_{\bar{i}}$,
$\alpha=h_{i}(x)\partial_{\sigma}x^{i}+h_{\bar{i}}(x)\partial_{\sigma}x^{\bar{i}}$,
the vertex algebra derivation being
 $$
 T=\rho(\partial_{\sigma}).\eqno{(3.5.8b)}
 $$

 (The twist that takes care of functions of $\sigma$ and imposed
 in the even case in Definition 2.7.2 has been tacitly assumed
 throughout.)

 One has, analogously to Proposition 2.6.1,

 \bigskip

 {\bf 3.5.1.1. Proposition.} {\it a) The set of isomorphism classes
 of sheaves of super-SVDOs on $M$ is identified with
 $H^{3}(M,\BR)$.

 b)If $\CV$ is a sheaf of super-SVDOs, then
 $$
 \text{Aut}\CV\iso H^{0}(M,\Omega_{M}^{2,cl}).
 $$}

 \bigskip

 Let
 $$
 \text{Lie}(\CV)=\CV/T(\CV).\eqno{(3.5.9)}
 $$
 Operation $_{(0)}$ makes $\text{Lie}\CV$ into a sheaf of Lie
 superalgebras. One has, cf. Proposition 2.7.3,

 \bigskip

 {\bf 3.5.2. Proposition.}
 {\it The algebra sheaves $\tilde{\CH}^{can}$ and
 $\text{Lie}(\Omega^{poiss}_{M})$ are isomorphic.}

 \bigskip

 {\bf 3.5.3.} Some of the constructions
 above are simplified when performed in the framework of
 of vertex Poisson superalgebras because some of the Lie algebras
 considered are the  value of the Lie-functor. For example, there
 are N=1,2 supersymmetric vertex Poisson algebras [K], $N1$ and
 $N2$, such that the N=1,2 supersymmetric Lie superalgebras, which
 appeared in (3.4.12), are
 $$
\CN1=\text{Lie}(C^{\infty}(\Sigma)\otimes N1),\;
\CN2=\text{Lie}(C^{\infty}(\Sigma)\otimes N2).\eqno{(3.5.10)}
$$

 The elements, see (3.4.45),
$$
\aligned
Q^{--}&=-x_{\bar{j}}\phi^{\bar{j}}+g_{i\bar{j}}\partial_{\sigma}x^{i}\phi^{\bar{j}},\\
Q^{-+}&=2(\partial_{\sigma}x^{\bar{j}}\phi_{\bar{j}}-g^{\bar{j}i}x_{i}\phi_{\bar{j}}+
g^{\bar{j}i}\Gamma_{i\alpha}^{s}\phi^{\alpha}\phi_{s}\phi_{\bar{j}})
\endaligned
\eqno{(3.5.11^{-})}
$$
$$
\aligned
Q^{++}&=x_{j}\phi^{j}+g_{\bar{j}i}\partial_{\sigma}x^{\bar{j}}\phi^{i}\\
Q^{+-}&=-2(\partial_{\sigma}x^{i}\phi_{i}+g^{i\bar{j}}x_{\bar{j}}\phi_{i}-
g^{i\bar{j}}\Gamma_{\bar{j}\bar{\alpha}}^{\bar{s}}\phi^{\bar{\alpha}}\phi_{\bar{s}}\phi_{i}).
\endaligned
\eqno{(3.4.45^{+})}
$$

define global sections of $\Omega^{poiss}_{M}$. By definition, the
following analogue of (2.8.12) and (2.9.25) holds true.

\bigskip

{\bf 3.5.3.1. Lemma.} {\it The two pairs of global sections
$(Q^{--},Q^{-+})$ and  $(Q^{++},Q^{+-})$ generate, inside
$H^{0}(M,\Omega^{poiss}_{M})$, two pairwise Poisson-commuting copies
of the vertex Poisson N=2 superalgebra:
$$
N2^{+}\hookrightarrow H^{0}(M,\Omega^{poiss}_{M})\hookleftarrow
N2^{-},\; (N2^{+}_{(n)}(N2^{-})=0\text{ if }n\geq 0. \eqno{(3.5.12)}
$$}

\bigskip

 A streamlined version of Witten's Lie algebra sheaf (3.4.39) is
 Witten's vertex Poisson algebra sheaf defined as follows.
 Relations (3.4.38) in the vertex algebra context imply that each
element of the quadruple
$\{Q^{\bullet,\bullet}_{(0)},\bullet=\pm\})$ and various linear
combination thereof are differentials of the sheaf
$\Omega^{poiss}_{M}$. Letting $Q_{(0)}$ be one such differential, we
obtain a cohomology sheaf
$$
H_{Q}(\Omega^{poiss}_{M})\buildrel\text{def}\over =
\frac{\text{Ker}Q_{(0)}}{\text{Im}Q_{(0)}}.\eqno{(3.5.13)}
$$
It is  a vertex Poisson algebra sheaf -- a well-known fact and an
immediate consequence of $_{(0)}$ being a derivation of all
$_{(n)}$-products (super-analogue of Jacobi identity, sect. 2.1,
II.2).

Of sheaves (3.5.13) the following 3 will be of interest to us:

\bigskip

{\bf 3.5.3.2. Definition.} (cf. [W2])
$$
\text{A-model sheaf: }
W_{A}=H_{Q^{--}+Q^{++}}(\Omega^{poiss}_{M}),\eqno{(3.5.14a)}
$$
$$
\text{B-model sheaf: }
W_{B}=H_{Q^{--}+Q^{+-}}(\Omega^{poiss}_{M}),\eqno{(3.5.14b)}
$$
$$
\text{half-twisted model sheaf: }
W_{1/2}=H_{Q^{--}}(\Omega^{poiss}_{M}),\eqno{(3.5.14c)}
$$
\bigskip
The relation of (3.5.13-14) to (3.4.39) is that
$$
\CW^{-}=\text{Lie}(W_{1/2}),\eqno{(3.5.15)}
$$
to give but one example.

The cohomology, $H^{*}(M,\CV)$, of a sheaf of vertex Poisson
algebras $\CV$ is a vertex Poisson algebra, of course. We are led
then, following [W2], to

\bigskip
{\bf 3.5.3.3. Definition.}
$$
\text{A-model vertex Poisson algebra: }
H^{*}(M,W_{A}),\eqno{(3.5.16a)}
$$
$$
\text{B-model  vertex Poisson algebra: }
H^{*}(M,W_{B}),\eqno{(3.5.16b)}
$$
$$
\text{half-twisted model  vertex Poisson algebra: }
H^{*}(M,W_{1/2}),\eqno{(3.5.16c)}
$$
\bigskip

\bigskip

{\bf 3.5.4. Theorem.} {\it Let $M$ be K\"ahlerian. Then

1) the following isomorphisms are valid:}
$$
H^{*}(M,W_{A})\iso H^{*}(M,\BC),\eqno{(3.5.17a)}
$$
$$
H^{*}(M,W_{B})\iso H^{*}(M,\Lambda^{*}\CT_{M}),\eqno{(3.5.17b)}
$$
$$
H^{*}(M,W_{1/2})\iso H^{*}(M,\Omega^{poiss,an}_{M}),\eqno{(3.5.17c)}
$$
{\it where $\Omega^{poiss,an}_{M}$ is a purely holomorphic version
of the sheaf $\Omega^{poiss}_{M}$ [MSV];

2) embedding $N2^{+}\hookrightarrow\Omega^{poiss}_{M}$, (3.5.12),
descends to an embedding
$N2^{+}\hookrightarrow\Omega^{poiss,an}_{M}$ whose image coincides
with N=2 superconformal structure introduced on
$\Omega^{poiss,an}_{M}$.}

\bigskip

{\bf 3.5.4.1.}{\it Remark.} 1) Of these, the first two are finite
dimensional supercommutative algebras and as such are trivial
examples of a vertex Poisson algebra with zero derivation $T$ as
noted in sect.2.3. Contrary to this, the last one is a full-fledged
infinite dimensional vertex Poisson algebra. Being infinite
dimensional it is characterized by its character ($q$-dimension),
which is closely related to the elliptic genus of $M$. The algebra
can be quantized, and the character of the quantum version has
provided some  insights into the elliptic genus [BL,MS,GM1, GM2].

2) This theorem, especially (3.5.17c) is a refined version of [Kap].
In fact, Kapustin deals with the quantum version of this result; we
will discuss quantization in the next section.

\bigskip

{\bf 3.5.4.2.} {\it Sketch of proof.}

Apply to $\Omega^{poiss}_{M}$ automorphism (3.5.6a-b) determined by
the K\"ahler 2-form $g_{i\bar{j}}dx^{i}\wedge dx^{\bar{j}}$. As a
result, $(Q^{--})_{(0)}$ will be replaced with a vertex analogue of
the $\bar{\partial}$-differential:
$$
\bar{\partial}_{vert}=(x_{\bar{j}}\phi^{\bar{j}})_{(0)}\eqno{(3.5.18)}.
$$
Essentially by definition,
$$
(\Omega^{poiss,an}_{M},0)\hookrightarrow(\Omega^{poiss}_{M},\bar{\partial}_{vert}),
$$
is a quasiisomorphism [MSV]. Indeed, a glance at (3.5.8a) convinces
that $x^{\bar{j}}$ $\phi^{\bar{j}}$ are not
$\bar{\partial}_{vert}$-cocycles, and $x_{\bar{j}}$ $\phi_{\bar{j}}$
are $\bar{\partial}_{vert}$-cohomologous to 0. Therefore
$\bar{\partial}_{vert}$ effectively kills all antiholomorphic
variables, leaving holomorphic ones intact. This defines a purely
holomorphic analogue of $\Omega^{poiss}_{M}$, that is,
$\Omega^{poiss, an}_{M}$.

Hence a quasiisomorphism
$$
(\Omega^{poiss,an}_{M},0)\hookrightarrow
(\Omega^{poiss}_{M},(Q^{--})_{(0)}),
$$
which proves (3.5.17c).

In (3.5.17a-b) one more  differential is turned on. Definition
(3.5.11) implies that  upon the  same shear by the K\"ahler form,
$$
Q^{++}=x_{j}\phi^{j}. \eqno{(3.5.19)}
$$
Therefore, $(Q^{--})_{(0)}+(Q^{++})_{(0)}$ is a vertex analogue of
total de Rham differential, and (3.5.17b) becomes essentially [MSV],
Theorem 2.4.

Similarly,  in the $(Q^{--})_{(0)}$-cohomology,
$$
Q^{+-}=-4\partial_{\sigma}x^{i}\phi_{i}, \eqno{(3.5.20)}
$$
and a simple analysis along the lines of [MSV], sect. 2.3-2.4, shows
that
$$
H_{\partial_{\sigma}x^{j}\phi_{j}}(\Omega^{poiss,an}_{M})\iso
H^{*}(M,\Lambda^{*}\CT_{M}),
$$
as desired.

Item 2) is a result of checking (3.5.19,20) against [MSV], (2.3b).

\bigskip

Next, we establish concrete complexes which compute vertex Poisson
algebras of the A-, B-, and half-twisted models.

\bigskip

{\bf 3.5.5. Corollary.}
$$
H^{*}(M,W_{A})\iso
H_{Q^{--}+Q^{++}}(\Gamma(M,\Omega^{poiss}_{M})),\eqno{(3.5.21a)}
$$
$$
H^{*}(M,W_{B})\iso
H_{Q^{--}+Q^{+-}}(\Gamma(M,\Omega^{poiss}_{M})),\eqno{(3.5.21b)}
$$
$$
H^{*}(M,W_{1/2})\iso
H_{Q^{--}}(\Gamma(M,\Omega^{poiss}_{M})),\eqno{(3.5.21c)}
$$

\bigskip

The sheaf $\Omega^{poiss}_{M}$ is a complex w.r.t.  the 3
differentials used above. Hence there arise 3 different
hypercohomology groups, $\BH_{A}(\Omega^{poiss}_{M})$,
$\BH_{B}(\Omega^{poiss}_{M})$, $\BH_{1/2}(\Omega^{poiss}_{M})$. Each
can be computed by any of the two spectral sequences. The
computation using one of them is the content of Theorem 3.5.4. It
says that the result is the vertex Poisson algebra of A-, B-, and
half-twisted models resp. The computation using another will then
prove the corollary, because the sheaf $\Omega^{poiss}_{M}$ being
flabby,  $H^{j}(M,\Omega^{poiss}_{M})=0$ if $j>0$.

{\it Remark.} In view of Theorem 3.5.4,  isomorphisms (3.5.21a,b)
are vertex Poisson algebra versions of the de Rham complex and
$\bar{\partial}$-resolution of the algebra of polyvector fields
resp., while (3.5.21c) is the $\bar{\partial}$-resolution of the
vertex Poisson de Rham complex.

\bigskip

{\bf 3.5.6.} {\it H-flux.}

Let us place ourselves in the situation of  sect. 2.8, (2.8.4).
Consider the $H$-twist of Lagrangian (3.4.2):
$$
L^{H}=\{L+\beta^{j}(D_{+}X,D_{-}X)du\wedge
dv\wedge[d\theta^{+}d\theta^{-}]\}.\eqno{(3.5.22)}
$$
analogous  to (2.8.5). The argument completely parallel to that
leading to (2.8.7) finds the following vertex algebra analogue of
this operation.

\bigskip

{\bf 3.5.6.1. Lemma.}
$$
\tilde{\CH}^{\omega_{L^{H}}}_{Sol_{L^{H}}^{o}}\iso\text{Lie}(\Omega^{poiss}_{M}\plus
H).
$$
{\it where $\Omega^{poiss}_{M}\plus H$ is defined as in (3.5.3).}

\bigskip

Therefore, all the constructions, originating in [GHR] and further
explored in papers such as [BLPZ,KL], translate into different
vertex Poisson subalgebras of $\Omega^{poiss}_{M}\plus H$ depending
on the choice of a generalized K\"ahler structure.

\bigskip\bigskip

{\bf 3.6. Quantization. B-model moduli.}

\bigskip

This section is an announcement. It will be assumed throughout that
the automorphism by the K\"ahler form has been performed so that
$$
Q^{--}=\phi^{\bar{j}}x_{\bar{j}},\;
Q^{++}=\phi^{i}x_{i},\eqno{(3.6.0)}
$$
cf. sect.3.5.4.2.
\bigskip

{\bf 3.6.1.} The differential graded sheaves of vertex Poisson
algebras, $(\Gamma(M,\Omega^{poiss}_{M}), Q_{(0)})$, where $Q$ is
any of the differentials appearing in (3.5.21a,b,c), can be
quantized. What we mean by this is that, first, there is a sheaf of
vertex algebras $\Omega^{vert}_{M}$ [MSV] whose quasiclassical limit
is $\Omega^{poiss}_{M}$ and, second, this sheaf carries quantum
analogues of each of the 3 differentials. In fact, quantum versions
of $(Q^{--}_{(0)})$ and $(Q^{++}_{(0)})$ are in [MSV], and
$(Q^{+-})_{(0)}$ has been recently proposed in [B-ZHS]; in what
follows the use of the latter is easy to avoid. Thus there arise  3
vertex algebra versions of A-, B-, and half-twisted models resp.:
$$
H^{*}(M,W_{A}^{quant})\iso
H_{Q^{--}+Q^{++}}(\Gamma(M,\Omega^{vert}_{M})),\eqno{(3.6.1a)}
$$
$$
H^{*}(M,W_{B}^{quant})\iso
H_{Q^{--}+Q^{+-}}(\Gamma(M,\Omega^{vert}_{M})),\eqno{(3.6.1b)}
$$
$$
H^{*}(M,W_{1/2}^{quant})\iso
H_{Q^{--}}(\Gamma(M,\Omega^{vert}_{M})),\eqno{(3.6.1c)}
$$

 The first two, $H^{*}(M,W_{A}^{quant})$ and $H^{*}(M,W_{B}^{quant})$, coincide
 with their quasiclassical limits (3.5.21a,b). The  3rd is quite different and
is the cohomology of the  chiral de Rham complex,
$H^{*}(M,\Omega^{ch}_{M})$ [MSV].

Relation of this naive quantization to the genuine quantum string
theory is  expressed by saying, in physics language, that the latter
 equals the former ``perturbatively'', [Kap]. But let us show that
both (3.6.1b,c) can be further deformed along the
Barannikov-Kontsevich moduli space [BK]. We will focus on the
half-twisted model (3.6.1c).

\bigskip

{\bf 3.6.2.} Recall that associated to any differential  Lie
superalgebra $(\fg,d)$ there is a deformation functor,
$\text{Def}_{\fg}$, with domain the category of Artin algebras and
range the category of sets [Kon,BK]. In order to define it,
introduce the space of solutions to the Maurer-Cartan equation with
values in an Artin algebra A:
$$
MC_{\fg}(A)=\{\gamma: d\gamma+\frac{1}{2}[\gamma,\gamma]=0,\;
\gamma\in(\fg\otimes A)^{1}\}. \eqno{(3.6.2)}
$$
The operation
$$
(\fg\otimes A)^{1}\ni\gamma\mapsto d\beta +[\gamma,\beta]\text{ if }
\beta\in(\fg\otimes A)^{0}\eqno{(3.6.3)}
$$
does not preserve the set $MC_{\fg}(A)$, but it does so
infinitesimally, see a lucid explanation in [M2], Ch.2, sect.9.
Exponentiating (3.6.3) gives a group action
$$
G(A)^{0}\times MC_{\fg}(A)\rightarrow MC_{\fg}(A).\eqno{(3.6.4)}
$$
Define
$$
\text{Def}_{\fg}(A)=MC_{\fg}(A)/G(A)^{0}. \eqno{(3.6.5)}
$$
The motivation behind this ([M2], Ch.2, sect.9) is that

(i) if $\gamma$ is a solution of the Maurer-Cartan equation, then
$d+[\gamma,.]$ is also a differential, and

(ii) the adjoint action of $\fg^{0}$ results in the action on
solutions of the Maurer Cartan equation defined in (3.6.3).

Barannikov and Kontsevich apply this functor in the case where
$$
\aligned
 &\fg_{BK}=\Gamma(M,\Omega^{0,*}_{M}\otimes
T^{*,0}_{M}),\;d=\bar{\partial}\\
& [.,.]\text{ is the Schouten-Nijenhuis bracket.}
\endaligned
\eqno{(3.6.6)}
$$
Our task is similar but somewhat different. We need, see (3.6.1c),
to deform $(Q^{--})_{(0)}$ within the class of differentials on the
vertex algebra $\Gamma(M,\Omega^{vert}_{M})$. Even though the latter
is not a Lie algebra,  this deformation problem is governed by the
differential Lie superalgebra
$$
(\hat{\fg},d,[.,.])\buildrel\text{def}\over
=\left(\Gamma(M,\text{Lie}(\Omega^{vert}_{M})),(Q^{--})_{(0)},_{(0)}\right),\eqno{(3.6.7a)}
$$
where
$$
\Gamma(M,\text{Lie}(\Omega^{vert}_{M}))=\Gamma(M,\Omega^{vert}_{M}/T(\Omega^{vert}_{M})).
\eqno{(3.6.7b)}
$$
Indeed,
$$
_{(0)}:\hat{\fg}\otimes \Gamma(M,\Omega^{vert}_{M})\rightarrow
\Gamma(M,\Omega^{vert}_{M})\eqno{(3.6.8)}
$$
makes $\Gamma(M,\Omega^{vert}_{M})$  a $\hat{\fg}$-module, on which
$\hat{\fg}$ operates by derivations. Furthermore,
$$
((Q^{--})_{(0)}+\gamma_{(0)})^{2}=(Q^{--}_{(0)}\gamma)_{(0)}+\frac{1}{2}(\gamma_{(0)}\gamma)_{(0)}.
\eqno{(3.6.9)}
$$
Hence, if $\gamma$ satisfies the Maurer-Cartan equation, then
$(Q^{--})_{(0)}+\gamma_{(0)}$ is a differential. Let us define then
$$
\text{Def}_{\Gamma(M,\Omega^{vert}_{M})}=\text{Def}_{\hat{\fg}}.\eqno{(3.6.10)}
$$

By definition,
$$
\fg_{BK}\subset\Gamma(M,\Omega^{vert}_{M}),\eqno{(3.6.11)}
$$
which, by virtue of (3.6.7b), gives a map,  an injection, in fact,
$$
\iota:\fg_{BK}\hookrightarrow\hat{\fg}.\eqno{(3.6.12)}
$$
It is not a differential Lie algebra homomorphism, but its twisted
version
$$
\iota_{Q^{++}}:\fg_{BK}\rightarrow\hat{\fg},\; a\mapsto
Q^{++}_{(0)}\iota(a)\eqno{(3.6.13)}
$$
is; here $Q^{++}$ is a vertex analogue of the
$\partial$-differential; it has appeared in (3.6.1) and is defined
by the same formula as its quasiclassical limit (3.6.0). Indeed, it
is a pleasing exercise to check that the Schouten-Nijenhuis bracket
can be expressed in purely vertex algebra terms
$$
[\iota(a),\iota(b)]=\iota(a)_{(0)}(Q^{++}_{(0)}\iota(b)).\eqno{(3.6.14)}
$$
Note that morphism (3.6.13) changes the parity, as it should,
because $\fg_{BK}$ is an odd Lie superalgebra.

This proves

\bigskip

{\bf 3.6.2.1. Lemma.} {\it Map (3.6.13) defines a morphism of}
functors
$$
\text{Def}_{\fg_{BK}}\rightarrow\text{Def}_{\Gamma(M,\Omega^{vert}_{M})}.
\eqno{(3.6.15)}
$$
\bigskip

 If $M$ is a Calabi-Yau manifold, then  $\text{Def}_{\fg_{BK}}$
 is represented by a formal scheme that is the formal neighborhood
 of 0 of the superspace $H^{*}(M,\Lambda^{*}\CT_{M})$ [BK].  In
 particular, there exists a generic formal solution of the
 Maurer-Cartan equation in variables chosen to be any basis of the
 dual space $(H^{*}(M,\Lambda^{*}\CT_{M}))^{*}$. Therefore,

 \bigskip

 {\bf 3.6.2.2. Corollary.} {\it If $M$ is a Calabi-Yau manifold, then there is
 a formal family of vertex algebras
 $$
H^{*}(M,W_{1/2}^{quant})_{t}\iso
H_{Q^{--}_{t}}(\Gamma(M,\Omega^{vert}_{M})),\eqno{(3.6.16)}
$$
with base the formal neighborhood of 0 of the superspace
$H^{*}(M,\Lambda^{*}\CT_{M})$.}

\bigskip

Some of these deformations are not so formal; for example,
$(Q^{--})_{(0)}$ itself  depends quite explicitly on the choice of a
complex structure, see (3.6.0); this can be extended by including
generalized complex structures [G]; and considerable work has been
done in order to interpret other points of the Barannikov-Kontsevich
moduli space.

\bigskip

{\bf 3.6.3.} {\it Vertex Frobenius manifolds?}

It appears that there is more than just that to this story. The
events unfolding in the conformal weight zero component of
$H_{Q^{--}_{t}}(\Gamma(M,\Omega^{vert}_{M}))$ is precisely the
Barannikov-Kontsevich construction of the Frobenius manifold
structure on $\text{Def}_{\fg_{BK}}$. Furthermore, it is plausible
that each line of [BK] has a vertex algebra analogue valid up to
homotopy. For example, operation $_{(-1)}$ makes each vertex algebra
into a homotopy associative commutative algebra [LZ].  Furthermore,
the order 2 differential operator $\Delta$ defined on $\fg$, which
is essential for [BK], has a vertex analogue; this analogue is
$(Q^{++})_{(1)}$, which is well defined precisely when $M$ is a
Calabi-Yau manifold [MSV]. It is also an order 2 differential
operator of sorts in that
$$
[(Q^{++})_{(1)},
a_{(-1)}]-(Q^{++}_{(1)}a)_{(-1)}=(Q^{++}_{(0)}a)_{(0)},\eqno{(3.6.17)}
$$
which is a derivation of all $_{(n)}$-multiplications -- a remark of
Lian and Zuckerman, [LZ], Lemma 2.1.

What all of this seems to indicate is that there is a reasonable
definition of a a vertex Frobenius manifold of which
$H_{Q^{--}_{t}}(\Gamma(M,\Omega^{vert}_{M}))$ is an important
example.

\bigskip

\bigskip

\bigskip

\centerline{{\bf References}}

\bigskip

[AG] Arkhipov, S.; Gaitsgory, D. Differential operators on the loop
group via chiral algebras, {\it  Int. Math. Res. Not.}  2002, no. 4,
165--210

[A-GF] Alvarez-Gaumé, L., Freedman, D.Z., Geometrical structure and
ultraviolet finiteness in the supersymmetric $\sigma $-model. {\it
Comm. Math. Phys.} {\bf  80} (1981), no. 3, 443--451

[BK] Barannikov, S., Kontsevich, M., Frobenius manifolds and
formality of Lie algebras of polyvector fields,{\it Internat. Math.
Res. Notices}, 1998, no. 4, 201--215

[BD] Beilinson, A., Drinfeld, V., Chiral algebras,{\it American
Mathematical Society Colloquium Publications}, {\bf 51}, American
Mathematical Society, Providence, RI, 2004

[B-ZHS] Ben-Zvi D., Heluani R., Szczesny, M, Supersymmetry of the
chiral de Rham complex, posted on the archive: QA/0601532

[BL] Borisov L., Libgober A., Elliptic genera of toric varieties and
applications to mirror symmetry, {\it Inv. Math.} {\bf 140} (2000)
no. 2, 453-485,

[BLPZ] Bredthauer, A., Lindstr\"om, U., Persson, J., Zabzine, M.,
Generalized K\"ahler geometry from supersymmetric sigma models,
posted on the archive: hep-th/0603130

[Bre] Bressler P., The first Pontryagin class, posted on the
archive: AT/0509563

[DF] Deligne, P., Freed, D., Classical field theory, in: Quantum
fields and strings: A course for mathematicians v.1, P.Deligne et
al, editors, {\it AMS}, 2000

[DM] Deligne, P., Morgan, J., Notes on supersymmetry (following
J.Bernstein), in: Quantum fields and strings: A course for
mathematicians v.1, P.Deligne et al, editors, {\it AMS}, 2000

[Di] Dickey, L. A. Soliton equations and Hamiltonian systems, Second
edition. {\it Advanced Series in Mathematical Physics}, {\bf 26},
World Scientific Publishing Co., Inc., River Edge, NJ, 2003

[Dor] Dorfman, I. Ya. Dirac structures of integrable evolution
equations {\it Phys. Lett.} A {\bf 125} (1987), no. 5, 240--246

[FP] Feigin, B. Parkhomenko, S. Regular representation of affine
Kac-Moody algebras; in: Algebraic and geometric methods in
mathematical physics (Kaciveli, 1993), 415--424, Math. Phys. Stud.,
19, Kluwer Acad. Publ., Dordrecht, 1996.

[F] Frenkel, E., private communication

[FL] Frenkel, E., Losev A., Mirror symmetry in two steps: A-I-B,
posted on the archive: hep-th/0505131

[FB-Z] E.Frenkel, D.Ben-Zvi, Vertex algebras and algebraic curves,
{\it Mathematical Surveys and Monographs} {\bf 88}, 2001,

[FS] Frenkel, I., Styrkas K., Modified regular representations of
affine and Virasoro algebras, VOA structure and semi-infinite
cohomology, posted on the archive: QA/0409117

[GHR] Gates, S. J., Jr.; Hull, C. M.; Ro\v cek, M. Twisted
multiplets and new supersymmetric nonlinear $\sigma$-models, {\it
Nuclear Phys. B} {\bf 248} (1984), no. 1, 157--186

[GW] Gepner, D., Witten, E., String theory on group manifolds. {\it
Nuclear Phys.} B 278 (1986), no. 3, 493--549

[GM1] Gorbounov, V., Malikov, F., Vertex algebras and the
Landau-Ginzburg/Calabi-Yau correspondence, {\it Moscow Math. J.},
{\bf 4} (2004), no.3

[GM2] Gorbounov, V., Malikov F., The chiral de Rham complex and the
positivity of the equivariant signature of the loop space, posted on
the archive: AT/0205132

[GMS1] V.Gorbounov, F.Malikov, V.Schechtman, Gerbes of chiral
differential operators. II. Vertex algebroids,  {\it Inv. Math.},
{\bf 155}, 605-680 (2004)

[GMS2] Gorbounov, V.; Malikov, F.; Schechtman, V., On chiral
differential operators over homogeneous spaces. {\it Int. J. Math.
Math. Sci.} {\bf 26} (2001), no. 2, 83--106

[GMS3] V.Gorbounov, F.Malikov, V.Schechtman, Gerbes of chiral
differential operators. III, posted on the archive: AG/0005201

[G] Gualtieri, M., Generalized complex geometry, posted on the
archive: DG/0401221

[HK] Heluani, R., Kac, V.G.,  Supersymmetric vertex algebras, posted
on the archive: QA/0603633

[K] V.Kac, Vertex algebras for beginners, 2nd edition, {\it AMS}, 1998

[Kap] Kapustin A., Chiral de Rham complex and the half-twisted
sigma-model, posted on the archive: hep-th/0504074

[KL] Kapustin, A., Li, Yi, Topological sigma-models with H-flux and
twisted generalized complex manifolds, posted on the archive:
hep-th/0407249

[Kon] Kontsevich, M., Deformation quantization of Poisson manifolds.
{\it Lett. Math. Phys}. {\bf 66} (2003), no. 3, 157--216

[L] Leites, D., Introduction to the theory of supermanifolds, {\it
Russ.Math.Surveys}, {\bf 35} (1980), no.1, 1-64

[LZ] Lian, B. H., Zuckerman, G. J., New perspectives on the
BRST-algebraic structure of string theory, {\it Comm. Math. Phys}.
{\bf 154} (1993), no. 3, 613--646

[LWX] Liu, Z.-J., Weinstein, A., Xu, P., Manin triples for Lie
bialgebroids, {\it J.Diff.Geom.} {\bf 45} (1997), 547-574

[M1] Manin Yu.I., Gauge field theory and complex geometry,
Grundlehren {\bf 289} (1988), Springer-Verlag

[M2] Manin, Yu.I., Frobenius manifolds, quantum cohomology, and
moduli spaces, {\it Colloquium Publications} {\bf 47}, American
Mathematical Society, Providence, Rhode Island, 1999

[MS] F.Malikov, V.Schechtman, Deformations of vertex algebras, quantum
cohomology of toric varieties, and
elliptic genus, {\it Comm. Math. Phys.}
{\bf 234} (2003), no.1 77-100,

[MSV] F. ~Malikov, V. ~Schechtman, A. ~Vaintrob, Chiral de Rham complex,
{\it Comm. Math. Phys.}, {\bf 204} (1999), 439 - 473,

[QFS] Quantum fields and strings: A course for mathematcians, v.1,
2, P.Deligne et al, editors, {\it AMS}, 2000

[T] Takens, F., A global version of the inverse problem of the
calculus of variations, {\it J. Differential Geom.}{\bf  14} (1979),
no. 4, 543--562

[V] Vinogradov, A. M. Cohomological analysis of partial differential
equations and secondary calculus. {\it Translations of Mathematical
Monographs}, 204. American Mathematical Society, Providence, RI,
2001

[W1] Witten, E., Nonabelian bosonization in two dimensions,{\it
Comm. Math. Phys.}{\bf 92} (1984), no. 4, 455--472

[W2] Witten, E., Mirror manifolds and topological field theories,
in: Essays on mirror symmetry, S.T.Yau ed., International Press,
Hong Kong, 1992

[W3] E. Witten, On the Landau-Ginzburg description of N=2 minimal
models, {\it Int.J.Mod.Phys.} {\bf A9} (1994), 4783-4800,
hep-th/9304026

[W4] Witten, E., Two-Dimensional Models With (0,2) Supersymmetry:
Perturbative Aspects, posted on the archive: hep-th/0504078

[Z] Zuckerman, Gregg J. Action principles and global geometry, in:
Mathematical aspects of string theory (San Diego, Calif., 1986),
259--284, Adv. Ser. Math. Phys., 1, World Sci. Publishing,

[Zu] Zumino, B., Supersymmetry and K\"ahler manifolds, {\it Phys.
Lett.}, {\bf 27B}(1979), 203

\bigskip\bigskip

F.M.: Department of Mathematics, University of Southern California,
Los Angeles, CA 90089, USA;\ fmalikov@math.usc.edu

\end{document}